\documentclass[10pt,reqno]{amsart}
\usepackage{import}
\usepackage{./preamble/Pakete}
\usepackage{./preamble/Befehle}

\begin{document}


\title[Standard extension algebras I]{Standard Extension Algebras I:\\
Perverse sheaves and Fukaya Calculus}

\author{Jens Niklas Eberhardt}
\email{mail@jenseberhardt.com}
\author{Catharina Stroppel}
\email{stroppel@math.uni-bonn.de}

\begin{abstract}
    In this first of a series of articles on standard extension algebras we study standard perverse sheaves on varieties with $\Gm$-actions. Based on Braden's hyperbolic localisation, we describe their extension algebra geometrically via a convolution structure on the intersections of attracting with repelling cells. We introduce a multiplicative structure on open Richardson varieties which provides a practical way to compose these extensions in case of flag varieties. For open Richardson varieties in Grassmannians we construct two explicit cell decompositions, of \emph{Gauss-} and of \emph{Deodhar-type}. It is shown that the latter is a stratification with the same combinatorics as Deodhar's decomposition. We introduce a calculus of Fukaya diagrams to encode the geometry of the decompositions. It provides a model for the cohomology of open Richardson varieties and thus for standard extensions. The calculus is motivated by the Mak--Smith Fukaya--Seidel category of a natural Lefschetz ﬁbration and 
 should allow to compute morphism spaces in there. We finally discuss the relation of our work to extensions of (parabolic) Verma modules in category $\cO$ as well as to the computation of $R$- and $R'$-polynomials. 
\end{abstract}
\maketitle
\section{Introduction} 
 
\subsection{Context} Studying varieties equipped with $\Gm$-actions and the geometry of associated stratifications is a recurring theme in geometric representation theory --- with Springer theory \cite{springerConstructionRepresentationsWeyl1978}, the (Kazhdan--)Lusztig conjectures \cite{carterSurveyWorkGeorge2006} and the geometric Satake equivalence \cite{mirkovicGeometricLanglandsDuality2007} as prominent examples. Typically, the geometry is studied via intersection cohomology of projective, of parity or of simple perverse sheaves. Our focus is however on standard perverse sheaves and their extensions with the main object of study: the \emph{standard extension algebra}.

The standard extension algebra encodes the gluing data between the categories of sheaves on the strata; extensions arise geometrically from the cohomology of normal cones `between' the strata.
We focus on Białynicki-Birula stratifications. There, one can replace normal cones with the intersections of attracting and repelling cells (these are open Richardson varieties \cite{richardsonIntersectionsDoubleCosets1992} in the case of flag varieties). 
The gluing data consists of the ($A_\infty$-) algebra structure on the extensions. It raises the natural questions of how to understand composition geometrically. We give an answer using the Drinfeld--Gaitsgory interpolation space \cite{drinfeldTheoremBraden2014}, \cite{drinfeldAlgebraicSpacesAction2015}, and concretize it in case of flag varieties via a multiplicative structure on open Richardson varieties.

For flag varieties, or equivalently category $\cO,$ computing the graded dimension of the standard extension algebra (encoded in the so-called $R'$-polynomials) is already a difficult problem. 
The famous \emph{Gabber--Joseph conjecture}, based on \cite{GabberJoseph}, predicts that $R'$-polynomials agree with Kazhdan--Lusztig's $R$-polynomials \cite{kazhdanRepresentationsCoxeterGroups1979}. The latter are well-understood and have an explicit geometric interpretation as point counting polynomials of open Richardson varieties.
However, the conjecture turned out to be incorrect, see \cite{Boe92}, so $R'$-polynomials remain mysterious and only known in special cases, see \cite{Abe} and \cite{koGradedExtensionsVerma2022}, \cite{dhillonExtensionsVermaModules2017} for some recent results. 
For Grassmannians, we provide both, a geometric and a homological interpretation for the $R'$-polynomials connecting work of Shelton \cite{sheltonExtensionsGeneralizedVerma1988} and Deodhar~\cite{deodharGeometricAspectsBruhat1985}.

The standard extension algebra is in general not formal; it carries non-trivial higher structure, see \cite{klamtExtAlgebrasParabolic2012}, \cite{makFukayaSeidelCategories2021}. Ideally one likes to have a concrete $dg$- or $A_\infty$-model encoding all this higher structure. Note that a complete knowledge of the higher data allows to recover the derived category and its highest weight structure, \cite{Keller}. With the example of Grassmannians, our article takes a first step. There are other approaches to package the higher structure via exact Borel subalgebras and boxes, \cite{KKO}, \cite{KM}. Apart from extreme cases, all examples coming from flag varieties or even Grassmannians are currently out of reach for these approaches, since they rely on a detailed knowledge of the involved algebras and on concrete data from the $A_\infty$-structure. This data is hoewever voluminous or simply unknown, see e.g. \cite{Strquiver}, \cite{klamtExtAlgebrasParabolic2012} for a taste of the complexity. 

The standard extension algebra is also relevant to categorification, where usually standard objects correspond to standard bases. Standard bases behave well under tensor products. As a consequence, a good understanding of the standard extension algebra is directly related to the notoriously difficult task of {\it taking} tensor products of categorifications, see \cite{manionHigherRepresentationsCornered2020} and for an overview and related questions  \cite{ICM}. Our explicit results for Grassmannians are related to tensor products of categorifications of quantum $\mathfrak{sl}_2$ representations underlying Khovanov homology.

Connections between categories of constructible sheaves and symplectic Fukaya--Seidel categories have been the focus of recent attention including the construction of $t$-structures or highest weight structures on the symplectic side,
see e.g. \cite{cotePerverseMicrosheaves2022}, \cite{makFukayaSeidelCategories2021} and references therein. Mak and Smith \cite{makFukayaSeidelCategories2021} give a Fukaya--Seidel category realization of the category of projective perverse sheaves on Grassmannians by establishing a symplectic realization of the extended arc algebra from \cite{StrSpringer}, \cite{BS3}. The $A_\infty$-structure of the symplectic realization of the standard extension algebra, namely the endomorphism algebra of Lefschetz thimbles, is unknown (in contrast to the arc algebra, 
\cite{abouzaidSymplecticArcAlgebra2016}, \cite{Barmeier}). Motivated by this, we introduce certain decorated Fukaya diagrams to give a $dg$-model for standard extensions. This Fukaya diagram calculus should allow explicit Floer homology calculations.
\vspace{-1mm}
\section*{Main results}
Let $Z$ be a smooth projective variety with a $\Gm$-action with discrete fixed point set $Z^0\subset Z.$ Then $Z$ admits a  stratification $Z^+=\bigsqcup_{w\in Z^0}Z^+_w\to Z$ into the attracting sets of the fixed points, see \cite{bialynicki-birulaTheoremsActionsAlgebraic1973}. Denote by $\D(Z)$ the derived category of constructible sheaves which are constant along the strata. The {\emph{(full) standard object} $\Delta\in \D(Z)$ is the direct sum of the constant sheaves on the strata. Our main actor is the \emph{standard extension algebra} $$\op{E}(\Delta)=\Hom_{\D(Z)}(\Delta,\Delta[\bullet]).$$
This algebra, equipped with its $A_\infty$-structure, completely determines the category $\D_{Z^+}(Z)\subset \D(Z)$ of sheaves constructible with respect to $Z^+$, see \cite{Keller}.
\subsection{Hyperbolic localization} 
Hyperbolic localization, see \cite{bradenHyperbolicLocalizationIntersection2003},  allows us to identify $\op{E}(\Delta)$ with the cohomology $H_c^\bullet(Z^-\times_ZZ^+)$ of the intersections of repelling and attracting sets.
Most importantly, in \Cref{thm:geometricdescriptionofcomposition} we realize the composition of extensions as a geometric convolution product on $Z^-\times_ZZ^+$ defined via the Drinfeld--Gaitsgory interpolation space $\widetilde{Z},$ see \cite{drinfeldTheoremBraden2014}, \cite{drinfeldAlgebraicSpacesAction2015}. We then consider specific examples.
\subsection{Partial flag varieties}
Let $Z=G/P$ be a partial flag variety\footnote{The same ideas apply to affine Grassmannians and affine flag varieties.} of a split reductive group equipped with the $\Gm$-action via a dominant cocharacter. The attracting and repelling sets are here the orbits of the Borel subgroup and its opposite, respectively. Their intersections $Z_x^-\times_Z Z_y^+=Z_x^-\cap Z_y^+$ are the famous \emph{open Richardon varieties}, see \cite{richardsonIntersectionsDoubleCosets1992} and \cite{speyerRichardsonVarietiesProjected2023} for a recent review.

We will define a `multiplication' on open Richardson varieties 
$$\psi: (Z^-\times_ZZ^+)\times_{Z^0}(Z^-\times_ZZ^+)\to Z^-\times_ZZ^+.$$
This induces a multiplication,  denoted $\star_\psi,$ on $H_c^\bullet(Z^-\times_ZZ^+)$ which provides an explicit geometric interpretation of compositions of extensions:
\begin{theorem*}[\Cref{Cor:concrete}]
    There is an isomorphism of algebras
    $$\op E(\Delta) \cong (H_c^\bullet(Z^-\times_ZZ^+),\star_\psi).$$
\end{theorem*}
\subsection{Grassmannians} By further specializing to Grassmannians, we will get a very concrete grip on the geometry of open Richardson varieties. Let $Z=\Gr(d,n)=\GL_n/P$ be the Grassmannian of $d$-planes in $n$-space. 
We label the fixed points $Z^0$ by $d$-element subsets $I,J\subseteq\{1,\dots,n\}$ and abbreviate the corresponding open Richardson variety by $\rich{I}{J}=Z^-_I\cap Z^+_J.$ We introduce a \emph{base change} map
$$\basechange{I}{J}: \rich{I}{J}\to \bG=\GL_d.$$
It allows to define two decompositions for $\rich{I}{J}$, the  \emph{\Gauss}\ and the \emph{\Deodhar}. They arise as the pullback via $\varphi$ of a stratification of $\bG$ (by $\opp\bB\times \bB$-orbits respectively $\bB\times \bB$-orbits).
Both decompose $\rich{I}{J}$ into cells parameterized via certain $d\times n$-matrices and each cell is isomorphic to some $(\Gm)^a\times \A^b$, see \Cref{sec:oppstratandfukaya} and \Cref{sec:stratandfukaya}.\footnote{These explicit parametrizations and decompositions extend the results of \cite{marshParametrizationsFlagVarieties2004b}.}

The \Gauss\ has the advantage that the parametrization of the cells are as explicit as possible: the entries in said $d\times n$-matrices are either zero, one, non-zero or arbitrary, following a simple rule. We will, however, mostly focus on the \Deodhar\ which is better suited for our cohomology calculations. The underlying combinatorics agrees with the original Deodhar decomposition from \cite{deodharGeometricAspectsBruhat1987a}, see \Cref{LA} for a precise statement. 
Our main result is the following:
\begin{theorem*}[\Cref{thm:richstrateverything}] The \Deodhar\ of $\rich{I}{J}$ is a stratification, that is, the closure of each stratum is a union of strata.  
\end{theorem*}
We note that a similar result for the Deodhar decomposition is not known, see \cite{marcottCombinatoricsDeodharDecomposition2020}. For general flag varieties, even in type $A$, the Deodhar decomposition sometimes fails to be a stratification, see \cite{dudasNoteDeodharDecomposition2008}, \cite{ speyerRichardsonVarietiesProjected2023}. 
\subsection{Cohomology computations} We use the \Deodhar\ to build a chain complex $\richmodel{I}{J}$ which serves as an explicit model for $H_c^\bullet(\rich{I}{J})$. 
Geometrically, $\richmodel{I}{J}$ is a direct sum of the cohomologies of the strata which is, algebraically, described in terms of (exterior) nil-Coxeter algebras (\Cref{def:exteriornilcoxeter}).
We construct a spectral sequence $E_1^{p,n-p}=\richmodel{I}{J}\Rightarrow H_c^\bullet(\rich{I}{J})$ with explicit boundary maps $\modeld$. Using weight arguments, we obtain:
\begin{theorem*}[\Cref{thm:modelrichardsongrassmannian}] 
    The spectral sequence $E_1^{p,n-p}$ degenerates on page two and there is an isomorphism
    $$H^\bullet(\richmodel{I}{J},d)\cong H_c^\bullet(\rich{I}{J}).$$
\end{theorem*}
We actually prove stronger results in \Cref{sec:cohomologygld} providing a  model for the cohomology of unions of Bruhat cells in  $\GL_d$. This may be of independent interest. 

\subsection{Fukaya calculus} We will describe the combinatorics, geometry and cohomology of the \Gauss\ and \Deodhar\ with certain decorated Fukaya diagrams. These are inspired by \cite{makFukayaSeidelCategories2021}. 
The diagrams are introduced in \Cref{THEguys}. They encode for instance the following:
\begin{itemize}
\item a labelling of the parts in the decompositions by elements of the Weyl group of  $\bG$ (providing also a new convenient way of labelling the parts in the original Deodhar decomposition),
\item a parametrization of the parts in the decompositions, 
\item a labelling of a basis in the complex computing the cohomology of the Richardson varieties and an explicit formula for the involved differential. 
\end{itemize}
Fukaya diagrams hopefully can be used to interpret our results in Floer homology.

\subsection{Category $\pazocal{O}$ and mixed Hodge polynomials} 
In \Cref{catO}, we consider the graded version $\hat\cO$ category $\cO$ from \cite{beilinsonKoszulDualityPatterns1996}, \cite{Strgradings}.
The geometry of open Richardson varieties is related to the extension algebra of Verma modules, see \cite{kazhdanRepresentationsCoxeterGroups1979}, \cite{richeModularKoszulDuality2014}. We establish the following for the graded parabolic category $\hat\cO^\mp$.
\begin{theorem*} [\Cref{thm:catO}] There is an isomorphism of  bigraded vector spaces
$$\Ext_{\hat\cO^\mp}^{i,j}(\Delta^\mp(x\cdot0),\Delta^\mp(y\cdot0))\cong H^{\ell(x)-\ell(y)+i,\frac{1}{2}(\ell(x)-\ell(y)-j)}_c(Z^-_{\psi(x)}\times_ZZ_{\psi(y)}^+)
$$
which gives rise to a bigraded isomorphism between the Lie-theoretic and geometric standard extension algebras $\stdE(\Delta^\mp)$ and $\op{E}(\Delta)$.  
\end{theorem*}
In nice settings, for instance if open Richardson admit a cellular stratification, the \emph{mixed Hodge polynomial} $\mixPol_c$ of the open Richardson variety can be expressed in terms of the bigraded Poincar\'e polynomial $\mixPol^\mp$ of Lie theoretic standard extensions 
$$\mixPol_c(Z^-_{\psi(x)}\times_ZZ_{\psi(y)}^+;q,t)=(q^{\frac{1}{2}}t)^{\ell(x)-\ell(y)}\mixPol^\mp(x,y;q^{-\frac{1}{2}},t).$$
In particular, we are in the beautiful setup, see \cite[Appendix]{hauselMixedHodgePolynomials2008}, where
specializations, at $t=-1$ and at $q=1$, already encode important information:
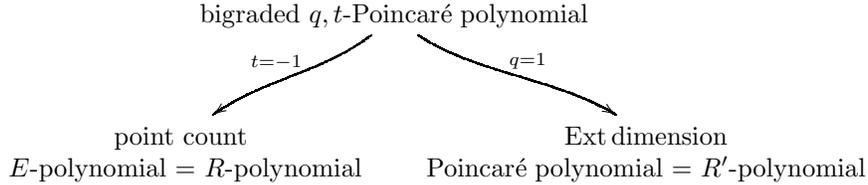
\begin{equation*}
\xymatrixcolsep{-6pc}\xymatrix{
&\text{bigraded $q,t$-Poincar\'e polynomial}
 \ar@(dl,ur)[dl]_{t=-1}  \ar@(dr,ul)[dr]^{q=1} \\
{\begin{array}[t]{c}
\text{point count        }\\
\text{$E$-polynomial = $R$-polynomial}
\end{array}}
&&{\begin{array}[t]{c}
 \op{Ext}\text{dimension}  \\
\text{Poincar\'e polynomial = $R'$-polynomial}
\end{array}}
}
\end{equation*}

The failure of the Gabber--Joseph conjecture can thus be explained by an asymmetry in $q$ and $t$ in the mixed Hodge polynomial of open Richardson varieties. 
There is a slightly different symmetry relating $q$ and $t$ in case of Koszul self-duality,  \cite{beilinsonKoszulDualityPatterns1996}; for instance in case of the full flag variety in type $A$ (or more generally for an appropriate direct sum of singular-parabolic blocks of $\cO$ which arises in the categorification of skew Howe duality, \cite{ICM}). For an interpretation of this symmetry in connection with knot homology and $q, t$-Catalan numbers we refer to \cite{galashinPositroidsKnotsCatalan2021}.

For Grassmannians, the mixed Hodge polynomials are determined recursively: 
\begin{theorem*}[\Cref{cor:recursiveformulabigradedpoincare} and \Cref{Es} for a more refined version]
Let $x,y\in\WPl$ and let $s$ be a simple reflection with $x>xs$. Then 
\begin{eqnarray}
    \mixPol^\mp(x,y)&=&
\begin{cases}
    \mixPol^\mp(xs,ys)&\text{if $y>ys$},\\
v^{-1}t \mixPol^\mp(xs,y)&\text{if $ys\notin\WPl$},\\
(v^{-1}t+v)\mixPol^\mp(xs,y)+\mixPol^\mp(xs,ys)&\text{otherwise.}\\
\end{cases}
\end{eqnarray}
\end{theorem*}

\subsection{Future Work} In a follow-up article we will discuss a $dg$-algebra upgrade and show that $\D_{Z^+}(Z)$ can be realized as category of $dg$-modules over $\richmodel{I}{J}$ using substantially the diagrammatics developed in this paper.

\subsection{Acknowledgments}
We like to thank Karin Erdmann, Yankı Lekili, Alexei Oblomkov and Jakob Scholbach for stimulating and useful discussions.  This work was supported by a DFG Leibniz Prize and the excellence cluster EXC-2047/1. 
\setcounter{tocdepth}{1} 
\tableofcontents

\addcontentsline{toc}{section}{Part I: Standard extension algebras: General constructions}
\section{Varieties with $\Gm$-action and extensions of standard objects}
We define the main player, the extension algebra of standard sheaves. 
\subsection{Varieties with $\Gm$-action after Drinfeld--Gaitsgory}
We first recall the approach of  Drinfeld--Gaitsgory \cite{drinfeldTheoremBraden2014}, \cite{drinfeldAlgebraicSpacesAction2015} to the Białynicki-Birula stratification \cite{bialynicki-birulaTheoremsActionsAlgebraic1973} and Braden's hyperbolic localization \cite{bradenHyperbolicLocalizationIntersection2003} for varieties $Z$ with $\Gm$-actions.

Throughout this article, we will make the following assumptions:
\begin{itemize}
	\item the ground field $k$ is algebraically closed,
	\item the variety $Z$ is smooth and projective,
	\item the $\Gm$-fixed points on $Z$ are isolated.
\end{itemize}
If not stated otherwise, coefficients are taken in an arbitrary commutative ring $\Lambda.$
\begin{remark}
We note that the results from \cite{drinfeldTheoremBraden2014}, \cite{drinfeldAlgebraicSpacesAction2015} work in greater generality, and by now, other generalizations of Braden's hyperbolic localization exist, see e.g. \cite{richarzSpacesMathbbGAction2018}, \cite{farguesGeometrizationLocalLanglands2021}. The standard extension algebra can be defined in a very general context, independent of perverse sheaves. We therefore do not use the perverse shifts of the standard objects. The results, however, do require some form of $\A^1$-homotopy invariance. For example, the assumptions hold for homotopy invariant algebraic $K$-theory and $K$-motives, see \cite{eberhardtKMotivesKoszulDuality2019}, \cite{eberhardtKtheorySoergelBimodules2022}.
\end{remark}

\begin{definition}{\rm
The \emph{space of $\Gm$-fixed points} $Z^0$, the \emph{attractor} $Z^+$ and the \emph{repeller} $Z^-$  are schemes representing the functors 
\begin{eqnarray*}
	&\Hom(S, Z^0)=\Hom_\Gm(S,Z),&\\
	&\Hom(S, Z^+)=\Hom_\Gm(\A^1_+\times S,Z)\text{ resp. }
	\Hom(S, Z^-)=\Hom_\Gm(\A^1_-\times S,Z)&
\end{eqnarray*}
for schemes $S$ over $k.$ Here, $\A^1_+$ and $\A^1_-$ denotes the affine line $\A^1$ equipped with the actions $\lambda\cdot x=\lambda x$ and $\lambda\cdot x=\lambda\inv x$ (for $\lambda\in \Gm$, $x\in \A^1$), respectively. Moreover, $\Hom_\Gm$ is the set of $\Gm$-equivariant morphisms and $S$ is equipped with the trivial $\Gm$-action.}
\end{definition}

Hence, $Z^0$ is the set of $\Gm$-fixed points, whereas $Z^+$ and $Z^-$ consists of pairs $(z,z')$ of a fixed point $z\in Z^0$ with a point $z'\in Z$ such that $\lim_{\lambda \to0}\lambda\cdot z'=z$ and  $\lim_{\lambda \to\infty}\lambda\cdot z'=z,$ respectively.

We have the obvious inclusions $i^\pm$ of fixed points and projection maps $q^\pm, p^\pm$:
\[\begin{tikzcd}
	{Z^0} & {Z^+} \\
	{Z^-} & Z.
	\arrow["{p^-}"', from=2-1, to=2-2]
	\arrow["{p^+}", from=1-2, to=2-2]
	\arrow["{i^+}"', shift right=1, from=1-1, to=1-2]
	\arrow["{i^-}", shift left=1, from=1-1, to=2-1]
	\arrow["{q^+}"', shift right=1, from=1-2, to=1-1]
	\arrow["{q^-}", shift left=1, from=2-1, to=1-1]
\end{tikzcd}\]
The fibers $Z_w^\pm=(q^\pm)\inv(\{w\})$ for $w\in Z^0$ are affine spaces and the maps $p^\pm$ are locally closed embeddings that define affine stratifications of $Z$
\begin{align*}
p^\pm:Z^\pm=\bigsqcup_{w\in Z^0}Z_w^\pm\hookrightarrow Z.
\end{align*}
In the following, for a variety $X$ denote by $\D(X)$ a formalism of triangulated categories (or stable $\infty$-categories) of constructible sheaves with a six-functor formalism evaluated at $X$. Prominent examples are for instance:
\begin{itemize}
	\item constructible sheaves on the complex points $\D_{\con}(X^{\an}(\C),\Lambda)$ for any ring $\Lambda$, 
	\item constructible $\ell$-adic sheaves $\D_{\con}(X/\overline{\F}_q,\Q_\ell).$
\end{itemize}
A sheaf $\Ff\in \D(Z)$ is called \emph{weakly $\Gm$-equivariant} if there is an isomorphism $a^*(\Ff)\cong p^*(\Ff)$ where $a,p:\Gm\times Z\to Z$ denote the action and projection map, respectively.
The category of $\Gm$-monodromic sheaves $\D_{\mon}(Z)$ is the full triangulated subcategory of $\D(Z)$ generated by weakly $\Gm$-equivariant sheaves. In this setting we obtain Braden's hyperbolic localization theorem,  \cite[Theorem 1]{bradenHyperbolicLocalizationIntersection2003}.
\begin{theorem}[Hyperbolic localisation]
	There is an equivalence of functors
	$$(q^+)_*(p^+)^!\simeq (q^-)_!(p^-)^*: \D_{\mon}(Z)\to \D(Z^0).$$
\end{theorem}
\subsection{Extension algebra of standard sheaves}
We will now study extensions of the \emph{standard object} $$\Delta=(p^+)_!(q^+)^*\Lambda\in \D(Z).$$ Hyperbolic localization allows us to interpret the extensions via the compactly supported cohomology of the intersections of attractors and repellers.
\begin{theorem}[Extspace]\label{Extspace}
	There is an isomorphism
	$$\Hom_{\D(Z)}(\Delta,\Delta[n])\stackrel{\sim}{\to} H^n_c(Z^-\times_ZZ^+).$$
\end{theorem}
\begin{proof}
Base change for the Cartesian diagram (with the obvious morphisms)
\[\begin{tikzcd}
	{Z^-\times_ZZ^+} & {Z^+} \\
	{Z^-} & Z
	\arrow["{k^-}", from=1-1, to=1-2]
	\arrow["{k^+}"', from=1-1, to=2-1]
	\arrow["{p^+}", from=1-2, to=2-2]
	\arrow["{p^-}"', from=2-1, to=2-2]
\end{tikzcd}\]
together with hyperbolic localization and adjunction gives isomorphisms
	\begin{align*}
		\Hom_{\D(Z)}((p^+)_!(q^+)^*\Q,(p^+)_!(q^+)^*\Q[n])&\cong\Hom_{\D(Z^0)}(\Q,(q^+)_*(p^+)_!(p^+)^!(q^+)^*\Q[n])\\
		&\cong \Hom_{\D(Z^0)}(\Q,(q^-)_!(p^-)_*(p^+)^!(q^+)^*\Q[n])\\
		&\cong \Hom_{\D(Z^0)}(\Q,(q^-)_!(k^+)_!(k^-)^*(q^+)^*\Q[n])\\
		&\cong \Hom_{\D(\pt)}(\Q,\pi_!\pi^*\Q[n])\\
		&= H^n_c(Z^-\times_ZZ^+)
	\end{align*}
	of vector spaces, where $\pi:Z^0\to \pt$ denotes the projection.
\end{proof}
\begin{remark}
The isomorphism from \Cref{Extspace} is in fact an isomorphism of $Z^0$-bigraded vector spaces by viewing $\Delta$ as a $Z^0$-graded object with $\Delta_w=(q^+_w)_*(p^+_w)^!\Lambda$ for $w\in Z^0.$ Then $H^n_c(Z^-\times_ZZ^+)=\Hom_{\D(Z)}(\Delta,\Delta[n])$ 
	is an object in the category of $Z^0$-bigraded vector spaces with components $\Hom_{\D(Z)}(\Delta_w,\Delta_{w'}[n])$ for $w,w'\in Z^0$. For example, the tensor product is defined as
		$$H^n_c(Z^-\times_ZZ^+)\otimes H^n_c(Z^-\times_ZZ^+)=\bigoplus_{w\in Z^0}H^n_c(Z^-\times_ZZ^+_w)\otimes H^n_c(Z^-_w\times_ZZ^+).$$
\end{remark}
The next goal is to understand the composition of extensions geometrically.
For this, we make use of the \emph{Drinfeld--Gaitsgory interpolation} space $\widetilde{Z}$, see \cite[Section 2.2]{drinfeldAlgebraicSpacesAction2015}, whose properties we recall briefly: 
\begin{itemize}
	\item $\widetilde Z$ is  a locally closed subvariety of $\A^1\times Z\times Z$. We denote the projection onto the first component by $p:\widetilde Z\to \A^1$ and set $\widetilde Z_V=p\inv(V)$ for $V\subset \A^1.$
	\item Over $\Gm,$ the space $\widetilde Z$ is the graph of the $\Gm$-action $$\widetilde Z_\Gm=\{(\lambda,z,\lambda\cdot z)| \lambda\in \Gm, z\in Z\}\cong \Gm\times Z.$$
	\item Over $\{0\},$ one has
	$\widetilde Z^0=\{0\}\times Z^+\times_{Z^0}Z^-\cong Z^+\times_{Z^0}Z^-.$
\end{itemize}
Intuitively speaking, the space $\widetilde Z$ interpolates between the variety $Z$ in the generic fiber and the disjoint union of tubular open neighborhoods of the fixed points in the special fiber. Namely, 
$$Z^+\times_{Z^0}Z^-=\bigsqcup_{w\in Z^0}Z_w^+\times Z_w^-.$$
Here, $Z_w^+\times Z_w^-\cong T_w^+Z\oplus T_w^-Z= T_wZ,$ where $T_w^\pm Z$ denotes the subspace of the tangent space $T_wZ$ at $w$ on which $\Gm$ acts with positive/negative weights.

\begin{remark}
	In the context of flag varieties, the families $\widetilde{Z}$ and $Z^-\times_Z \widetilde Z\times_ZZ^+$ were considered before in the context of standard monomial theory and Frobenius splitting in the work of Brion--Lakshmibai \cite{brionGeometricApproachStandard2003} and Brion--Polo \cite{brionLargeSchubertVarieties2000} to construct for instance degenerations of the diagonal for Richardson varieties.
\end{remark}

There are two natural projection maps $\widetilde{Z}\to Z$ and we consider the  fiber product $Z^-\times_Z \widetilde Z\times_ZZ^+.$ This yields a family over $\A^1$ with special and generic fibers
$$Z^-\times_Z Z^+\times_{Z^0}Z^-\times_Z Z^+ \text{ and }Z^-\times_ZZ^+, \text{ respectively.}$$
We can use this to define a product
\[\begin{tikzcd}[row sep=0pt]
	{\star:H^\bullet_c(Z^-\times_ZZ^+)\otimes H^\bullet_c(Z^-\times_ZZ^+)} & {H^\bullet_c(Z^-\times_ZZ^+\times_{Z^0}Z^-\times_ZZ^+)} \\
	{\phantom{\star:H^n_c(Z^-\times_ZZ^+)\otimes H^\bullet_c(Z^-\times_ZZ^+)}} & {H^\bullet_c(Z^-\times_ZZ^+)\phantom{\times_{Z^0}Z^-\times_ZZ^+)}.}
	\arrow["\boxtimes", from=1-1, to=1-2]
	\arrow["\cosp", from=2-1, to=2-2]
\end{tikzcd}\]
where $\cosp$ denotes the {\it cospecialization map}, in the sense of \cite[6.2.7]{schiederMonodromyVinbergFusion2019}, associated to the family $Z^-\times_Z \widetilde Z\times_ZZ^+\to \A^1.$
\begin{remark}\label{rem:cospezialisation}
	In the case $k=\C$ and singular cohomology on complex points, cospecialization admits the following explicit definition. For a family $p: X\to \A^1$ which trivializes to $X_1\times \Gm$ over $\Gm$, the map $\cosp$ is the composition
\[\begin{tikzcd}
	{H_c^n(X_0)} & {H^{n+1}_c(X_{\R_{>0}})} & {H^{n}_c(X_1)}
	\arrow["\delta", from=1-1, to=1-2]
	\arrow["{\pi_!}", from=1-2, to=1-3]
	\arrow["\cosp"', shift right=1, curve={height=6pt}, from=1-1, to=1-3]
\end{tikzcd}\]
where $\delta$ is the boundary map in the long exact sequence associated to the open and closed decomposition $X_{\R_{>0}}\sqcup X_0=X_{\R_{\geq}0}$ and $\pi_!$ denotes the pushforward along the projection $X_{\R_{>0}}=X_1\times \R_{>0}\to X_1.$
	This is the linear dual of the \emph{specialization} map in Borel--Moore homology
	$\spe: H^{\BM}_n(X_1)\to H^{\BM}_n(X_0)$
	which, intuitively speaking, maps a cycle 
	$\sigma\in H^{\BM}_n(X_1)$ to the limit $\lim_{\lambda \to 0} \lambda\cdot \sigma,$ see \cite[Section 2.6.30]{chrissRepresentationTheoryComplex2010}.
\end{remark}

This product describes composition of extensions of the standard object $\Delta$:
\begin{theorem}[Standard extension algebra]\label{thm:geometricdescriptionofcomposition}
	The following diagram commutes
\[\begin{tikzcd}
	{\Hom_{\D(Z)}(\Delta,\Delta[\bullet])\otimes \Hom_{\D(Z)}(\Delta,\Delta[\bullet])} & {\Hom_{\D(Z)}(\Delta,\Delta[\bullet])} \\
	{H^\bullet_c(Z^-\times_ZZ^+)\otimes H^\bullet_c(Z^-\times_ZZ^+)} & {H^\bullet_c(Z^-\times_ZZ^+)}
	\arrow["\star", from=2-1, to=2-2]
	\arrow["\circ^{\operatorname{op}}", from=1-1, to=1-2]
	\arrow["\wr", from=1-1, to=2-1]
	\arrow["\wr", from=1-2, to=2-2]
\end{tikzcd}\]
where $f\circ^{\operatorname{op}}g=g\circ f$ denotes the opposite composition.
\end{theorem}
\begin{proof}
	To prove \Cref{thm:geometricdescriptionofcomposition},
	 one follows Drinfeld--Gaitsgory's proof of the hyperbolic localisation theorem \cite{drinfeldAlgebraicSpacesAction2015}, \cite{ drinfeldTheoremBraden2014}. There, the interpolation space $\widetilde Z$ is used to construct the (co-)units of an adjunction between the functors 
	 $(p^+)_!(q^+)^*$ and $(q^-)^!(p^-)_*.$ The theorem follows then from a standard yoga of six functors and these (co-)units.
\end{proof}
\begin{remark}
	There is also a differential graded (or $\infty$-categorical) version of the above discussion which yields a homotopy equivalence
	$$\Maps_{\D(Z)}(\Delta,\Delta)\stackrel{\simeq}{\to} R\Gamma_c(Z^-\times_ZZ^+)$$
	between the mapping space of the standard object and the complex of compactly supported cochains on $Z^-\times_ZZ^+.$
\end{remark}
\begin{definition}\label{Extalg}{\rm
We call $\stdE(\Delta):={\Hom_{\D(Z)}(\Delta,\Delta[\bullet])}$ the \emph{algebra of extensions of standards} or the \emph{standard extension algebra}.}
\end{definition}

\subsection{Bigradings on cohomology}\label{sec:bigradings}
It is very useful to equip the cohomology groups considered here with an additional grading. In fact there is (in many cases) a $\mathbb{Z}$-bigraded enhancement $H^{\bullet,\bullet}(X)$ of $H^\bullet(X)$, where the first grading is the usual cohomological grading, and the second, additional, grading comes from a notion of weights. For example, $H^{\bullet,\bullet}(\Proj^1)$ lives in bidegrees $(0,0)$ and $(2,1).$ We write $H^{\bullet,\bullet}(\Proj^1)=\Z\oplus \Z(-1)[-2].$

For varieties $X/\C$ and coefficients in $\C,$ the additional grading comes from the weight filtration of the associated mixed Hodge structure.
For $X/\F_p$ and coefficients in $\overline\Q_\ell$, it  comes from (generalized) eigenspaces of the Frobenius action. 

For $X/S$ over a general base $S$ with a cellular stratification (and so strata of the form $(\Gm)^n\times \A^m$) with integral coefficients, one obtains a bigraded theory using reduced motivic cohomology, which is, roughly speaking, the motivic cohomology of $X$ modulo the higher motivic cohomology of the base $S.$ More precisely, 
$$H^{p,q}(X)=\Hom_{\operatorname{DM}_{\operatorname{r}}(S)}(M(X/S),\Z(q)[p])$$
where $\operatorname{DM}_{\operatorname{r}}(S)$ is the category of reduced motivic sheaves introduced in \cite{eberhardtIntegralMotivicSheaves2023}, $M(X/S)$ is the motive of $X/S$ and $(i)$ denotes the $i$-th Tate twist. We also refer to \cite{soergelPerverseMotivesGraded2018} and \cite{eberhardtMixedMotivesGeometric2019} for the cases $\Lambda=\Q$ and $\Lambda=\F_p,$ respectively.
The main examples we discuss here, such as flag varieties or open Richardson varieties for Grassmannians, admit cellular stratifications, see \Cref{sec:flagvarieties} and \Cref{thm:richstrateverything}.

The same considerations apply to cohomology with compact support, homology and Borel--Moore homology, for instance $H_c^{\bullet,\bullet}(\A^1)= \Z(-1)[-2]$ is concentrated in bidegree $(2,1)$ and $H_c^{\bullet,\bullet}(\Gm)=\Z[-1]\oplus \Z(-1)[-2]$ is lives in bidegree $(2,1)$ and $(1,0)$.  As a consequence the standard extension algebra from \Cref{Extalg} is in fact a  $\Z^2$-graded algebras. 

Moreover, for a variety $X/S$ we can consider the bigraded Poincar\'e polynomial
$$\mixPol_c(X;q,t)=\sum_{i,j}\dim(H_c^{i,j}(X,\Q))q^jt^i.$$
If $X$ is defined over $\Z$ (or, some subring of $\C$) and has a cellular stratification, the polynomial $\mixPol_c(X;q,t)$ coincides with the (compactly supported) \emph{mixed Hodge polynomial}, see \cite{hauselMixedHodgePolynomials2008}, and we call it like that. In particular, by substituting $t=-1,$ we obtain the \emph{point counting polynomial}
$P_X(q)=\mixPol_c(X;q,-1)=|X(\F_q)|.$
\addcontentsline{toc}{section}{Part II: Standard extension algebras: Flag varieties}
\section{Composition of extensions for flag Varieties}\label{sec:flagvarieties}
In this section we consider the special case of partial flag varieties and give an alternative description of the multiplication in the extension algebra of standards.
\subsection{Notations}\label{sec:notationsflagvarieties}
Let $G\supset B\supset T$ be a split reductive algebraic group over $k$ with a Borel subgroup $B$ and maximal torus $T.$ Let $U\subset B$ be the unipotent radical.
Let $\pW=N_G(T)/T\supset \pS$ be the Weyl group with the set $\pS$ of simple reflections with respect to the choice of $B.$ Let $w_0\in \pW$ be the longest element. Denote by $\opp B=B^{w_0}$ the opposite Borel subgroup, similarly $\opp U=U^{w_0},$ etc.

Denote by $X(T)=\Hom_{grp}(T,\Gm)$ and  $Y(T)=\Hom_{grp}(\Gm,T)$ the character and cocharacter lattice. Let $X(T)\supset\Phi\supset\Phi^\pm, \Delta$ be the set of roots, of positive/negative roots and of simple roots.

We fix a standard parabolic subgroup $B\subset P\subset G$ and consider the flag variety $Z=G/P.$ Denote by $\pS_P\subset\pW_P\subset \pW$ the corresponding simple reflections and parabolic subgroup of the Weyl group. Denote by $\pW^P$ and $\opp \pW{}^P$ the set of shortest and longest coset representatives of $\pW/\pW_P,$ respectively.

Let $\eta\in Y(T)$ be a generic and dominant cocharacter, that is, $\langle\eta, \alpha\rangle > 0$ for all $\alpha\in \Phi^+.$ We consider the $\Gm$-action on $Z$ induced by $\eta.$ Then the fixed points, attractor and repeller sets can be described as
$$Z^0=\,\bigsqcup_{\mathclap{w\in \pW/\pW_P}}\, wP/P,\, Z^+=\,\bigsqcup_{\mathclap{w\in \pW/\pW_P}}\, BwP/P\text{ and }Z^-=\,\bigsqcup_{\mathclap{w\in \pW/\pW_P}}\,\opp BwP/P.$$
In the following, we will identify $Z^0$ with $\pW/\pW_P.$ 
The intersections $Z^-_x\cap Z^+_y$ are called \emph{open Richardson varieties} \cite{richardsonIntersectionsDoubleCosets1992} and are non-empty if and only if $x\leq y.$ 

The question of finding explicit parametrizations was initiated in   \cite{marshParametrizationsFlagVarieties2004b} for full flag varieties and studied in connection with total positivity questions.  
\subsubsection{Parametrizations}
For $w\in \pW$ let $U_w=U\cap w\opp{U}w^{-1}$ and $\opp U_w= (U_{w_0w})^{w_0}=\opp U\cap wUw^{-1}.$ In the following, by abuse of notation, for $w\in \pW/\pW_P$ we denote $U_w=U_{\dot w}$ for the shortest coset representative $\dot w\in \pW^P$ of $w$, and $\opp{U}_w=\opp{U}_{\dot w}$ for the longest coset representative $\dot w\in {\opp\pW}{}^P$ of $w.$ 
We obtain affine parametrizations 
\begin{eqnarray*}
 U_w\stackrel\sim\to U_wwP/P=Z^+_w&\text{and} &\opp U_w\stackrel\sim\to \opp U_wwP/P=Z^-_w
 \end{eqnarray*}
of $Z_w^\pm$. They provide two parametrizations of an open neighborhood $V_w$ of $wP/P$
\begin{equation}\label{eq:swapmap}\begin{tikzcd}
	{\opp U_w\times U_w} & {V_w\defi \opp U_w  U_wwP/P= U_w  \opp U_wwP/P} & { U_w\times \opp U_w}
	\arrow["\sim", from=1-1, to=1-2]
	\arrow["\sim"', from=1-3, to=1-2]
	\arrow["{\sigma_w}"', shift right=2, curve={height=6pt}, from=1-1, to=1-3]
\end{tikzcd}\end{equation}
which induce the indicated isomorphism $\sigma_w.$ We note that this is in general (in fact already for $\GL_3/B$) not the obvious isomorphism, since $\sigma_w(\opp u, u)=(u,\opp u)$ if and only if $u$ and $\bar u$ commute.  
For each pair $u\in U_w$ and $\opp u\in \opp U_w$ there are however unique $u'\in U_w$ and $\opp u'\in \opp U_w$ such that $\opp u'uwP/P=u'\opp uwP/P.$ Namely, they are defined via $\sigma_w(\bar u, u\inv)=((u')\inv, \bar u').$

Abbreviate $V=\bigsqcup_{w\in Z^0} V_w$. The inclusion map $\iota: V\to Z$ is open, $\Gm$-equivariant and induces an isomorphism on fixed points, attractors and repellers.
\subsection{Composition of extensions}\label{sec:conjectureofgeometricalcompositionformula}
The explicit parametrizations of the cells in the flag variety allows for another way to compose elements in $H^\bullet_c(Z^-\times_ZZ^+)$ which works without using the Drinfeld--Gaitsgory interpolation space $\widetilde Z$ and cospecialization. It is much more feasible for explicit computations. 
\begin{definition}[Multiplicative Structure]{\rm
	For $w\in \pW/\pW_P,$ we define the map
\begin{equation}\label{eq:multiplicationofrichardsons}
	\begin{tikzcd}[row sep=0pt]
	{\psi_w\colon Z^-\times_ZZ_w^+\times Z_w^-\times_ZZ^+} & {Z^-\times_ZZ^+} \\
	{\phantom{\psi_w:}(z^-,uwP/P,\opp uwP/P,z^+)} & {(\opp u'z^-,u'z^+)}.
	\arrow[from=1-1, to=1-2]
	\arrow[maps to, from=2-1, to=2-2]
\end{tikzcd}
\end{equation}
Here, $(u,\opp u)\in U_w\times\opp U_w$ and $(u',\opp u')\in U_w\times\opp U_w$ is the unique pair such that $\opp u'uwP/P=u'\opp uwP/P.$ These maps assemble into
\begin{equation}\label{multspaces}
\psi\colon Z^-\times_ZZ^+\times_{Z^0} Z^-\times_ZZ^+\to Z^-\times_ZZ^+
\end{equation}
which induces a product
\begin{eqnarray*}
\begin{tikzcd}[row sep=0pt]
	{\star_{\psi}\colon H^\bullet_c(Z^-\times_ZZ^+)\otimes H^\bullet_c(Z^-\times_ZZ^+)} & {H^\bullet_c(Z^-\times_ZZ^+\times_{Z^0}Z^-\times_ZZ^+)} \\
	{\phantom{\star:H^n_c(Z^-\times_ZZ^+)\otimes H^n_c(Z^-\times_ZZ^+)}} & {H^\bullet_c(Z^-\times_ZZ^+)\phantom{\times_{Z^0}Z^-\times_ZZ^+}.}
	\arrow["{\psi_!}", from=2-1, to=2-2]
	\arrow["\boxtimes", from=1-1, to=1-2]
	\arrow["\sim"', from=1-1, to=1-2]
\end{tikzcd}
\end{eqnarray*}
}
\end{definition}
\begin{theorem}[Product formula]\label{thm:twocompositionsagree} The products $\star$ and $\star_{\psi}$ on $H^\bullet_c(Z^-\times_ZZ^+)$ agree.
\end{theorem}
\begin{proof}
	Consider the Drinfeld interpolation space $\widetilde{V}\subset \A^1\times V\times V$ associated to $V.$ 
	
	Then $\iota$ induces a map $\widetilde{\iota}:\widetilde{V}\to  \widetilde{Z}$ which restricts in the special fiber to the equality $\widetilde{V}_0=Z^+\times_{Z^{0}}Z^-=\widetilde{Z}_0$ and over $\Gm$ to a map $\Gamma(\Gm,V)\to \Gamma(\Gm,Z)$ between the graphs of the $\Gm$-action. It can be trivialized to a map $\id\times \iota: \Gm\times V\to \Gm\times Z.$

	By passing to the fiber product $Z^-\times_Z-\times_ZZ^+$ we hence obtain the diagram
	\[\begin{tikzcd}
		{H^n_c(Z^-\times_ZZ^+\times_{Z^0}Z^-\times _ZZ^+)} & {H^{n}_c(Z^-\times _ZZ^+)} \\
		{H^n_c(Z^-\times_ZZ^+\times_{Z^0}Z^-\times _ZZ^+)} & {H^{n}_c(Z^-\times_ZV\times_ZZ^+).}
		\arrow["{\iota'_!}"', from=2-2, to=1-2]
		\arrow[Rightarrow, no head, from=2-1, to=1-1]
		\arrow["\cosp'", from=2-1, to=2-2]
		\arrow["\cosp", from=1-1, to=1-2]
	\end{tikzcd}\]
	This diagram is commutative. In the complex analytic case, this follows easily from the definition of the cospecialization map, see \Cref{rem:cospezialisation}, and the fact that the boundary map in the long exact sequence commutes with pushforwards along open maps. Note that the trivializations over $\Gm$ of $\widetilde{Z}$ and $\widetilde{V}$ are compatible.
	
	In contrast to $Z^-\times_Z\widetilde{Z}\times_Z Z^+$, the family $Z^-\times_Z\widetilde{V}\times_ZZ^+$ admits a trivialization over $\A^1$ which we construct now.
	Let $w\in \pW/\pW_P,$ $u\in U_w$ and $\opp u\in \opp U_w.$ For $\lambda\in\A^1,$ there are unique $u_\lambda\in U_w$ and $\opp u_\lambda\in \opp U_w$ such that 
	$$\opp u_\lambda u^{\eta(\lambda)} yP/P=u_\lambda\opp uyP/P\in Z.$$
	This is well-defined for $\lambda=0,$ since $\lim_{\lambda\to 0}u^{\eta(\lambda)}=1.$
	Similarly, it follows that $\lim_{\lambda\to 0}\opp u_\lambda^{\eta(\lambda\inv)}=1$ using $\lim_{\lambda\to 0}\bar u^{\eta(\lambda\inv)}=1.$
	In this way, we interpolate between $\opp u_0=\opp u$ and $u_0=1$ and $u_1=u'$ and $\opp u_1=\opp u'$ with $\opp u'u wP/P= u'\opp u wP/P.$ 
	
	Consider $\varphi: \A^1\times Z^-\times_Z Z^+\times_{Z^0} Z^-\times_Z Z^+\to Z^-\times_Z\widetilde{V}\times_ZZ^+$ defined via 
	\begin{eqnarray*}
	\A^1\times Z^-\times_Z Z^+_w\times Z^-_w\times_Z Z^+&\to& Z^-\times_Z\widetilde{V}\times_ZZ^+ \\
	(\lambda, z^-, uyP/P, \bar uyP/P,z^+)&\mapsto& (\opp u_\lambda^{\eta(\lambda\inv)}z^-,\lambda,\opp u_\lambda^{\eta(\lambda\inv)}uyP/P, u_\lambda \bar uyP/P, u_\lambda z^+).
	\end{eqnarray*}
	The above discussion implies that $\varphi$ is a well-defined isomorphism. It induces an isomorphism $\psi':Z^-\times_ZZ^+\times_{Z^0}Z^-\times _ZZ^+\to Z^-\times_ZV\times _ZZ^+$ between the fibers over $0$ and $1$ of $Z^-\times_Z\widetilde{V}\times_ZZ^+.$ The cospecialization map is hence simply given by $\cosp'=\psi'_!.$ The statement follows since $\iota'\psi'=\psi.$
\end{proof}
\begin{corollary}[Concrete extension algebra]\label{Cor:concrete}
Via \Cref{thm:geometricdescriptionofcomposition} we obtain in case of partial flag varieties $Z$ a canonical isomorphism of graded algebras
\begin{eqnarray*}
 \stdE(\Delta)&\cong&(H^\bullet_c(Z^-\times_ZZ^+),\star_{\psi}).
\end{eqnarray*}
\end{corollary}

\addcontentsline{toc}{section}{Part III: Standard extension algebras: Grassmannians}
\section{Grassmannians and normal forms for Schubert cells} 
From now on we restrict ourselves to Grassmannians. As preparation for the subsequent sections we recall concrete normal forms for elements in Bruhat cells. 
\subsection{Notations}
Let $0\leq d\leq n\in \Z$ and let now $Z=\Gr(d,n)$ be the Grassmannian of $d$-dimensional subspaces in $k^{1\times n}.$ 

We let $G=\GL_n$ act on $k^{1\times n}$ from the left via $g\cdot v=vg\tr.$  
The stabilizer of the subspace $k^d\subset k^n$ spanned by the first $d$ coordinate vectors is the standard parabolic subgroup $P$ containing the Borel subgroup $B\subset G$ of upper-triangular matrices. We identify $\Gr(d,n)$ with  $G/P$ via this action. 

The corresponding parabolic subgroup $\pW_P\subset \pW$ of the Weyl group $\pW$ is isomorphic to $S_d\times S_{n-d}\subset S_n.$ Let $\pS=\{s_1=(1,2),s_2=(2,3),\dots, s_{n-1}=(n-1,n)\}\subset \pW$ be the set of simple reflections.
We still denote by $\pW^P$ the set of representatives of shortest length for the cosets ${\pW/\pW_P}$.
\subsection{$\Gm$-action and combinatorics of fixed points}\label{fixedpoints}
We let $\Gm$ act on $k^{1\times n}$ via $\lambda\cdot(a_i)=(\lambda^{-i}a_i).$ This induces a $\Gm$-action on $Z=\Gr(d,n)$. The fixed points are exactly the subspaces spanned by coordinate vectors.

We abbreviate $\firstn{k}=\{1,\dots,k\}$ and denote by $\setbinom{\firstn{n}}{d}$ the set of subsets of $\firstn{n}$ of cardinality $d$. For $I=\{i_1<\dots<i_d\}\in \setbinom{[n]}{d}$ we denote by $e_I=(I_n)_{I\times [n]}\in k^{d\times n}$ the matrix whose $k$-th row is the $i_k$-th standard basis vector and by $[e_I]\in Z$ the subspace spanned by the rows of $e_I.$

We will make use of the identifications
\begin{equation}\label{identifications}
\begin{tikzcd}[row sep=0]
	\pW^P&{\pW/\pW^P} & {\setbinom{\firstn{n}}{d}} & {Z^0} \\
	w&{w\pW^P} & {w(\firstn{d})} & {w\cdot[e_{\firstn{d}}]=[e_{w(\firstn{d})}]}
	\arrow[maps to, from=2-1, to=2-2]
	\arrow[maps to, from=2-2, to=2-3]
    \arrow[maps to, from=2-3, to=2-4]
	\arrow["\sim",from=1-1, to=1-2]
	\arrow["\sim",from=1-2, to=1-3]
    \arrow["\sim",from=1-3, to=1-4]
\end{tikzcd}
\end{equation}
and label fixed points by subsets $I\in \setbinom{\firstn{n}}{d}.$ Via \eqref{identifications}, $\setbinom{\firstn{n}}{d}$ inherits a length function and an ordering from the Bruhat order on $\pW^P$. Concretely, if $I=\{i_1<\dots<i_d\},\, J=\{j_1<\dots<j_d\}\in \setbinom{\firstn{n}}{d},$ then $I\leq J$ if and only if $i_k\leq j_k$ for all $k.$
In this case, $\ell(J)-\ell(I)=\sum_k(j_k-i_k).$

\subsection{Normal forms for Schubert cells} \label{sec:parabruhat}
In the following, it will be convenient to describe subspaces in $W\in Z$ as the row span $W=[X]$ of matrices $X\in k^{d\times n}$. 
This description is not unique since performing any row operation on $X$ yields the same subspace. However, for some subsets of $Z,$ such as Bruhat cells or open Richardson varieties, we have preferred choices for $X$ depending on the specific use we have in mind.  In the following we will present several such  normal forms.

Let $I=\{i_1<\dots<i_d\}\in \setbinom{\firstn{n}}{d}$ and $w\in \pW$ with $I=w(\firstn{d}).$ The attractor sets, repeller sets  and open neighborhoods of the fixed point labeled by $I$ are given by
\begin{equation}\label{eq:parabruhat}Z_I^+=B\cdot[e_I],\, Z_I^-=\opp B\cdot[e_I] \text{ and }V_I=(\opp B B)\cdot[e_I],\text{ respectively}.\end{equation}
Here, $B, \opp B \subset G$ are the Borel subgroups of upper respectively lower diagonal matrices in $G$. The sets $Z_I^\pm$ are the (positive and negative) \emph{Bruhat cells} and their closures yield \emph{Schubert varieties}.

Consider the following sets of matrices
\begin{align}
	M_I&=\setbuild{X\in k^{d\times n}}{X_{\firstn{d}\times I}=I_d},\\
	M_I^+&=\setbuild{X\in M_I}{X_{k,l}=0\text{ for all } i_k<l}\text{ and }\label{eq:MIplus}\\
	M_I^-&=\setbuild{X\in M_I}{X_{k,l}=0\text{ for all } l<i_k}\label{eq:MIminus}.
\end{align}
Here, by $X_{\firstn{d}\times I}$ we denote the submatrix given by the columns in $X$ indexed by $I$. For $X\in M_I$ we call the corresponding entries $X_{i_k,k}=1$ \emph{leading entries}. Then $M_I^+$ (or $M_I^-$) consists of matrices $X\in M_I$ such all entries in $X$ to the right (or left) of the leading entries in $X$ are zero, see \Cref{sec:fukayadiagrams} for an explicit example.  There, the $*$ indicate an affine parametrization.  
The subsets $M_I^\pm$ and $M_I$ yield normal forms for the Bruhat cells $Z^\pm_I$ and the open neighborhood $V_I,$ respectively.

Formally, denote by $V_{d,n}\subset k^{d\times n}$ the subset of full-rank matrices.
Then the map
$$[_-]: V_{d,n}\rightarrow Z, X\mapsto [X]$$
is a torsor for the group $\bG=\GL_d$ which acts on $k^{d\times n}$ by left multiplication. 

\begin{definition}{\rm
For a subset $Z'\subset Z,$ {\it a normal form} is a section of $[_-]$ over $Z',$ or equivalently, a subset $V'\subset V_{d,n}$ for which $[_-]:V'\to Z'$ becomes an isomorphism.}
\end{definition}

\subsection{More sets of matrices}To construct normal forms, we will need explicit descriptions of certain subsets of matrices derived from $M_I^\pm.$ 

Denote by $\bT\subset \bB,\opp\bB\subset \bG=\GL_d$ the maximal torus of diagonal matrices, the Borel subgroup of upper-triangular matrices and its opposite inside $\bG=\GL_d$. Let $\bW=S_d$ be the Weyl group. Denote by $\bU\subset\bB$ and $\opp\bU\subset\opp\bB$ the unipotent radicals. Moreover, let 
\begin{align*}
	\bU_w&=\bU  \cap w\inv\opp \bU w=\setbuild{R\in \bU}{R_{k,l}=0\text{ if }w\inv(k)<w\inv(l)}\text{ and}\\
	\opp \bU_w&=\opp \bU  \cap w\inv\opp \bU w=\setbuild{L\in \opp\bU}{L_{k,l}=0\text{ if }w\inv(k)<w\inv(l)}.
\end{align*}
Then we obtain the following descriptions
\begin{equation}
	\label{eq:wUcapoppUw}
	\begin{split}
	w\bU_w&=w\bU\cap \opp \bU w\\
	&=\setbuild{A\in \bG}{A_{k,w\inv(k)}=1\text{ and } A_{k,l}=0\text{ if }l<w\inv(k) \text{ or }k<w(l)},
	\end{split}
\end{equation}
\begin{equation}
	\label{eq:woppUcapoppUw}
	\begin{split}
	w\opp \bU_w&=w \opp \bU  \cap \opp \bU w\\
	&=\setbuild{A\in \bG}{A_{k,w\inv(k)}=1\text{ and } A_{k,l}=0\text{ if }w\inv(k)<l \text{ or }k<w(l)}.
	\end{split}
\end{equation}
\begin{proposition}\label{prop:elementsofUMIminus}
Let $I=\{i_1<\dots<i_d\}$ and
	 $X\in k^{d\times n}.$ Then 
	 \begin{eqnarray}\label{eq:elementsofUMIminus}
	 X\in\bU M_I^- &\Leftrightarrow&
X_{k,l} =\begin{cases}
	0, & \text{if $l<i_k$;}\\
	1, & \text{if $l=i_k$.}\\
 \end{cases}
\end{eqnarray}
\end{proposition}
\begin{proposition}\label{prop:elementsofoppUMIminus}
	Let $J=\{j_1<\dots<j_d\}$ and $X\in k^{d\times n}.$ Then  
\begin{eqnarray}\label{eq:elementsofoppUMJplus}
X\in w\bT\opp\bU_w M_J^+&\Leftrightarrow&
X_{k,l} \begin{cases}
	\neq 0, & \text{if $l=j_{w\inv(k)}$;}\\
	=0, & \text{if $l>j_{w\inv (k)}$;}\\
	=0, & \text{if $l=j_{w\inv(k')}$ for some $k'>k$}.\\ 
 \end{cases}
\end{eqnarray}
The set $w\bT\opp\bU M_J^+$ arises by removing the third condition from \eqref{eq:elementsofoppUMJplus}.
\end{proposition}
\begin{remark}
The overlap of the second and third case in \eqref{eq:elementsofoppUMJplus} could be avoided by 
adding the condition $w\inv(k')<w\inv(k)$ to the third case.
\end{remark}

We provide no proofs of \Cref{prop:elementsofUMIminus,prop:elementsofoppUMIminus}, but instead prove the following more complicated description of the set $w\bT\bU_w M_J^+.$
\begin{proposition}\label{prop:elementsofUMJplus} 
	Let $J=\{j_1<\dots<j_d\}$ and $X\in k^{d\times n}.$ Then 
\begin{eqnarray}\label{eq:elementsofUMJplus}
X\in w\bT\bU_w M_J^+&\Leftrightarrow&
X_{k,l} \begin{cases}
	\neq 0, & \text{if $l=j_{w\inv(k)}$;}\\
	=\eqref{eq:thedeterminedentries}, & \text{if $l>j_{w\inv(k)}$ and $l\not\in J$;}\\
	=0, &\text{if $l=j_{l'}$ with $l'<w\inv(k)$;}\\
	=0, &\text{if $l=j_{l'}$ with $k<w(l')$.}
 \end{cases}
\end{eqnarray}
Here, by $X_{k,l}=(\star)$ we denote that the entry $X_{k,l}$ is the value of a polynomial involving entries of the form $X_{k',l}$ for $k'<k$ and entries in $X_{\firstn{d}\times J};$
explicitly,
\begin{align}
	\tag{$\star$}\label{eq:thedeterminedentries}
	X_{k,l}&=\sum_{k'<k} ({L_{(k)}})_{k,k'}X_{k',l},
\end{align}
where $L_{(k)}$ is a certain matrix depending polynomially of $X_{\firstn{d}\times J}.$
The set $w\bT\bU M_J^+$ arises by removing the fourth condition from \eqref{eq:elementsofUMJplus}.
\end{proposition}
\begin{proof} It suffices to describe the set $w\bU_w M_J^+$. The additional factor of $\bT$ simply multiplies each row by  non-zero scalar and changes the first condition in \eqref{eq:elementsofUMJplus} from $=1$ to $\neq 0.$ Now, let $X=AX'$ for $X'\in M_J^+$ and $A=w\bU_w.$ Moreover, write $A=wR=Lw$ for $R\in \bU$ and $L\in \opp\bU.$ We consider the entry $X_{k,l}.$

{\it Case  $l\in J$:} This corresponds to the first, third and fourth condition in \eqref{eq:elementsofUMJplus}.
Since $(AX')_{\firstn{d}\times J}=A$ these conditions come from the description of $w\bU_w$ in \eqref{eq:wUcapoppUw}. 

{\it Case $l\not\in J$:} Now, let $l\not\in J.$ The entry $X_{k,l}$ has the following form, where we use the description of $w\bU_w$ in \eqref{eq:elementsofUMJplus} and of $M_J^+$ in \eqref{eq:MIplus}:
\begin{align}\label{eq:complicatedentries}
	X_{k,l}=\sum_{l'}A_{k,l'}X'_{l',l}=X'_{w\inv (k),l}+\sum_{\substack{w\inv(k)<l'\\ w(l')<k\\l<j_{l'}}}A_{k,l'}X'_{l',l}.
\end{align}

{\it Subcase $l<j_{w\inv(k)}$:} Here, $X'_{w\inv (k),l}$ can take arbitrary values. Consequently, the same holds true for $X_{k,l}$ and there is no condition on this entry.

{\it Subcase $l>j_{w\inv(k)}$:} This corresponds to the second condition in \eqref{eq:elementsofUMJplus}.
 Write $X''=wX'$ such that $X=LX''.$ Denote by $v=X_{\firstn{d}\times \{l\}}$ and $v''=X''_{\firstn{d}\times \{l\}}$ the $l$-th column of $X$ and $X'',$ respectively. Factorize $L=L_{(d)}L_{(d-1)}\cdots L_{(1)}$ where $L_{(m)}\in \opp\bU$ only has non-trivial entries in row $m.$ Write $L_{\underline{m}}=L_{(m)}L_{(m-1)}\cdots L_{(1)}$
For any column vector $x\in k^{d\times 1}$ and index $m\in \firstn{d},$ denote by $x_{\leq m}$ the vector obtained from $x$ by setting all entries with index $>m$ to zero. Using the lower-triangular shape, one may show that $(Lx)_{\leq m}=L_{\underline{m}}x_{\leq m}$ and $(L\inv x)_{\leq m}=L\inv_{\underline{m}}x_{\leq m}.$
Then we get
$	v_{\leq k}=(Lv'')_{\leq k}=L_{\underline{k}}v''_{\leq k}=L_{\underline{k}}v''_{\leq k-1}$. 
On the other hand, $X''_{k,l}=X'_{w\inv (k),l}=0$, since $l>j_{w\inv(k)}$ and $X'\in M_J^+$. This implies $v''_{\leq k}=v''_{\leq k-1}$ and we can rewrite 
	$v_{\leq k}=L_{\underline{k}}(L\inv v)_{\leq k-1}=L_{\underline{k}}L_{\underline{k-1}}\inv v_{\leq k-1}=L_{(k)}v_{\leq k-1}$ and condition \eqref{eq:thedeterminedentries} follows. (In case there is a non-trivial factor of $\bT$ one simply takes $L\in \bT\opp \bU$ analogously.)
\end{proof}
\section{Open Richardson varieties: Decompositions}\label{THEguys}
We introduce two decompositions of open Richardson varieties of Grassmannians. 
\subsection{Richardson varieties} 
Fix $I=\{i_1<\dots<i_d\}$, $J=\{j_1<\dots<j_d\} \in \setbinom{\firstn{n}}{d}$. The associated {\it open Richardson variety} is defined as
$$\rich{I}{J}\defi Z_I^-\cap Z_J^+=Z_I^-\times_Z Z_J^+\subset Z.$$

Then $\rich{I}{J}\neq \emptyset$ if and only if $I\leq J.$ In this case $\rich{I}{J}$ is a variety of dimension $\ell(J)-\ell(I).$ It is a single point if $I=J$.

\begin{remark}\label{infty}
	Note that $\rich{I}{J}$ does not change if $n$ is replaced by some integer $n'\geq n$.  We tacitly identify all these spaces and use the notation $I,J\in\setbinom{\infty}{d}$ and $\rich{I}{J}\subset\Gr(d,\infty)$ in case we do not want to specify $n\geq\op{max}\{i_d,j_d\}$. 
\end{remark}
\subsection{Fukaya diagrams}\label{sec:fukayadiagrams} We diagrammatically represent $I,J$ in terms of strands aligned to a $d\times n$ matrix:
\begin{enumerate}
\item[(F-1)] a vertical strand that extends from the bottom to the top of the $j_k$-th column (for each $j_k\in J),$
	\item[(F-2)]  a strand that extends vertically from the bottom of the $i_k$-th column to the entry $(k,i_k)$ and then horizontally to the right of the $k$-th row (for each $i_k\in I$).
\end{enumerate} 
We colour the first type of strands in black and the second type in red or blue, and call the resulting diagram the \emph{Fukaya diagram} of shape $(I,J)$. We read these diagrams from the bottom (indicating $I$) to the top (indicating $J$). 

For example, to the sets $I=\{1,2,4\}$ and $J=\{2,5,7\}$ for $n=7$, we associate the following diagram (for comparison we display also $M_I^-\cong Z_I^-$ and $M_J^+\cong Z_J^+$)

\noindent\begin{minipage}[t]{4cm}
$\begin{bmatrix}
1&\cdot&*&\cdot&*&*&*\\
0&1&*&\cdot&*&*&*\\
0&0&0&1&*&*&*\\
\end{bmatrix}
$
\end{minipage}
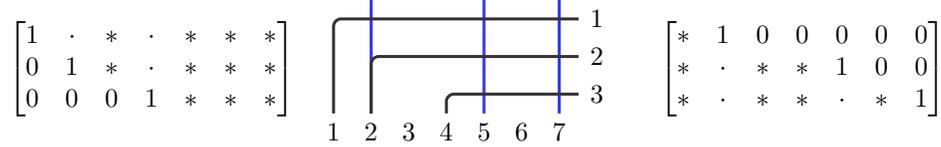
\begin{minipage}[t]{4.5cm}
\begin{tikzpicture}[baseline={(0, 0)}]
	\def\n{7}
	\def\d{3}
	\pgfmathtruncatemacro{\dpp}{\d+1}
	\pgfmathtruncatemacro{\npp}{\n+1}

	\let\mymatrixcontent\empty
		\foreach \row in {1,...,\dpp}{
			\foreach \col in {0,...,\n}{%
				\xappto\mymatrixcontent{ \expandonce{\&}}
	}%
			\gappto\mymatrixcontent{\\}
	}
	\matrix (m) [matrix of nodes, nodes in empty cells, ampersand replacement=\&,
				matrixstyle ] (m){
	\mymatrixcontent 
	};

	\foreach \x in {1,...,\n}
	   \node at (m-\dpp-\x.center) {$\x$};
	\foreach \y in {1,...,\d}
	   \node at (m-\y-\npp.center) {$\y$};
	
	\foreach \x in {2,5,7}
	  \draw[verticallinesopp] (m-1-\x.north) -- (m-\d-\x.south);
	
	\foreach \ik [count=\k] in {1,2,4} 
		\draw[horizontallines]  (m-\d-\ik.south) -- (m-\k-\ik.center) -- (m-\k-\n.east);

	

\end{tikzpicture}
\end{minipage}
\begin{minipage}[t]{4cm}
$\begin{bmatrix}
*&1&0&0&0&0&0\\
*&\cdot&*&*&1&0&0\\
*&\cdot&*&*&\cdot&*&1\\
\end{bmatrix}$
\end{minipage}

We will often omit the numbering of the strands, see also \Cref{infty}, but always count from left to right and top to bottom. The leading  entries of $M_I^-$ and $M_J^+$ correspond to the turning points of the black lines, respectively, the intersection points of the $k$th black with the $k$th blue strand.
\subsection{A naive normal form}
The normal form for the Bruhat cell $Z_I^+$ discussed in \Cref{sec:parabruhat} yields a \emph{naive normal form} for $\rich{I}{J}$ via the set of matrices
\begin{equation}
	\label{eq:naiveparameterrichardson}\begin{split}
	\richparanaive{I}{J}&\defi M^-_I\cap \bG M_J^+\\
	&=\setbuild{X\in k^{d\times n}}{X\in M_I^-, \det X_{\firstn{d}\times J}\neq 0,(X_{\firstn{d}\times J})\inv X\in M_J^+}.
\end{split}
\end{equation}
In particular, the map $[_-]$ from \Cref{sec:parabruhat} sending a matrix to its row span, yields an isomorphism $\richparanaive{I}{J}\iso \rich{I}{J}.$ 
There is one particular example where the naive normal formal allows to directly describe the open Richardson variety:
\begin{example}\label{ex:richardsonfordisjointintervals}
	Let $I=\firstn{d}$ and $J=a+d+\firstn{d}$ for $a\geq 0.$ Then there is an isomorphism (which can easily be read off the Fukaya diagram)
	\begin{eqnarray*}
	\richparanaive{I}{J}\iso k^{d\times a}\times \bG, &&X\mapsto (X_{\firstn{d}\times (d+[a])},X_{\firstn{d}\times J}).
	\end{eqnarray*}
\end{example}
In general the parametrization \eqref{eq:naiveparameterrichardson} is of limited use for us, since the equations involved in the definition of $\richparanaive{I}{J}$ are quite unwieldy. Also, note that the definition of $\richparanaive{I}{J}$ involves a choice making it asymmetric---we could instead use $\bG M^-_I\cap M_J^+.$ With help of the following base change we will present more useful parametrizations.

\subsection{Base change}\label{basechange} We will study the open Richardson variety $\rich{I}{J}$ via the \emph{base change map}
\[\begin{tikzcd}[row sep=0]
	{\basechange{I}{J}:\rich{I}{J}} & {\basechangeimage{I}{J}\subset \bG=\GL_d}\\
	{\phantom{\basechange{I}{J}:}W=[X]} & {(X_{\firstn{d}\times I})\inv(X_{\firstn{d}\times J})\phantom{\subset \bG:=\GL_d}}
	\arrow[from=1-1, to=1-2]
	\arrow[maps to, from=2-1, to=2-2]
\end{tikzcd}\]
Here, $X\in k^{d\times n}$ is any matrix representing $W\in\rich{I}{J}$. The value of $\basechange{I}{J}$ does not depend on this choice, and the submatrices $X_{\firstn{d}\times I}$ and $X_{\firstn{d}\times J}$ are invertible, since $W\in Z^-_I$ and $W\in Z^+_J$. Thus we obtain a well-defined map $\basechange{I}{J}:\rich{I}{J}\to \bG$.  Let $\basechangeimage{I}{J}\subset \bG$ be its image. We will prove in \Cref{thm:reductiontobasechangeimage} that the base change map $\basechange{I}{J}:\rich{I}{J}\to \basechangeimage{I}{J}$ is a vector bundle.

\subsection{Decompositions of open Richardson varieties}\label{decsec}
Recall that $\bB,\opp\bB\subset \bG=\GL_d$ denote the  Borel subgroup of upper-triangular matrices and its opposite and $\bW=S_d$ the Weyl group. We will consider two decompositions of $\rich{I}{J}$ in terms of two Bruhat stratifications of $\bG.$ For $w\in \bW$, we define

\begin{minipage}[t]{6cm}
\fcolorbox{blue}{white}{
\parbox{0.3\textwidth}{
\begin{align*}
\richstratopp{I}{J}{w}&=\basechange{I}{J}\inv(\bB w\opp\bB)\subset \rich{I}{J},\\
\basechangeimagestratopp{I}{J}{w}&=\basechange{I}{J}(\richstratopp{I}{J}{w}),
\end{align*}
}}
\end{minipage}
\begin{minipage}[t]{6cm}
\fcolorbox{red}{white}{
\parbox{0.3\textwidth}{
\begin{align*}
\richstrat{I}{J}{w}&=\basechange{I}{J}\inv(\bB w\bB)\subset \rich{I}{J}, \\
\basechangeimagestrat{I}{J}{w}&=\basechange{I}{J}(\richstrat{I}{J}{w}).
\end{align*}
}}
\end{minipage}

\begin{definition}{\rm We have two resulting decompositions 
\begin{eqnarray*}
\rich{I}{J}=\bigcup_{\mathclap{w\in \basechangeweylopp{I}{J}}} \richstratopp{I}{J}{w}&&
\rich{I}{J}=\bigcup_{\mathclap{w\in \basechangeweyl{I}{J}}} \richstrat{I}{J}{w},
\end{eqnarray*}
we call the  \emph{\Gauss} respectively 
\emph{\Deodhar}.}
\end{definition}
\begin{example}\label{easyex}
	 Consider $Z=\Gr(2,3),$ $I=\{1,2\}$ and $J=\{2,3\}.$ Denote by $W_{x,y}\in Z$ the subspace spanned by rows of the matrix $$\begin{bmatrix}
		1 & 0 & x\\
		0 & 1 & y
	\end{bmatrix}.$$
	The \Gausstype\ and \Deodhar\ of $\rich{I}{J}$ are given by
		\begin{align*}
	\richstratopp{I}{J}{e}&=\setbuild{W_{x,y}}{x,y\neq 0}\cong (\Gm)^2,&&\richstrat{I}{J}{e}=\emptyset,\\
		\richstratopp{I}{J}{s}&=\setbuild{W_{x,y}}{x\neq 0,y= 0}
		\cong \Gm,&&
		\richstrat{I}{J}{s}=\setbuild{W_{x,y}}{x\neq 0}\cong \Gm\times\A^1.
		\end{align*}
We see that the non-empty parts are parametrized by $\basechangeweylopp{I}{J}=\{e,s\}$ and $\basechangeweyl{I}{J}=\{s\}$, respectively. They are all of the form  $(\Gm)^\alpha\times \A^\beta$  for some $\alpha,\beta$.
\end{example}
We observe that the \Gauss\ and the \Deodhar\ of $\rich{I}{J}$ are quite different and might have a different number of parts. We will provide in \Cref{thm:richstratoppeverything,thm:richstrateverything}  a labelling set for the non-empty strata $\richstratopp{I}{J}{w}$ and  $\richstrat{I}{J}{w}$ respectively by explicit sets $\bW\supset \basechangeweylopp{I}{J}\supset \basechangeweyl{I}{J}$. Namely, $\basechangeweylopp{I}{J}$ consists of all elements satisfying the {\it monotonicity condition}
\begin{align}
	\tag{$\leq$}\label{eq:monotonicityWIJ} 	
	\basechangeweylopp{I}{J}&= \{w\in\bW\mid i_{w(k)}\leq j_{k} \text{ for all }k\in \firstn{d}\},
	\end{align}
whereas the elements in $\basechangeweyl{I}{J}$ satisfy additionally the {\it equality condition}
\begin{align}
	\tag{$=$}\label{eq:equalityWIJ}
     \basechangeweyl{I}{J}&=\{w\in\basechangeweylopp{I}{J}\mid k=w(k') \text{ if } i_k=j_{k'}\in I\cap J\}.
\end{align}

 In particular, $\basechangeweyl{I}{J}=\basechangeweylopp{I}{J}$ if $I\cap J=\emptyset$ and $\basechangeweyl{I}{J}=\basechangeweylopp{I}{J}=\{e\}$ if $I=J$. 
 Moreover,  $\basechangeweylopp{I}{J}\subset\bW$ is a lower set for the Bruhat order:
\begin{proposition}\label{makedisjoint} Let $I'=I-J$ and $J'=J-I.$
	\begin{enumerate}[(a)]
	\item Let $w\in \basechangeweylopp{I}{J}$ and $x\in \bW$ with $x\leq w.$ Then $x\in \basechangeweylopp{I}{J}.$
		\item The natural restriction map yields an isomorphism 
		\begin{eqnarray*}
		\basechangeweyl{I}{J}&\iso&\basechangeweyl{I'}{J'}=\basechangeweylopp{I'}{J'}.
		\end{eqnarray*}
	\end{enumerate}
\end{proposition}
\begin{proof}
	For the first part, let $w\in \basechangeweyl{I}{J}$ and $a<b$ such that $w(b)<w(a).$ Let $w'=w(a,b).$ Then $w'\in \basechangeweyl{I}{J}$ since
	$i_{w'(c)}=i_{w(c)}\leq j_c \text{ for all }c\neq a,b,$ $i_{w'(a)}=i_{w(b)}<i_{w(a)}\leq j_a$, and $
		i_{w'(b)}=i_{w(a)}\leq j_a < j_b$.
\end{proof}
\subsection{Marked and strongly marked Fukaya diagrams}
In  \Cref{sec:fukayadiagrams} we introduced the Fukaya diagram of shape $(I,J)$ with strands labeled by $I$ and $J.$ We now define matchings on them to present elements in $\basechangeweylopp{I}{J}$ and $ \basechangeweyl{I}{J}$ graphically and use the colors blue and red to help to distinguish the two situations.
\begin{definition}\label{def:matching} {\rm A \emph{matching} on a Fukaya diagram of shape $(I,J)$ is a collection of $d$ circles marking intersection points\footnote{We use the convention that the intersection point of two overlapping strands $i_k=j_{k'}$ is the point $(k,i_k)$ where the strand corresponding to $i_k$ bends.} of the strands corresponding to $I$ and $J$ such that each line contains exactly one circle. The resulting diagram is called a \emph{matched Fukaya diagram.} It is \emph{strongly matched} if any pair of strands which start at the same point are matched. This requires that $k=w(k')$ if  $i_k=j_{k'}\in I\cap J$. }
\end{definition}
\begin{remark}
	To the Fukaya diagram one might associate a bipartite graph with vertices $I\sqcup
	J$ and an edge between $i_k\in I$ and $j_{k'}\in J$ if the strands intersect, that is, if $i_k\leq j_{k'}.$
	Then, a matching of a Fukaya diagram is the same as a matching of the corresponding graph in the sense of graph theory.
\end{remark}
\begin{proposition}[Matched and strongly matched Fukaya diagrams]
	Mapping $w\in \basechangeweylopp{I}{J}$ to the matched Fukaya diagram of shape $(I,J)$ where the intersection points of the strands $i_{w(k)}$ and $j_k$ are marked by a circle, defines a bijection of sets
	\begin{eqnarray*}
	\label{match}
 \basechangeweylopp{I}{J}&\iso&
	\left\{\text{matched Fukaya diagrams of shape }(I,J)\right\}.
	\end{eqnarray*}
	which restricts to a bijection 
	\begin{eqnarray*}
	\label{strongmatch}
 \basechangeweyl{I}{J}&\iso&
	\left\{\text{strongly matched Fukaya diagrams of shape }(I,J)\right\}.
	\end{eqnarray*}
\end{proposition}
\begin{proof}
In the matched Fukaya diagram associated with $w\in \bW$, the intersection points of the strands $i_{w(k)}$ and $j_{k}$ are marked. These strands intersect precisely if $i_{w(k)}\leq j_{k}.$ Hence, matched Fukaya diagrams correspond to elements in $w\in\basechangeweylopp{I}{J}.$ This proves the first statement and then the second is obvious.
\end{proof}
 In the situation of \Cref{easyex}, the element $s$ corresponds to a strongly matched Fukaya diagram, whereas $e$ does not: 

$\basechangeweylopp{I}{J}\iso\left \{ \begin{array}{l} 
\begin{minipage}{6cm}
\hspace{-0.3cm}
\begin{tikzpicture}[baseline={(0, 0)}]
	\def\n{3}
	\def\d{2}
	\pgfmathtruncatemacro{\dpp}{\d+1}
	\pgfmathtruncatemacro{\npp}{\n+1}

	\let\mymatrixcontent\empty
		\foreach \row in {1,...,\dpp}{
			\foreach \col in {0,...,\n}{%
				\xappto\mymatrixcontent{ \expandonce{\&}}
	}%
			\gappto\mymatrixcontent{\\}
	}
	\matrix (m) [matrix of nodes, nodes in empty cells, ampersand replacement=\&,
				matrixstyle ] (m){
	\mymatrixcontent 
	};

	\foreach \x in {1,...,\n}
	   \node at (m-\dpp-\x.center) {$\x$};
	\foreach \y in {1,...,\d}
	   \node at (m-\y-\npp.center) {$\y$};
	
	\foreach \x in {2,3}
	  \draw[verticallinesopp] (m-1-\x.north) -- (m-\d-\x.south);
	
	\foreach \ik [count=\k] in {1,2} 
		\draw[horizontallines]  (m-\d-\ik.south) -- (m-\k-\ik.center) -- (m-\k-\n.east);

	 Add matching 
	 \foreach \x/\y in {1/2,2/3}
    	\node[matchingnodes] at (m-\x-\y.center) {};
%
\end{tikzpicture}
\hspace{-0.5cm},
\begin{tikzpicture}[baseline={(0, 0)}]
	\def\n{3}
	\def\d{2}
	\pgfmathtruncatemacro{\dpp}{\d+1}
	\pgfmathtruncatemacro{\npp}{\n+1}

	\let\mymatrixcontent\empty
		\foreach \row in {1,...,\dpp}{
			\foreach \col in {0,...,\n}{%
				\xappto\mymatrixcontent{ \expandonce{\&}}
	}%
			\gappto\mymatrixcontent{\\}
	}
	\matrix (m) [matrix of nodes, nodes in empty cells, ampersand replacement=\&,
				matrixstyle ] (m){
	\mymatrixcontent 
	};

	\foreach \x in {1,...,\n}
	   \node at (m-\dpp-\x.center) {$\x$};
	\foreach \y in {1,...,\d}
	   \node at (m-\y-\npp.center) {$\y$};
	
	\foreach \x in {2,3}
	  \draw[verticallinesopp] (m-1-\x.north) -- (m-\d-\x.south);
	
	\foreach \ik [count=\k] in {1,2} 
		\draw[horizontallines]  (m-\d-\ik.south) -- (m-\k-\ik.center) -- (m-\k-\n.east);

	 Add matching 
	 \foreach \x/\y in {1/3,2/2}
    	\node[matchingnodes] at (m-\x-\y.center) {};
	
\end{tikzpicture}
\end{minipage}
\end{array}\hspace{-2cm}\right\}, $
$\basechangeweyl{I}{J}\iso\left \{ \begin{array}{l} 
\begin{minipage}{4cm}
\hspace{-0.3cm}
\begin{tikzpicture}[baseline={(0, 0)}]
	\def\n{3}
	\def\d{2}
	\pgfmathtruncatemacro{\dpp}{\d+1}
	\pgfmathtruncatemacro{\npp}{\n+1}

	\let\mymatrixcontent\empty
		\foreach \row in {1,...,\dpp}{
			\foreach \col in {0,...,\n}{%
				\xappto\mymatrixcontent{ \expandonce{\&}}
	}%
			\gappto\mymatrixcontent{\\}
	}
	\matrix (m) [matrix of nodes, nodes in empty cells, ampersand replacement=\&,
				matrixstyle ] (m){
	\mymatrixcontent 
	};

	\foreach \x in {1,...,\n}
	   \node at (m-\dpp-\x.center) {$\x$};
	\foreach \y in {1,...,\d}
	   \node at (m-\y-\npp.center) {$\y$};
	
	\foreach \x in {2,3}
	  \draw[verticallines] (m-1-\x.north) -- (m-\d-\x.south);
	
	\foreach \ik [count=\k] in {1,2} 
		\draw[horizontallines]  (m-\d-\ik.south) -- (m-\k-\ik.center) -- (m-\k-\n.east);

	 Add matching 
	 \foreach \x/\y in {1/3,2/2}
    	\node[matchingnodes] at (m-\x-\y.center) {};
\end{tikzpicture}
\end{minipage}
\end{array}\hspace{-2.3cm}\right\}.$
Note that the identification $\basechangeweyl{I}{J}\cong\basechangeweylopp{I-J}{J-I}=\basechangeweyl{I-J}{J-I}$ from \Cref{makedisjoint} diagrammatically means removing pairs of strands which start at the same point. The result still encodes all previous matchings, but matchings are now automatically strong. This trick will later allow reductions to the case $I\cap J=\emptyset$. By convention, the empty diagram for $d=0$ has a unique matching which is strong corresponding to the unique (also strong) matching in case $I=J$.

\section{\Gauss\ and Fukaya diagrams}\label{sec:oppstratandfukaya}
We now describe the \Gauss\ in more detail. We illustrate the geometry and combinatorics of the strata in terms of matched Fukaya diagrams.
\subsection{Hybrid normal form}\label{sec:hybridnormalformstratopp} Let $w\in \bW.$ The naive normal form $\richparanaive{I}{J}$ of the open Richardson variety $\rich{I}{J}$ yields a \emph{naive normal form} 
$$\richstratoppparanaive{I}{J}{w}\defi M_I^+\cap \bB w\opp\bB M_J^+$$
for the part $\richstratopp{I}{J}{w}$ of the \Gauss.

It will be much more convenient to `put the factor $\bB$ on the left'. Namely, consider the following set, called {\it hybrid normal form},  of  matrices 
\begin{align}\label{hybrid}
\richstratopppara{I}{J}{w}&\defi\bU M_I^-\cap w\bT\opp\bU_w M_J^+.
\end{align}
We will use this now to construct a kind of hybrid normal form for  $\richstratopp{I}{J}{w}$. It  has the great advantage that the equations describing the subsets $\bU M_I^-$ and $w\bT\opp\bU_w M_J^+$ are simply conditions on entries in the matrices to be either one, zero, non-zero or arbitrary. This then directly shows that $\richstratopp{I}{J}{w}$ is a product of copies of $\A^1$ and $\Gm.$ The multiplicities of the factors will be encoded diagrammatically.

\subsection{An illustrative example}\label{sec:illustrativeexampleoppstrat} We consider the example $I=\{1,2,4\}$ and $J=\{2,5,7\}.$ Denote $s=(1,2),t=(2,3)\in \bW=S_3.$
We want to describe $\richstratopp{I}{J}{w}$ for different choices of $w$ via the hybrid normal form $\richstratopppara{I}{J}{w}.$

We start with the case $w=e$.
\Cref{prop:elementsofUMIminus} and~\Cref{prop:elementsofoppUMIminus} yield the following description for the sets of matrices $\bU M_I^-$ and $e\bT\opp \bU_e M_J^+=\bT\opp \bU M_J^+$
\[
\bU M_I^-=\!\!
	\vcenter{\hbox{
	\begin{tikzpicture}
	\def\n{7}
	\def\d{3}
	\pgfmathtruncatemacro{\dpp}{\d}
	\pgfmathtruncatemacro{\npp}{\n}
	\pgfmathtruncatemacro{\nmm}{\n-1}

	\let\mymatrixcontent\empty
		\foreach \row in {1,...,\d}{
			\foreach \col in {1,...,\nmm}{%
				\xappto\mymatrixcontent{ \expandonce{\&}}
	}%
			\gappto\mymatrixcontent{\\}
	}
	\matrix (m) [matrix of nodes, nodes in empty cells, ampersand replacement=\&,
				matrixstyle, left delimiter={[},right delimiter={]}] (m){
	 \mymatrixcontent 
	};

	\foreach \x/\y in {1/1,2/2,3/4}
		{
		\node[onenodes] at (m-\x-\y.center) {$1$};
		\pgfmathtruncatemacro{\ypp}{\y+1}
		\foreach \a in {\ypp,...,\n}
    	{
			\node[aonenodes] at (m-\x-\a.center) {};}
		}

\end{tikzpicture}
	}}\text{and }
	\bT\opp \bU M_J^+=\!\!
	\vcenter{\hbox{
	\begin{tikzpicture}
	\def\n{7}
	\def\d{3}
	\pgfmathtruncatemacro{\dpp}{\d}
	\pgfmathtruncatemacro{\npp}{\n}
	\pgfmathtruncatemacro{\nmm}{\n-1}

	\let\mymatrixcontent\empty
		\foreach \row in {1,...,\d}{
			\foreach \col in {1,...,\nmm}{%
				\xappto\mymatrixcontent{ \expandonce{\&}}
	}%
			\gappto\mymatrixcontent{\\}
	}
	\matrix (m) [matrix of nodes, nodes in empty cells, ampersand replacement=\&,
				matrixstyle, left delimiter={[},right delimiter={]}] (m){
	 \mymatrixcontent 
	};

	\foreach \x/\y in {1/2,2/5,3/7}
		{
		\node[aonenodes] at (m-\x-\y.center) {};
		\node[gmnodes] at (m-\x-\y.center) {};
		\pgfmathtruncatemacro{\ymm}{\y-1}
		\foreach \a in {1,...,\ymm}
    	{
			\node[aonenodes] at (m-\x-\a.center) {};}
		}

\end{tikzpicture}
	}}\!\!,
\]
the latter denotes the set of all matrices that have the given shape and arbitrary entries  
at crosses and non-zero entries at circled crosses. 
\begin{definition}{\rm
We refer to crosses and circled crosses as {\it $\A^1$- resp. $\Gm$-nodes}}.
\end{definition}
It is now easy to describe their intersection. Their hybrid set is just
\begin{eqnarray}\label{cupoftea}
	\richstratopppara{I}{J}{e}=\!\!
	\vcenter{\hbox{
	\begin{tikzpicture}
	\def\n{7}
	\def\d{3}
	\pgfmathtruncatemacro{\dpp}{\d}
	\pgfmathtruncatemacro{\npp}{\n}
	\pgfmathtruncatemacro{\nmm}{\n-1}

	\let\mymatrixcontent\empty
		\foreach \row in {1,...,\d}{
			\foreach \col in {1,...,\nmm}{%
				\xappto\mymatrixcontent{ \expandonce{\&}}
	}%
			\gappto\mymatrixcontent{\\}
	}
	\matrix (m) [matrix of nodes, nodes in empty cells, ampersand replacement=\&,
				matrixstyle, left delimiter={[},right delimiter={]}] (m){
	 \mymatrixcontent 
	};
	\foreach \x in {2,5,7}
	  \draw[verticallinesopp,draw opacity=0.3] (m-1-\x.north) -- (m-\d-\x.south);
	
	\foreach \ik [count=\k] in {1,2,4} 
		\draw[horizontallines,draw opacity=0.3]  (m-\d-\ik.south) -- (m-\k-\ik.center) -- (m-\k-\n.east);
	\foreach \x/\y/\xx/\yy in {1/1/1/2,2/2/2/5,3/4/3/7}
		{
		\node[onenodes] at (m-\x-\y.center) {$1$};
		\node[aonenodes] at (m-\xx-\yy.center) {};
		\node[gmnodes] at (m-\xx-\yy.center) {};
		\pgfmathtruncatemacro{\ypp}{\y+1}
		\pgfmathtruncatemacro{\yymm}{\yy-1}
		\foreach \a in {\ypp,...,\yy}
    	{
			\node[aonenodes] at (m-\x-\a.center) {};}
		}

\end{tikzpicture}
	}}\!\!.
\end{eqnarray}
By counting the $\A^1$- and $\Gm$-nodes  we see that $\richstratopp{I}{J}{e}\cong \richstratopppara{I}{J}{e}\cong (\Gm)^3\times \A^4.$

We overlayed in \Cref{sec:fukayadiagrams} the matrix with the Fukaya diagram of shape $(I,J)$. To the matched Fukaya diagram for $w=e$ 
we assign diagrams with extra data describing the geometry and cohomology with compact support of $\richstratopp{I}{J}{e}$:
\[\begin{tikzpicture}
	\def\n{7}
	\def\d{3}
	\pgfmathtruncatemacro{\dpp}{\d}
	\pgfmathtruncatemacro{\npp}{\n}
	\pgfmathtruncatemacro{\nmm}{\n-1}

	\let\mymatrixcontent\empty
		\foreach \row in {1,...,\d}{
			\foreach \col in {1,...,\nmm}{%
				\xappto\mymatrixcontent{ \expandonce{\&}}
	}%
			\gappto\mymatrixcontent{\\}
	}
	\matrix (m) [matrix of nodes, nodes in empty cells, ampersand replacement=\&,
				matrixstyle] (m){
	 \mymatrixcontent 
	};
	\foreach \x in {2,5,7}
	  \draw[verticallinesopp] (m-1-\x.north) -- (m-\d-\x.south);
	
	\foreach \ik [count=\k] in {1,2,4} 
		\draw[horizontallines]  (m-\d-\ik.south) -- (m-\k-\ik.center) -- (m-\k-\n.east);
	\foreach \x/\y/\xx/\yy in {1/1/1/2,2/2/2/5,3/4/3/7}
		{
		\node[matchingnodes] at (m-\xx-\yy.center) {};
		}
	\end{tikzpicture}\hspace{10pt}
	\begin{tikzpicture}
	\def\n{7}
	\def\d{3}
	\pgfmathtruncatemacro{\dpp}{\d}
	\pgfmathtruncatemacro{\npp}{\n}
	\pgfmathtruncatemacro{\nmm}{\n-1}

	\let\mymatrixcontent\empty
		\foreach \row in {1,...,\d}{
			\foreach \col in {1,...,\nmm}{%
				\xappto\mymatrixcontent{ \expandonce{\&}}
	}%
			\gappto\mymatrixcontent{\\}
	}
	\matrix (m) [matrix of nodes, nodes in empty cells, ampersand replacement=\&,
				matrixstyle] (m){
	 \mymatrixcontent 
	};
	\foreach \x in {2,5,7}
	  \draw[verticallinesopp] (m-1-\x.north) -- (m-\d-\x.south);
	
	\foreach \ik [count=\k] in {1,2,4} 
		\draw[horizontallines]  (m-\d-\ik.south) -- (m-\k-\ik.center) -- (m-\k-\n.east);

	\foreach \x/\y/\xx/\yy in {1/1/1/2,2/2/2/5,3/4/3/7}
	{
		\node[matchingnodes] at (m-\xx-\yy.center) {};
	}
	\foreach \x/\y/\xx/\yy in {1/1/1/2,2/2/2/5,3/4/3/7}
		{
		\node[aonenodes] at (m-\xx-\yy.center) {};
		\node[gmnodes] at (m-\xx-\yy.center) {};
		\pgfmathtruncatemacro{\ypp}{\y+1}
		\pgfmathtruncatemacro{\yymm}{\yy-1}
		\foreach \a in {\ypp,...,\yy}
    	{
			\node[aonenodes] at (m-\x-\a.center) {};}
		}
\end{tikzpicture}\hspace{10pt}
\begin{tikzpicture}
	\def\n{7}
	\def\d{3}
	\pgfmathtruncatemacro{\dpp}{\d}
	\pgfmathtruncatemacro{\npp}{\n}
	\pgfmathtruncatemacro{\nmm}{\n-1}

	\let\mymatrixcontent\empty
		\foreach \row in {1,...,\d}{
			\foreach \col in {1,...,\nmm}{%
				\xappto\mymatrixcontent{ \expandonce{\&}}
	}%
			\gappto\mymatrixcontent{\\}
	}
	\matrix (m) [matrix of nodes, nodes in empty cells, ampersand replacement=\&,
				matrixstyle] (m){
	 \mymatrixcontent 
	};
	\foreach \x in {2,5,7}
	  \draw[verticallinesopp] (m-1-\x.north) -- (m-\d-\x.south);
	
	\foreach \ik [count=\k] in {1,2,4} 
		\draw[horizontallines]  (m-\d-\ik.south) -- (m-\k-\ik.center) -- (m-\k-\n.east);
	\foreach \x/\y/\xx/\yy in {1/1/1/2,2/2/2/5,3/4/3/7}
	{
		\node[matchingnodes] at (m-\xx-\yy.center) {};
	}
	\foreach \x/\y/\xx/\yy in {2/2/2/5,3/4/3/7}
		{
		\pgfmathtruncatemacro{\ypp}{\y+1}
		\pgfmathtruncatemacro{\yymm}{\yy-1}
		\foreach \a in {\ypp,...,\yymm}
    	{
			\node[aonenodes] at (m-\x-\a.center) {};}
		}
		\node[gmnodes] at (m-1-2.center) {};
		\node[aonenodes] at (m-2-5.center) {};
		\node[gmnodes] at (m-3-7.center) {};
\end{tikzpicture}
\]
We call the middle diagram a \emph{\Gausstype\ Fukaya diagram}. It encodes the positions in matrices in $\richstratopppara{I}{J}{w}$ which are in $\A^1$ and $\Gm$ by marking them with $\A^1$-nodes and $\Gm$-nodes, respectively. It hence describe the variety $\richstratopp{I}{J}{w}.$ 
The right-most diagram is called a \emph{cohomological \Gausstype\ Fukaya diagram}. These cohomological diagrams are obtained from \Gausstype\ Fukaya diagram by replacing the $\Gm$-nodes by either a circle or cross.
The result stands for a basis element in the cohomology with compact support $H_c^\bullet(\richstratopp{I}{J}{w})$ which is a tensor product of the cohomology with compact support of the factors of $\Gm$ and $\A^1.$ A node marked by a cross corresponds to the fundamental class in $H^2_c(\A^1)=H^2_c(\Gm)$  and is of bidegree $(2,1)$ following the conventions from \Cref{sec:bigradings}. A node marked by a circle corresponds to the canonical class in $H^1_c(\Gm)$ which is of bidegree $(1,0)$.  

 In case $w\notin\basechangeweylopp{I}{J}$, the matched Fukaya diagram, and thus all three diagrams, does not exist and the set $\richstratopp{I}{J}{w}$ is empty.

Let us consider now $w=st=(1,2,3)$, which exhibits two important features that did not occur in case $w=e$. We have 
\[st\bT\opp \bU_{st} M_J^+=
	\vcenter{\hbox{
		\begin{tikzpicture}
			\def\n{7}
			\def\d{3}
			\pgfmathtruncatemacro{\dpp}{\d}
			\pgfmathtruncatemacro{\npp}{\n}
			\pgfmathtruncatemacro{\nmm}{\n-1}
		
			\let\mymatrixcontent\empty
				\foreach \row in {1,...,\d}{
					\foreach \col in {1,...,\nmm}{%
						\xappto\mymatrixcontent{ \expandonce{\&}}
			}%
					\gappto\mymatrixcontent{\\}
			}
			\matrix (m) [matrix of nodes, nodes in empty cells, ampersand replacement=\&,
						matrixstyle, left delimiter={[},right delimiter={]}] (m){
			 \mymatrixcontent 
			};
		
			\foreach \x/\y in {2/2,1/7,3/5}
				{
				\node[aonenodes] at (m-\x-\y.center) {};
				\node[gmnodes] at (m-\x-\y.center) {};
				\pgfmathtruncatemacro{\ymm}{\y-1}
			\foreach \x/\y in {1/1,1/3,1/4,1/6,2/1,3/1,3/2,3/3,3/4}
				\node[aonenodes] at (m-\x-\y.center) {};}
		\end{tikzpicture}
	}}
.\]
Working with $\opp\bU_w$ instead of $\opp\bU$ removes the $\A^1$-nodes above the $\Gm$-nodes and prevents a redundant description of $\richstratopp{I}{J}{w}.$
By intersecting with $\bU M_I^+$ we obtain
\[\richstratopppara{I}{J}{st}=
	\vcenter{\hbox{
	\begin{tikzpicture}
	\def\n{7}
	\def\d{3}
	\pgfmathtruncatemacro{\dpp}{\d}
	\pgfmathtruncatemacro{\npp}{\n}
	\pgfmathtruncatemacro{\nmm}{\n-1}

	\let\mymatrixcontent\empty
		\foreach \row in {1,...,\d}{
			\foreach \col in {1,...,\nmm}{%
				\xappto\mymatrixcontent{ \expandonce{\&}}
	}%
			\gappto\mymatrixcontent{\\}
	}
	\matrix (m) [matrix of nodes, nodes in empty cells, ampersand replacement=\&,
				matrixstyle, left delimiter={[},right delimiter={]}] (m){
	 \mymatrixcontent 
	};
	\foreach \x in {2,5,7}
	  \draw[verticallinesopp,draw opacity=0.3] (m-1-\x.north) -- (m-\d-\x.south);
	
	\foreach \ik [count=\k] in {1,2,4} 
		\draw[horizontallines,draw opacity=0.3]  (m-\d-\ik.south) -- (m-\k-\ik.center) -- (m-\k-\n.east);
	
	\foreach \x/\y in {1/1,2/2,3/4}
		\node[onenodes] at (m-\x-\y.center) {$1$};
		
	\foreach \x/\y in {1/7,3/5}
	{
		\node[aonenodes] at (m-\x-\y.center) {};
		\node[gmnodes] at (m-\x-\y.center) {};
	}
	\foreach \x/\y in {1/3,1/4,1/6}
				\node[aonenodes] at (m-\x-\y.center) {};
\end{tikzpicture}
	}}.
\]
We see that for each $k$ with $i_{w(k)}=j_k,$ we get one less copy of $\Gm.$ In total $\richstratopp{I}{J}{st}\cong \richstratopppara{I}{J}{st}\cong (\Gm)^2\times \A^3.$

To summarize, the following four \Gausstype\ Fukaya diagrams can be used to explicitly parametrize the open Richardson variety $\rich{I}{J}$.
\begin{center}
	\begin{tabular}{c@{\hspace{.6cm}}c}
		\begin{tikzpicture}
			\def\n{7}
			\def\d{3}
			\pgfmathtruncatemacro{\dpp}{\d}
			\pgfmathtruncatemacro{\npp}{\n}
			\pgfmathtruncatemacro{\nmm}{\n-1}
		
			\let\mymatrixcontent\empty
				\foreach \row in {1,...,\d}{
					\foreach \col in {1,...,\nmm}{%
						\xappto\mymatrixcontent{ \expandonce{\&}}
			}%
					\gappto\mymatrixcontent{\\}
			}
			\matrix (m) [matrix of nodes, nodes in empty cells, ampersand replacement=\&,
						matrixstyle] (m){
			 \mymatrixcontent 
			};
			\foreach \x in {2,5,7}
			  \draw[verticallinesopp] (m-1-\x.north) -- (m-\d-\x.south);
			
			\foreach \ik [count=\k] in {1,2,4} 
				\draw[horizontallines]  (m-\d-\ik.south) -- (m-\k-\ik.center) -- (m-\k-\n.east);
			\foreach \x/\y/\xx/\yy in {1/1/1/2,2/2/2/5,3/4/3/7}
				{
				\node[matchingnodes] at (m-\xx-\yy.center) {};
				\node[aonenodes] at (m-\xx-\yy.center) {};
				\node[gmnodes] at (m-\xx-\yy.center) {};

				\pgfmathtruncatemacro{\ypp}{\y+1}
				\pgfmathtruncatemacro{\yymm}{\yy-1}
				\foreach \a in {\ypp,...,\yy}
				{
					\node[aonenodes] at (m-\x-\a.center) {};}
				}
		
		\end{tikzpicture}
	&
	\begin{tikzpicture}
		\def\n{7}
		\def\d{3}
		\pgfmathtruncatemacro{\dpp}{\d}
		\pgfmathtruncatemacro{\npp}{\n}
		\pgfmathtruncatemacro{\nmm}{\n-1}
	
		\let\mymatrixcontent\empty
			\foreach \row in {1,...,\d}{
				\foreach \col in {1,...,\nmm}{%
					\xappto\mymatrixcontent{ \expandonce{\&}}
		}%
				\gappto\mymatrixcontent{\\}
		}
		\matrix (m) [matrix of nodes, nodes in empty cells, ampersand replacement=\&,
					matrixstyle] (m){
		 \mymatrixcontent 
		};
		\foreach \x in {2,5,7}
		  \draw[verticallinesopp] (m-1-\x.north) -- (m-\d-\x.south);
		
		\foreach \ik [count=\k] in {1,2,4} 
			\draw[horizontallines]  (m-\d-\ik.south) -- (m-\k-\ik.center) -- (m-\k-\n.east);
		
		\foreach \x/\y in {1/5,2/2,3/7}
		{
			\node[matchingnodes] at (m-\x-\y.center) {};
		}
		\foreach \x/\y in {1/5,3/7}
		{
			\node[aonenodes] at (m-\x-\y.center) {};
			\node[gmnodes] at (m-\x-\y.center) {};
		}
		\foreach \x/\y in {1/3,1/4,3/5,3/6}
					\node[aonenodes] at (m-\x-\y.center) {};
		\end{tikzpicture}
	\\
	$\richstratopp{I}{J}{e}\cong (\Gm)^3\times \A^4$ &
	$\richstratopp{I}{J}{s}\cong (\Gm)^2\times \A^4$  \\[.3cm]
	\begin{tikzpicture}
		\def\n{7}
		\def\d{3}
		\pgfmathtruncatemacro{\dpp}{\d}
		\pgfmathtruncatemacro{\npp}{\n}
		\pgfmathtruncatemacro{\nmm}{\n-1}
	
		\let\mymatrixcontent\empty
			\foreach \row in {1,...,\d}{
				\foreach \col in {1,...,\nmm}{%
					\xappto\mymatrixcontent{ \expandonce{\&}}
		}%
				\gappto\mymatrixcontent{\\}
		}
		\matrix (m) [matrix of nodes, nodes in empty cells, ampersand replacement=\&,
					matrixstyle] (m){
		 \mymatrixcontent 
		};
		\foreach \x in {2,5,7}
		  \draw[verticallinesopp] (m-1-\x.north) -- (m-\d-\x.south);
		
		\foreach \ik [count=\k] in {1,2,4} 
			\draw[horizontallines]  (m-\d-\ik.south) -- (m-\k-\ik.center) -- (m-\k-\n.east);
			
		\foreach \x/\y in {1/2,2/7,3/5}
		{
			\node[matchingnodes] at (m-\x-\y.center) {};
			\node[aonenodes] at (m-\x-\y.center) {};
			\node[gmnodes] at (m-\x-\y.center) {};
		}
		\foreach \x/\y in {2/3,2/4,2/6}
					\node[aonenodes] at (m-\x-\y.center) {};
		\end{tikzpicture}
	&
		\begin{tikzpicture}
		\def\n{7}
		\def\d{3}
		\pgfmathtruncatemacro{\dpp}{\d}
		\pgfmathtruncatemacro{\npp}{\n}
		\pgfmathtruncatemacro{\nmm}{\n-1}
	
		\let\mymatrixcontent\empty
			\foreach \row in {1,...,\d}{
				\foreach \col in {1,...,\nmm}{%
					\xappto\mymatrixcontent{ \expandonce{\&}}
		}%
				\gappto\mymatrixcontent{\\}
		}
		\matrix (m) [matrix of nodes, nodes in empty cells, ampersand replacement=\&,
					matrixstyle] (m){
		 \mymatrixcontent 
		};
		\foreach \x in {2,5,7}
		  \draw[verticallinesopp] (m-1-\x.north) -- (m-\d-\x.south);
		
		\foreach \ik [count=\k] in {1,2,4} 
			\draw[horizontallines]  (m-\d-\ik.south) -- (m-\k-\ik.center) -- (m-\k-\n.east);
		
		\foreach \x/\y in {1/7,2/2, 3/5}
		{
			\node[matchingnodes] at (m-\x-\y.center) {};
		}
			
		\foreach \x/\y in {1/7,3/5}
		{
			\node[aonenodes] at (m-\x-\y.center) {};
			\node[gmnodes] at (m-\x-\y.center) {};
		}
		\foreach \x/\y in {1/3,1/4,1/6}
					\node[aonenodes] at (m-\x-\y.center) {};
		\end{tikzpicture}
	\\
	$\richstratopp{I}{J}{t}\cong (\Gm)^3\times \A^3$  &
	$\richstratopp{I}{J}{st}\cong (\Gm)^2\times \A^3$ \\
	\end{tabular}
	\end{center}
\subsection{Formal definitions and results} To formalize the observations from the example, 
let $I=\{i_1<\dots<i_d\}, J=\{j_1<\dots<j_d\}\in \setbinom{\firstn{n}}{d}$ such that $I\leq J.$
\begin{theorem}[\Gauss]\label{thm:richstratoppeverything} Let $w\in \bW.$ Then
	\begin{enumerate}[(a)]
		\item\label{enum:stratoppparaworks} the map $[_-]: \richstratopppara{I}{J}{w}\to \richstratopp{I}{J}{w}$ is an isomorphism,
		\item\label{enum:basechangeweylopparam} $\richstratopp{I}{J}{w}$ is non-empty if and only if $w\in \basechangeweylopp{I}{J}$, and
		\item\label{enum:sizeoppstrat} we have $\richstratopp{I}{J}{w}\cong (\Gm)^\alpha\times \A^\beta$ where $\alpha=\left|\setbuild{k\in \firstn{d}}{i_{w(k)}\neq j_k}\right|$ and $\beta=\ell(J)-\ell(I)-\alpha-\ell_\bW(w)$.
		\end{enumerate}	
\end{theorem}
 \begin{proof} The first part follows from the definitions by similar, but in fact easier,  arguments to the proof of \Cref{thm:richstrateverything}. The explicit parametrization of $\richstratopp{I}{J}{w}$ via matrices in $\richstratopppara{I}{J}{w}\subset k^{d\times n}$ gives directly a parametrization of $\richstratopp{I}{J}{w}\cong (\Gm)^\alpha\times \A^\beta$; the copies of $\Gm$ and $\A^1,$ correspond to entries in the matrices in $\richstratopppara{I}{J}{w}\subset k^{d\times n}.$ The check of the explicit formulas is straight-forward using the equalities $ \beta=\sum_{k\in\firstn{d}}\left|\setbuild{l\in \firstn{n}}{i_k<l<j_{w\inv(k)}, l\neq j_{k'} \text{ for } w\inv(k)<w\inv(k')}\right|=\ell(J)-\ell(I)-\alpha-\ell_\bW(w)$. Thus the theorem holds.
\end{proof}

The isomorphisms from \Cref{thm:richstratoppeverything} defines a normal form for each non-empty $\richstratopp{I}{J}{w}$ (that means $w\in \basechangeweylopp{I}{J}$); we call it the \emph{hybrid normal form} of $\richstratopp{I}{J}{w}$.

We now give the formal definition of a \Gausstype\ Fukaya diagram:
\begin{definition}{\rm 
	To a matched Fukaya diagram corresponding to $w\in \basechangeweylopp{I}{J}$, we associate a \emph{\Gausstype\ Fukaya diagram}. The diagram is constructed in steps:
	\begin{enumerate}[(G-1)]
		\item\label{Gi} add a circled cross (called \emph{$\Gm$-node}) on each marked intersection point on strands starting at different points (so $i_k\neq j_{w\inv (k)}$),
		\item \label{Gii}for each row indexed by $k\in \firstn{d},$ put crosses (called \emph{$\A^1$-nodes}) at all positions strictly between $(k,i_k)$ (the position where the strand corresponding to $i_k$ bends) and $(k,j_{w\inv(k)})$ (the marked intersection point in the matched Fukaya diagram),
		\item \label{Giii}remove each cross that has a marked intersection point anywhere in the column below it.
	\end{enumerate}
	}
\end{definition}
In case $I=J$, we have $\basechangeweylopp{I}{J}=\{e\}$ and $\richstratopp{I}{J}{e}\cong  (\Gm)^0\times \A^0$ is a point. 
\begin{theorem}  For $w\in \basechangeweylopp{I}{J},$ there is an isomorphism of varieties $\richstratopp{I}{J}{w}\cong (\Gm)^\alpha\times \A^\beta$, where $\alpha$ and $\beta$ are the number of $\Gm$-nodes and $\A^1$-nodes in the corresponding \Gausstype\ Fukaya diagram.
\end{theorem}
\begin{proof}
Steps \ref{Gi}, \ref{Giii} in the definition of a \Gausstype\ Fukaya diagram is just a reformulation of the hybrid set \eqref{hybrid}. Step \ref{Gii} removes the redundant  $\A^1$. 
\end{proof}

\section{\Deodhar\ and Fukaya diagrams}\label{sec:stratandfukaya}
In this section, in parallel to \Cref{sec:oppstratandfukaya}, we describe the parts of the \Deodhar\ explicitly.
\subsection{Hybrid normal form} Let $w\in \bW.$ The naive normal form $\richparanaive{I}{J}$ of the open Richardson variety $\rich{I}{J}$ yields a \emph{naive normal form} 
$$\richstratparanaive{I}{J}{w}\defi M_I^+\cap \bB w\bB M_J^+$$
for the part $\richstrat{I}{J}{w}$ of the \Deodhartype\ decomposition.
Similarly to the discussion in \Cref{sec:hybridnormalformstratopp}, we consider the \emph{hybrid normal form}
\begin{eqnarray}\label{hybridset}
\richstratpara{I}{J}{w}\defi\bU M_I^-\cap w\bT\bU_w M_J^+
\end{eqnarray}
for the part or the \Deodhartype\ stratification $\richstrat{I}{J}{w}.$ The equations describing the subsets $\bU M_I^-$ and $w\bT\opp\bU_w M_J^+$ are conditions on entries in the matrix to be either one, zero, non-zero, arbitrary or completely determined by other entries in the matrix. We encode this information diagrammatically.
\subsection{An illustrative example} 
As in~\Cref{sec:illustrativeexampleoppstrat},
we consider $I=\{1,2,4\}$ and $J=\{2,5,7\}.$ Denote $s=(1,2),t=(2,3)\in \bW=S_3.$

We start with the case $w=s=(1,2).$ 
Propositions~\ref{prop:elementsofUMIminus} and~\ref{prop:elementsofoppUMIminus} give the following description for the sets of matrices
\[	\bU M_I^-=\!\!
\vcenter{\hbox{
\begin{tikzpicture}
\def\n{7}
\def\d{3}
\pgfmathtruncatemacro{\dpp}{\d}
\pgfmathtruncatemacro{\npp}{\n}
\pgfmathtruncatemacro{\nmm}{\n-1}

\let\mymatrixcontent\empty
	\foreach \row in {1,...,\d}{
		\foreach \col in {1,...,\nmm}{%
			\xappto\mymatrixcontent{ \expandonce{\&}}
}%
		\gappto\mymatrixcontent{\\}
}
\matrix (m) [matrix of nodes, nodes in empty cells, ampersand replacement=\&,
			matrixstyle, left delimiter={[},right delimiter={]}] (m){
 \mymatrixcontent 
};

\foreach \x/\y in {1/1,2/2,3/4}
	{
	\node[onenodes] at (m-\x-\y.center) {$1$};
	\pgfmathtruncatemacro{\ypp}{\y+1}
	\foreach \a in {\ypp,...,\n}
	{
		\node[aonenodes] at (m-\x-\a.center) {};}
	}

\end{tikzpicture}
}}\text{and }
	s\bT\bU_s M_J^+=\!\!
	\vcenter{\hbox{
		\begin{tikzpicture}
			\def\n{7}
			\def\d{3}
			\pgfmathtruncatemacro{\dpp}{\d}
			\pgfmathtruncatemacro{\npp}{\n}
			\pgfmathtruncatemacro{\nmm}{\n-1}
		
			\let\mymatrixcontent\empty
				\foreach \row in {1,...,\d}{
					\foreach \col in {1,...,\nmm}{%
						\xappto\mymatrixcontent{ \expandonce{\&}}
			}%
					\gappto\mymatrixcontent{\\}
			}
			\matrix (m) [matrix of nodes, nodes in empty cells, ampersand replacement=\&,
						matrixstyle, left delimiter={[},right delimiter={]}] (m){
			 \mymatrixcontent 
			};
		
			\foreach \x/\y in {1/5,2/2,3/7}
				{
				\node[aonenodes] at (m-\x-\y.center) {};
				\node[gmnodes] at (m-\x-\y.center) {};
				}
			\foreach \x/\y in {1/1,1/3,1/4,2/1,2/5,3/1,3/3,3/4,3/6}
				\node[aonenodes] at (m-\x-\y.center) {};
			\foreach \x/\y in {2/3,2/4}
				\node[onenodes] at (m-\x-\y.center) {$\star$};
		\end{tikzpicture}
	}}\!\!.
\]

Here, $\star$ denotes that the corresponding entry is uniquely determined by the other entries in the matrix. Note that matrices in $\bU M_I^-$ have no conditions on the $\star$-entries in $s\bT\bU_s M_J^+.$
This yields the following simple description of the intersection
\[
	\richstratpara{I}{J}{s}=\!\!
	\vcenter{\hbox{
		\begin{tikzpicture}
			\def\n{7}
			\def\d{3}
			\pgfmathtruncatemacro{\dpp}{\d}
			\pgfmathtruncatemacro{\npp}{\n}
			\pgfmathtruncatemacro{\nmm}{\n-1}
		
			\let\mymatrixcontent\empty
				\foreach \row in {1,...,\d}{
					\foreach \col in {1,...,\nmm}{%
						\xappto\mymatrixcontent{ \expandonce{\&}}
			}%
					\gappto\mymatrixcontent{\\}
			}
			\matrix (m) [matrix of nodes, nodes in empty cells, ampersand replacement=\&,
						matrixstyle, left delimiter={[},right delimiter={]}] (m){
			 \mymatrixcontent 
			};
	\foreach \x in {2,5,7}
	\draw[verticallines,draw opacity=0.3] (m-1-\x.north) -- (m-\d-\x.south);
  
  \foreach \ik [count=\k] in {1,2,4} 
	  \draw[horizontallines,draw opacity=0.3]  (m-\d-\ik.south) -- (m-\k-\ik.center) -- (m-\k-\n.east);
			\foreach \x/\y in {1/5,3/7}
				{
				\node[aonenodes] at (m-\x-\y.center) {};
				\node[gmnodes] at (m-\x-\y.center) {};
				}
			\foreach \x/\y in {1/1,2/2,3/4}
				\node[onenodes] at (m-\x-\y.center) {$1$};
			\foreach \x/\y in {1/3,1/4,2/5,3/6}
				\node[aonenodes] at (m-\x-\y.center) {};
			\foreach \x/\y in {2/3,2/4}
				\node[onenodes] at (m-\x-\y.center) {$\star$};
		\end{tikzpicture}
	}}\!\!.
\]

We see that $\richstrat{I}{J}{s}\cong \richstratpara{I}{J}{s}\cong (\Gm)^2\times \A^4$, illustrated by two $\Gm$- and one $\A^1$-node. As in \Cref{sec:illustrativeexampleoppstrat}, we associate three types of decorated Fukaya diagrams to this situation. Namely, we overlay the strongly matched Fukaya diagram associated to $w\in \bW$ with additional information.

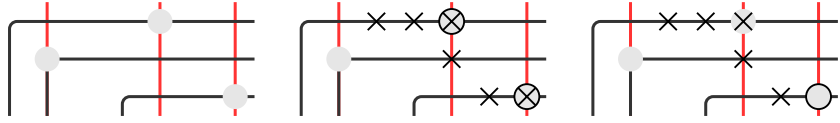
\[
			\begin{tikzpicture}
				\def\n{7}
				\def\d{3}
				\pgfmathtruncatemacro{\dpp}{\d}
				\pgfmathtruncatemacro{\npp}{\n}
				\pgfmathtruncatemacro{\nmm}{\n-1}
		
				\let\mymatrixcontent\empty
				\foreach \row in {1,...,\d}{
					\foreach \col in {1,...,\nmm}{%
						\xappto\mymatrixcontent{ \expandonce{\&}}
				}%
					\gappto\mymatrixcontent{\\}
				}
				\matrix (m) [matrix of nodes, nodes in empty cells, ampersand replacement=\&,
						matrixstyle] (m){
			 	\mymatrixcontent 
				};
				\foreach \x in {2,5,7}
					\draw[verticallines] (m-1-\x.north) -- (m-\d-\x.south);

  				\foreach \ik [count=\k] in {1,2,4} 
	 				 \draw[horizontallines]  (m-\d-\ik.south) -- (m-\k-\ik.center) -- (m-\k-\n.east);
				\foreach \x/\y in {1/5,2/2,3/7}
				{
					\node[matchingnodes] at (m-\x-\y.center) {};
				}
		\end{tikzpicture}\hspace{10pt}
		\begin{tikzpicture}
			\def\n{7}
			\def\d{3}
			\pgfmathtruncatemacro{\dpp}{\d}
			\pgfmathtruncatemacro{\npp}{\n}
			\pgfmathtruncatemacro{\nmm}{\n-1}
	
			\let\mymatrixcontent\empty
			\foreach \row in {1,...,\d}{
				\foreach \col in {1,...,\nmm}{%
					\xappto\mymatrixcontent{ \expandonce{\&}}
			}%
				\gappto\mymatrixcontent{\\}
			}
			\matrix (m) [matrix of nodes, nodes in empty cells, ampersand replacement=\&,
					matrixstyle] (m){
			 \mymatrixcontent 
			};
			\foreach \x in {2,5,7}
				\draw[verticallines] (m-1-\x.north) -- (m-\d-\x.south);

			  \foreach \ik [count=\k] in {1,2,4} 
				  \draw[horizontallines]  (m-\d-\ik.south) -- (m-\k-\ik.center) -- (m-\k-\n.east);
				  \foreach \x/\y in {1/5,2/2,3/7}
				  {
					  \node[matchingnodes] at (m-\x-\y.center) {};
				  }
			\foreach \x/\y in {1/5,3/7}
			{
			\node[aonenodes] at (m-\x-\y.center) {};
			\node[gmnodes] at (m-\x-\y.center) {};
			}
		\foreach \x/\y in {2/5,1/3,1/4,3/6}
			\node[aonenodes] at (m-\x-\y.center) {};
	\end{tikzpicture}\hspace{10pt}
	\begin{tikzpicture}
		\def\n{7}
		\def\d{3}
		\pgfmathtruncatemacro{\dpp}{\d}
		\pgfmathtruncatemacro{\npp}{\n}
		\pgfmathtruncatemacro{\nmm}{\n-1}

		\let\mymatrixcontent\empty
		\foreach \row in {1,...,\d}{
			\foreach \col in {1,...,\nmm}{%
				\xappto\mymatrixcontent{ \expandonce{\&}}
		}%
			\gappto\mymatrixcontent{\\}
		}
		\matrix (m) [matrix of nodes, nodes in empty cells, ampersand replacement=\&,
				matrixstyle] (m){
		 \mymatrixcontent 
		};
		\foreach \x in {2,5,7}
			\draw[verticallines] (m-1-\x.north) -- (m-\d-\x.south);

		  \foreach \ik [count=\k] in {1,2,4} 
			  \draw[horizontallines]  (m-\d-\ik.south) -- (m-\k-\ik.center) -- (m-\k-\n.east);
			  \foreach \x/\y in {1/5,2/2,3/7}
			  {
				  \node[matchingnodes] at (m-\x-\y.center) {};
			  }
		\foreach \x/\y in {3/7}
		{
		\node[gmnodes] at (m-\x-\y.center) {};
		}
		\foreach \x/\y in {1/5}
		{
		\node[aonenodes] at (m-\x-\y.center) {};
		}
	\foreach \x/\y in {2/5,1/3,1/4, 3/6}
		\node[aonenodes] at (m-\x-\y.center) {};
\end{tikzpicture}
\]

The middle diagram is called a \emph{\Deodhartype\ Fukaya diagram}. It encodes the positions in matrices in $\richstratpara{I}{J}{w}$ from $\A^1$ and $\Gm$ by marking them with $\A^1$-nodes and $\Gm$-nodes, respectively. It hence describes the variety $\richstrat{I}{J}{w}.$ This diagram can be constructed from the strongly matched Fukaya diagram by a simple set of rules, which we will discuss below.

The right-most diagram is called a \emph{cohomological \Deodhartype\ Fukaya diagram} and represents a basis element in $H^\bullet_c(\richstrat{I}{J}{w}),$ similar to the cohomological \Gausstype\ Fukaya diagrams.

In our example we obtain the following two \Deodhartype\ Fukaya diagrams, which can be used to explicitly parametrize the open Richardson variety $\rich{I}{J}$:
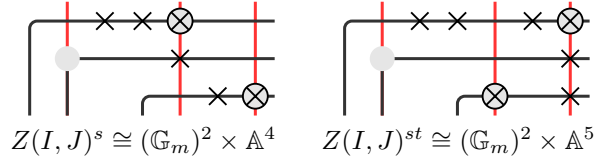
\begin{center}
	\begin{tabular}{c@{\hspace{.6cm}}c}
		\begin{tikzpicture}
			\def\n{7}
			\def\d{3}
			\pgfmathtruncatemacro{\dpp}{\d}
			\pgfmathtruncatemacro{\npp}{\n}
			\pgfmathtruncatemacro{\nmm}{\n-1}
			\let\mymatrixcontent\empty
			\foreach \row in {1,...,\d}{
				\foreach \col in {1,...,\nmm}{%
					\xappto\mymatrixcontent{ \expandonce{\&}}
			}%
				\gappto\mymatrixcontent{\\}
			}
			\matrix (m) [matrix of nodes, nodes in empty cells, ampersand replacement=\&,
					matrixstyle] (m){
			 \mymatrixcontent 
			};
			\foreach \x in {2,5,7}
				\draw[verticallines] (m-1-\x.north) -- (m-\d-\x.south);

			  \foreach \ik [count=\k] in {1,2,4} 
				  \draw[horizontallines]  (m-\d-\ik.south) -- (m-\k-\ik.center) -- (m-\k-\n.east);
				  \foreach \x/\y in {1/5,2/2,3/7}
				  {
					  \node[matchingnodes] at (m-\x-\y.center) {};
				  }
			\foreach \x/\y in {1/5,3/7}
			{
			\node[aonenodes] at (m-\x-\y.center) {};
			\node[gmnodes] at (m-\x-\y.center) {};
			}
		\foreach \x/\y in {1/3,1/4,2/5,3/6}
			\node[aonenodes] at (m-\x-\y.center) {};
	\end{tikzpicture}
	&
	\begin{tikzpicture}
		\def\n{7}
		\def\d{3}
		\pgfmathtruncatemacro{\dpp}{\d}
		\pgfmathtruncatemacro{\npp}{\n}
		\pgfmathtruncatemacro{\nmm}{\n-1}

		\let\mymatrixcontent\empty
		\foreach \row in {1,...,\d}{
			\foreach \col in {1,...,\nmm}{%
				\xappto\mymatrixcontent{ \expandonce{\&}}
		}%
			\gappto\mymatrixcontent{\\}
		}
		\matrix (m) [matrix of nodes, nodes in empty cells, ampersand replacement=\&,
				matrixstyle] (m){
		 \mymatrixcontent 
		};
		\foreach \x in {2,5,7}
			\draw[verticallines] (m-1-\x.north) -- (m-\d-\x.south);

		  \foreach \ik [count=\k] in {1,2,4} 
			  \draw[horizontallines]  (m-\d-\ik.south) -- (m-\k-\ik.center) -- (m-\k-\n.east);
			  \foreach \x/\y in {1/7,2/2,3/5}
			  {
				  \node[matchingnodes] at (m-\x-\y.center) {};
			  }
		\foreach \x/\y in {1/7,3/5}
		{
		\node[aonenodes] at (m-\x-\y.center) {};
		\node[gmnodes] at (m-\x-\y.center) {};
		}
	\foreach \x/\y in {1/3,1/4,1/6,2/7,3/7}
		\node[aonenodes] at (m-\x-\y.center) {};
\end{tikzpicture}
	\\
	$\richstratopp{I}{J}{s}\cong (\Gm)^2\times \A^4$ &
	$\richstratopp{I}{J}{st}\cong (\Gm)^2\times \A^5$\\
	\end{tabular}
	\end{center}
We see that the \Deodhar\ is more economical than the \Gausstype\ decomposition, which consists of four parts.

\begin{remark}
The methods introduced in \Cref{sec:reduction} show that in this example $\rich{I}{J}$ is a vector bundle of rank $3$ over $\GL_2$, and the \Deodhar\ corresponds to the Bruhat stratification of $\GL_2.$
\end{remark}
\subsection{Formal definitions and results} We formalize now our observations. For that let $I=\{i_1<\dots<i_d\}, J=\{j_1<\dots<j_d\}\in \setbinom{\firstn{n}}{d}$ such that $I\leq J.$
\begin{theorem}[\Deodhartype\ decomposition]\label{thm:richstrateverything} Let $w\in \bW.$ Then the following statements hold:
	\begin{enumerate}[(a)]
		\item\label{enum:stratparaworks} The map $[_-]: \richstratpara{I}{J}{w}\to \richstrat{I}{J}{w}$ is an isomorphism.
		\item\label{enum:basechangeweylparam} The part $\richstrat{I}{J}{w}$ is non-empty if and only if $w\in \basechangeweyl{I}{J}.$
		\item\label{enum:sizestrat} If $w\in\basechangeweyl{I}{J}$ then there is an isomorphism $\richstrat{I}{J}{w}\cong (\Gm)^\alpha\times \A^\beta$ where
		\begin{align*}  
			\alpha&=\left|\setbuild{k\in \firstn{d}}{i_{w(k)}\neq j_k}\right|=d-|I\cap J|\text{ and }\\
			\beta&= \ell(J)-\ell(I)-\alpha-|\setbuild{(k,k')}{i_k<j_{k'}<j_{w\inv(k)}}|+\ell_\bW(w).
		\end{align*}
		\item\label{enum:stratification} The \Deodhartype\ decomposition
			$$\rich{I}{J}=\bigcup_{w\in \basechangeweyl{I}{J}}\richstrat{I}{J}{w}$$
			is a stratification, that is, the closure of each stratum is a union of strata.
		\item\label{enum:orderofstrata} The closure relation of the \Deodhar\ comes from the Bruhat order of $\bW$ restricted to $\basechangeweyl{I}{J}.$
	\end{enumerate}
\end{theorem}
\begin{proof}
	For \ref{enum:stratparaworks}, consider the following commutative diagram
	\[\begin{tikzcd}[row sep=0]
		{X=L_1wDL_2Y} && {L_1\inv X=wDL_2Y} \\
		{M_I^-\cap \bB w\bB M_J^+} && {\bU M_I^-\cap  w \bT\bU_w M_J^+} \\
		& {\richstratopp{I}{J}{w}}
		\arrow["{[_-]}"', from=2-1, to=3-2]
		\arrow["\sim", from=2-1, to=2-3]
		\arrow[maps to, from=1-1, to=1-3]
		\arrow["{[_-]}", from=2-3, to=3-2]
		\arrow["\sim", from=2-1, to=3-2]
		\arrow["\sim"', from=2-3, to=3-2]
	\end{tikzcd}\]
	where $L_1\in \bU,$ $D\in \bT,$ $L_2\in \bU_w,$ $X\in M_I^-$ and $Y\in M_J^+.$ The horizontal map is a well-defined isomorphism since $G$ acts freely on $M_J^+$ and multiplication yields an isomorphism $\bU\times w\bT\bU_w\to \bB w\bB.$ Since the left diagonal map in the diagram is an isomorphism, the statement follows by definition, see \eqref{hybridset}. 

	For \ref{enum:basechangeweylparam} we claim that if $X\in \richstratpara{I}{J}{w}\cong \richstrat{I}{J}{w}$ then $w\in \basechangeweyl{I}{J},$ that is, $w$ fulfills \eqref{eq:monotonicityWIJ} and \eqref{eq:equalityWIJ} from \Cref{decsec}. 

	To see this let $k\in \firstn{d}.$ Then $X_{k,l}=0$ for all $l<i_k$ by \eqref{eq:elementsofUMIminus} and $X_{k,j_{w\inv(k)}}\neq 0$ by \eqref{eq:elementsofUMJplus}. Hence $i_k\leq j_{w\inv(k)}$ and thus \eqref{eq:monotonicityWIJ} holds. Assume now $i_k=j_{k'}$ for some $k'\in \firstn{d}.$ By \eqref{eq:elementsofUMIminus}, $X_{k,i_k}=1$ and  $i_k\leq j_{w\inv k}$ by the above discussion.  If \eqref{eq:equalityWIJ} does not hold, that means $k\neq w(k')$, then $i_k< j_{w\inv k}$ and hence $X_{k,i_k}=0$ by the third condition of \eqref{eq:elementsofUMJplus}. This is a contradiction. 
	
The reverse direction in \ref{enum:basechangeweylparam} is clear from  \ref{enum:sizestrat}. 
	For \ref{enum:sizestrat}, let $w\in \basechangeweyl{I}{J}.$
	Note that almost all equations in \eqref{eq:elementsofUMIminus} and \eqref{eq:elementsofUMJplus} (except for \eqref{eq:thedeterminedentries} in \eqref{eq:elementsofUMJplus}) are simply stating that a particular entry of a matrix in $\richstratpara{I}{J}{w}\cong \richstrat{I}{J}{w}$ is either arbitrary, nonzero, equal to zero or equal to one. Since $w\in \basechangeweyl{I}{J},$ \eqref{eq:elementsofUMIminus} puts no conditions on entries which are given by equation \eqref{eq:thedeterminedentries} in \eqref{eq:elementsofUMJplus}. Hence, one may determine $\richstrat{I}{J}{w}$ entry-by-entry using \eqref{eq:elementsofUMIminus} and \eqref{eq:elementsofUMJplus} which amounts to $\richstrat{I}{J}{w}\cong \Gm^\alpha\times \A^\beta$ where $\alpha=d-|I\cap J|$ and, with $\beta'= \sum_{k\in \firstn{d}} (\max(j_{w\inv(k)}-i_{k}-1,0)$, 
		 \begin{align*}
			 \beta&=\beta'-|\setbuild{k'}{i_k<j_{k'}<j_{w\inv(k)}}|
			 +|\setbuild{k'<k}{w\inv(k')>w\inv(k)}|)\\
			 &=\ell(J)-\ell(I)-(d-|I\cap J|)-|\setbuild{(k,k')}{i_k<j_{k'}<j_{w\inv(k)}}|+\ell_\bW(w).
		 \end{align*}
Now, we consider \ref{enum:stratification} and \ref{enum:orderofstrata}. In the special case that $I$ and $J$ are disjoint intervals, $\rich{I}{J}$ is a trivial affine bundle over the $\bG=\GL_d.$ The \Deodhar\ corresponds to the Bruhat stratification of $\bG$ and is hence a stratification. The general case requires reduction arguments which we will develop in \Cref{sec:reduction}. 
\end{proof}

Via the explicit parametrization of $\richstrat{I}{J}{w}$ by matrices in $\richstratpara{I}{J}{w}\subset k^{d\times n}$ we have shown that $\richstrat{I}{J}{w}\cong (\Gm)^\alpha\times \A^\beta$ in \Cref{thm:richstrateverything}. These copies of $\Gm$ and $\A^1,$ correspond to entries in the matrices in $\richstratpara{I}{J}{w}\subset k^{d\times n}.$ We will visualise them via the following decoration of the Fukaya diagram:
\begin{definition}\label{def:multiplicativefukayadiagram}{\rm 
	To a strongly matched Fukaya diagram corresponding to $w\in \basechangeweyl{I}{J}$, we associate a \emph{\Deodhartype\ Fukaya diagram}  constructed as follows:
	\begin{enumerate}[(D-1)]
		\item \label{Di}add a circled cross (called \emph{$\Gm$-node}) on each marked intersection point on strands starting at different points (so $i_k\neq j_{w\inv (k)}$),		\item\label{Dii} for each row indexed by $k\in \firstn{d},$ put crosses (called \emph{$\A^1$-nodes}) at all positions strictly between $(k,i_k)$ (the position where the strand corresponding to $i_k$ bends) and $(k,j_{w\inv(k)})$ (the marked intersection point in the matched Fukaya diagram), but omit those positions which are on red vertical strands,
		\item \label{Diii}put a cross on each intersection point of a vertical line with a horizontal line each containing a marked point, which is strictly above respectively strictly to the left of the intersection point. 
	\end{enumerate}
}	
\end{definition}
\begin{proposition}\label{Deodhartypestrata}
For $w\in \basechangeweyl{I}{J}$ let  $\alpha$ be the number of  $\Gm$-nodes and  $\beta$ be the number of $\A^1$-nodes in the corresponding \Deodhartype\ Fukaya diagram. Then we have an isomorphism of varieties
$\richstrat{I}{J}{w}\cong (\Gm)^\alpha\times \A^\beta$.
\end{proposition}
\begin{proof}
We argue by induction on $\ell(J)$. Assume there exists $l\notin J$, $l+1\in J$ and let $J'$ be obtained from $J$ by swaping $l+1$ with $l$ and vice versa. Similarly define $I'$ from $I$. The following basic cases imply then the statement.

\emph{Case (i): $l\notin I$, $l+1\in I$.} Then one easily verifies  $\richstrat{I}{J}{w}=\richstrat{I'}{J'}{w}$. Diagrammatically this means that two strands starting at $l+1$ get shifted to start at $l$. There are no $\A^1$-nodes on the horizontal part of the strands and the ones on the vertical part get just shifted, thus the $\alpha$ and $\beta$ values do not change. 

\emph{Case (ii): $l,l+1\notin I$ or  $l,l+1\in I$.} Then $\richstrat{I}{J}{w}=\richstrat{I'}{J'}{w}\times \A^1$ and the $\alpha$-value differs by $1$.

\emph{Case (iii): $l\notin I$,  $l+1\in I$.} Locally we have a situation as displayed in the first two diagrams. If we compare this with the situation in the second two diagrams  
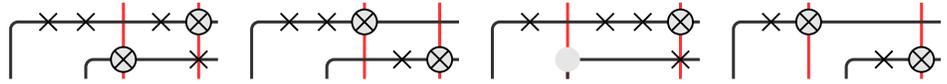
\[
	\begin{tikzpicture}[scale=0.7]
			\def\n{6}
			\def\d{2}
			\pgfmathtruncatemacro{\dpp}{\d}
			\pgfmathtruncatemacro{\npp}{\n}
			\pgfmathtruncatemacro{\nmm}{\n-1}
	
			\let\mymatrixcontent\empty
			\foreach \row in {1,...,\d}{
				\foreach \col in {1,...,\nmm}{%
					\xappto\mymatrixcontent{ \expandonce{\&}}
			}%
				\gappto\mymatrixcontent{\\}
			}
			\matrix (m) [matrix of nodes, nodes in empty cells, ampersand replacement=\&,
					matrixstyle] (m){
			 \mymatrixcontent 
			};
			\foreach \x in {4,6}
				\draw[verticallines] (m-1-\x.north) -- (m-\d-\x.south);

			  \foreach \ik [count=\k] in {1,3} 
				  \draw[horizontallines]  (m-\d-\ik.south) -- (m-\k-\ik.center) -- (m-\k-\n.east);
				  \foreach \x/\y in {1/6,2/4}
				  {
					  \node[matchingnodes] at (m-\x-\y.center) {};
				  }
			\foreach \x/\y in {1/6,2/4}
			{
			\node[aonenodes] at (m-\x-\y.center) {};
			\node[gmnodes] at (m-\x-\y.center) {};
			}
		\foreach \x/\y in {1/2,1/3,2/6, 1/5}
			\node[aonenodes] at (m-\x-\y.center) {};
	\end{tikzpicture}\hspace{5pt}
	\begin{tikzpicture}[scale=0.7]
			\def\n{6}
			\def\d{2}
			\pgfmathtruncatemacro{\dpp}{\d}
			\pgfmathtruncatemacro{\npp}{\n}
			\pgfmathtruncatemacro{\nmm}{\n-1}
	
			\let\mymatrixcontent\empty
			\foreach \row in {1,...,\d}{
				\foreach \col in {1,...,\nmm}{%
					\xappto\mymatrixcontent{ \expandonce{\&}}
			}%
				\gappto\mymatrixcontent{\\}
			}
			\matrix (m) [matrix of nodes, nodes in empty cells, ampersand replacement=\&,
					matrixstyle] (m){
			 \mymatrixcontent 
			};
			\foreach \x in {4,6}
				\draw[verticallines] (m-1-\x.north) -- (m-\d-\x.south);

			  \foreach \ik [count=\k] in {1,3} 
				  \draw[horizontallines]  (m-\d-\ik.south) -- (m-\k-\ik.center) -- (m-\k-\n.east);
				  \foreach \x/\y in {1/4,2/6}
				  {
					  \node[matchingnodes] at (m-\x-\y.center) {};
				  }
			\foreach \x/\y in {1/4,2/6}
			{
			\node[aonenodes] at (m-\x-\y.center) {};
			\node[gmnodes] at (m-\x-\y.center) {};
			}
		\foreach \x/\y in {1/2,1/3,2/5}
			\node[aonenodes] at (m-\x-\y.center) {};
	\end{tikzpicture}\hspace{5pt}
	\begin{tikzpicture}[scale=0.7]
			\def\n{6}
			\def\d{2}
			\pgfmathtruncatemacro{\dpp}{\d}
			\pgfmathtruncatemacro{\npp}{\n}
			\pgfmathtruncatemacro{\nmm}{\n-1}
	
			\let\mymatrixcontent\empty
			\foreach \row in {1,...,\d}{
				\foreach \col in {1,...,\nmm}{%
					\xappto\mymatrixcontent{ \expandonce{\&}}
			}%
				\gappto\mymatrixcontent{\\}
			}
			\matrix (m) [matrix of nodes, nodes in empty cells, ampersand replacement=\&,
					matrixstyle] (m){
			 \mymatrixcontent 
			};
			\foreach \x in {3,6}
				\draw[verticallines] (m-1-\x.north) -- (m-\d-\x.south);

			  \foreach \ik [count=\k] in {1,3} 
				  \draw[horizontallines]  (m-\d-\ik.south) -- (m-\k-\ik.center) -- (m-\k-\n.east);
				  \foreach \x/\y in {1/6,2/3}
				  {
					  \node[matchingnodes] at (m-\x-\y.center) {};
				  }
			\foreach \x/\y in {1/6}
			{
			\node[aonenodes] at (m-\x-\y.center) {};
			\node[gmnodes] at (m-\x-\y.center) {};
			}
		\foreach \x/\y in {1/2,2/6,1/4,1/5}
			\node[aonenodes] at (m-\x-\y.center) {};
	\end{tikzpicture}\hspace{5pt}
	\begin{tikzpicture}[scale=0.7]
			\def\n{6}
			\def\d{2}
			\pgfmathtruncatemacro{\dpp}{\d}
			\pgfmathtruncatemacro{\npp}{\n}
			\pgfmathtruncatemacro{\nmm}{\n-1}
	
			\let\mymatrixcontent\empty
			\foreach \row in {1,...,\d}{
				\foreach \col in {1,...,\nmm}{%
					\xappto\mymatrixcontent{ \expandonce{\&}}
			}%
				\gappto\mymatrixcontent{\\}
			}
			\matrix (m) [matrix of nodes, nodes in empty cells, ampersand replacement=\&,
					matrixstyle] (m){
			 \mymatrixcontent 
			};
			\foreach \x in {3,6}
				\draw[verticallines] (m-1-\x.north) -- (m-\d-\x.south);

			  \foreach \ik [count=\k] in {1,4} 
				  \draw[horizontallines]  (m-\d-\ik.south) -- (m-\k-\ik.center) -- (m-\k-\n.east);
				  \foreach \x/\y in {1/3,2/6}
				  {
					  \node[matchingnodes] at (m-\x-\y.center) {};
				  }
			\foreach \x/\y in {1/3,2/6}
			{
			\node[aonenodes] at (m-\x-\y.center) {};
			\node[gmnodes] at (m-\x-\y.center) {};
			}
		\foreach \x/\y in {1/2,2/5}
			\node[aonenodes] at (m-\x-\y.center) {};
	\end{tikzpicture}
\]
we check easily that the geometry differs by a $\Gm$ respectively an $\A^1$. 
\end{proof}
Since the cohomology (with compact support) of $\Gm$ is two-dimensional and of $\A^1$ is one-dimensional, the following encodes a basis of  $H_c(\richstrat{I}{J}{w})$.
\begin{definition}\label{def:cohomologicalmultiplicativefukayadiagram} {\rm A \emph{cohomological \Deodhartype\ Fukaya diagram} is a \Deodhartype\ Fukaya diagram with the choice of  either a cross or a circle for each $\Gm$-node.}
\end{definition}
Namely, to each \Deodhartype\ Fukaya diagram one can associate $2^\alpha$ many cohomological diagrams. The circles represent the generators of $H^1_c(\Gm)$ while the crosses correspond to the fundamental classes in $H^2_c(\A^1)=H^2_c(\Gm)$.

\subsection{Comparison with Deodhar decomposition}
In this section we compare the combinatorics of the \Deodhar\ and the original Deodhar's decomposition from \cite{deodharGeometricAspectsBruhat1987a} which we recall briefly. 
\begin{definition} {\rm 
Let $x,y\in \pW^P.$ Let $[y]=[t_{\ell(y)},\dots ,t_{1}]\in \pS^{\ell(y)}$ be a reduced expression of $y.$ We denote a subexpression $\gamma$ of $[y]$ as sequence $\gamma=[\gamma_{\ell(y)},\dots ,\gamma_{1}]$ where $\gamma_i\in\{e,t_i\}.$ Moreover, we set $\gamma^i=\gamma_i\cdots\gamma_1.$ A subexpression $\gamma$ is called \emph{distinguished} if for all $i=0,\dots, \ell(y)$ the following two conditions hold
\begin{eqnarray*}
\text{Parabolic condition:}& &\gamma^i\in \pW^P, \\
\text{Jump condition:} &&t_i\gamma^{i-1}<\gamma^{i-1}\implies \gamma_i=t_i.
\end{eqnarray*}
We denote by  $\dist{x}{[y]}\subset\distnot{x}{[y]}$ the set of all distinguished (respectively all) subexpression $\gamma$ of $[y]$ such that $\gamma^{\ell(y)}=x.$ Note that $\dist{y}{[y]}=\{[y]\}.$}
\end{definition}

The behaviour of a distinguished subexpression is encoded by the {\it n-m--data}:
\begin{equation}\label{eq:sizedeodharcells}
	\begin{aligned}
	n_1=n_1(\gamma)&=|\setbuild{i\in \firstn{\ell(y)}}{\gamma^i=\gamma^{i-1}, t_i\gamma^{i-1}\in \bW^P}|\text{ and}\\
	n_2=n_2(\gamma)&=|\setbuild{i\in \firstn{\ell(y)}}{\gamma^i=\gamma^{i-1}, t_i\gamma^{i-1}\not\in \bW^P}|.\\
		m=m(\gamma)&=|\setbuild{i\in \firstn{\ell(y)}}{\gamma^i<\gamma^{i-1}}|,\\
\end{aligned}
\end{equation}
\begin{example}
Let $n=4$ and $d=2$ such that $\pW^P\cong{\setbinom{\firstn{4}}{2}}$. Let $[y]=[s_2,s_1,s_3,s_2]$. There are two subexpression $\gamma$ with $\gamma^{\ell(y)}=s_2$, namely $[s_2,e,e,e,e]$ and $[e,e,e,e,s_2]$.  Both satisfy the parabolic condition; but only the first is distinguished. $\dist{e}{[y]}=\{[e,e,e,e],[s_2,e,e,s_2]\}$ with $(n_1,n_2,m)=(2,2,0)$ and $(n_1,n-2,m)=(2,0,1)$. 
\end{example}

This data is used to label the parts in \emph{Deodhar's decomposition} from \cite{deodharGeometricAspectsBruhat1985}:
\begin{theorem}[Deodhar's decomposition]
There is a decomposition of the open Richardson variety 
\begin{align*}
	\rich{I}{J}&=\bigsqcup_{\gamma\in \dist{x}{[y]}} \richstrat{I}{J}{\gamma}
\end{align*}
into parts $\richstrat{I}{J}{\gamma}$, called \emph{Deodhar cells},  isomorphic to $(\Gm)^{m+n_2(\gamma)}\times \A^{n_1(\gamma)}$.
\end{theorem}
Distinguished subexpression can be computed recursively from the base case $\dist{e}{[]}=\{[\,]\}$ and the following formulas.
\begin{proposition}[D-Recursion] \label{prop:recursiondistexpression} 
	Let $x,y\in \pW^P$ and $[y]=[t_{\ell(y)},\dots ,t_{1}]\in \pS^{\ell(y)}$ be a reduced expression. Let $s=t_{\ell(y)}$ and $[sy]=[t_{\ell(y)-1},\dots ,t_{1}]$ the corresponding reduced expression of $sy<y.$
	\begin{enumerate}[(i)]
		\item \label{Dist1} If $sx<x,$ then $\dist{x}{[y]}=s\dist{sx}{[sy]}.$
		\item \label{Dist2} If $sx\notin \pW^P$ then $\dist{x}{[y]}=e\dist{x}{[sy]}.$
		\item \label{Dist3} If $sx>x,$ then $\dist{x}{[y]}=e\dist{x}{[sy]}\sqcup s\dist{sx}{[sy]}.$
	\end{enumerate}
\end{proposition}
\begin{proof} We obviously always have
$\distnot{x}{[y]}=e\distnot{x}{[sy]}\sqcup s\distnot{sx}{[sy]}.$ 
In case \ref{Dist1}, elements in the first summand fail the jump condition. In case \ref{Dist2}, elements in the second summand fail the parabolic condition.
 For \ref{Dist3}, the second summand only exists if $sx<sy$
  The extra element $s$ or $e$ does not change the {n-m--data} in  \ref{Dist1},  increases $n_2$ by 1 in \ref{Dist1} and finally increases in case \ref{Dist3} $n_1$ by 1 and if $sx<sy$ also $m$ by 1.
\end{proof}

In the same way, the set $\basechangeweyl{I}{J}$ can be computed using $\basechangeweyl{I}{I}=\{e\}$ and the following formulas.
\begin{proposition}[W-Recursion]\label{prop:recursionW}  Let $I=\{i_1<\dots<i_d\} and  J=\{j_1<\dots<j_d\}\in \setbinom{\firstn{n}}{d}.$ Let $s=(l,l+1)\in \pS$ such that $s(J)<J,$ that is, $l+1\in J$ and $l\notin J.$
	\begin{enumerate}[(i)]
		\item If $s(I)<I$ then $\basechangeweyl{I}{J}=\basechangeweyl{s(I)}{s(J)}.$
		\item  If $s(I)=I$ then $\basechangeweyl{I}{J}\iso \basechangeweyl{I}{s(J)}$ where $\sigma\mapsto(\sigma(k),\sigma(k)-1)\sigma.$
		\item If $s(I)>I$ then $\basechangeweyl{I}{J}=\basechangeweyl{I}{s(J)}\sqcup \basechangeweyl{s(I)}{s(J)}.$
	\end{enumerate}
\end{proposition}
Let $x,y$ correspond to $I,J$ under the bijection \eqref{identifications} from \Cref{fixedpoints}.
By comparing the recursive formulas, we obtain directly the following.
\begin{lemma}[Comparison Lemma]\label{prop:comparisoncombinatoricsdeodharstyledeodhar}
	There is a bijection $\dist{x}{[y]}\iso \basechangeweyl{I}{J}$ (in fact unique), which is compatible with the formulas in \Cref{prop:recursiondistexpression,prop:recursionW}.
\end{lemma}
\begin{example} Let $x=e$ such that $I=\firstn{d}.$ Then $[e,\dots,e]\in \dist{x}{[y]}$ corresponds to $\sigma_{0}\in \basechangeweyl{I}{J},$ where
	$\sigma_{0}(k)=\max\setbuild{k'\in \firstn{d}-\sigma_0([k+1,d])}{i_{k'}\leq j_k}$.
\end{example}

By comparing the n-m--data (see the proof of  \Cref{prop:recursiondistexpression}) with the 
description of the strata of the \Deodhar\ (see \Cref{Deodhartypestrata}), we directly obtain the following \emph{combinatorial comparison result}:
\begin{theorem}\label{LA}
	There is an isomorphism $\richstrat{I}{J}{\sigma}\cong\richstrat{I}{J}{\gamma}$ where $\sigma$ corresponds to $\gamma$ with respect to the bijection from \Cref{prop:comparisoncombinatoricsdeodharstyledeodhar}.
\end{theorem}
\begin{remark} We view the \Deodhar\ as a concrete incarnation of Deodhar's decomposition, since it shares all felicitous properties of the latter needed for our purposes. There are, apart from Deodhar's original construction, other combinatorial labellings of the parts in \Cref{LA} as well as parametrizations of the Deodhar's parts, see e.g. \cite{kodamaKPSolitonsTotal2014}, \cite{talaskaNetworkParametrizationsGrassmannian2013}, and its closure relation, \cite{marcottCombinatoricsDeodharDecomposition2020}. Our approach has the advantage that it arises from the Bruhat decomposition of $\bG$ with the closure relation given naturally by the Bruhat ordering.  
\end{remark}
\section{Reduction techniques} \label{sec:reduction}
We now discuss various reduction techniques to describe $\rich{I}{J}$ for general $I\leq J$ in terms of the special case where $I$ and $J$ are disjoint intervals, in which case $\rich{I}{J}$ is simply a trivial bundle over $\bG$, see \Cref{ex:richardsonfordisjointintervals}. This is used to finish the proof  that the \Deodhar\ is a stratification (\Cref{thm:richstrateverything}). Moreover, we will prove that $\rich{I}{J}$ is a vector bundle over $\basechangeimage{I}{J}.$

\begin{lemma} \label{prop:affinebundlesvialinearequations}
	Denote by $V_{r,b}\subset k^{a\times b}$ the subset of rank $r$ matrices. Let $Q_{r,b}=\setbuild{(x,v)}{x\in X, v\in \ker x}.$ Then $\pi:Q_{r,b}\to V_{r,b}$ is a vector bundle of rank $b-r.$
\end{lemma}
\begin{proof}
	Denote by $\pazocal{Q}\to \Gr(r,b)$ the quotient bundle which has the fiber $k^b/W$ over $W\in \Gr(r,b).$
	The map $\pi$ is the pullback of $\pazocal{Q}$ along $V_{r,b}\to \Gr(r,b), X\mapsto [X].$
\end{proof}

\subsection{Reduction to $\basechangeimage{I}{J}$}\label{sec:proofofreductiontobasechangeimage}
\begin{theorem}\label{thm:reductiontobasechangeimage}
	The map $\basechange{I}{J}:\rich{I}{J}\to \basechangeimage{I}{J}$ is a  vector bundle of rank 
		\begin{align*}
			r=|\setbuild{(k,l)\in \firstn{d}\times\firstn{n}}{l\not\in I\cup J, i_k<l<j_k}|.
		\end{align*}
\end{theorem}
\begin{proof}
	We first describe $\basechange{I}{J}$ via the parametrization \eqref{eq:naiveparameterrichardson}
\[\begin{tikzcd}[row sep=0pt]
	{\richparanaive{I}{J}} & {\rich{I}{J}} & {\basechangeimage{I}{J}} \\
	X && {X_{\firstn{d}\times J}}.
	\arrow["{\basechange{I}{J}}", from=1-2, to=1-3]
	\arrow["\sim", from=1-1, to=1-2]
	\arrow[maps to, from=2-1, to=2-3]
\end{tikzcd}\]
For $A\in \basechangeimage{I}{J}$ the preimage $\basechange{I}{J}\inv(A)$ is parametrized by all matrices $X\in k^{d\times n}$
with $X_{\firstn{d}\times J}=A$ such that  $X\in M^+_I$ and $A\inv X\in M^+_J.$ 
For each row in $X$ indexed by $l\in \firstn{n}-(I\cup J)$ these conditions amount to the solutions of the linear equation $(A\inv)_{F_l\times E_l}X_{E_l\times\{l\}}=0$ where 
$E_l=\setbuild{k}{i_k<l}$ and $F_l=\setbuild{k}{j_k<l}.$ 

{\it Claim:} The matrix $(A\inv)_{F_l\times E_l}$ has full rank $|F_l|$,  that is nullity $|E_l|-|F_l|.$ 

This shows that $\basechange{I}{J}$ is isomorphic to the pullback of the product of the vector bundles $Q_{|F_l|\times |E_l|}\to V_{|F_l|\times |E_l|}$ along the maps $A\mapsto (A\inv)_{F_l\times E_l}$ for $l\in \firstn{n}-(I\cup J),$ see \Cref{prop:affinebundlesvialinearequations}. 
Hence, $\basechange{I}{J}$ is a vector bundle of rank 
\begin{align*}
	r&=\sum_{\substack{l\in \firstn{n} \\ l \notin I\cup J}}(|E_l|-|F_l|)=\sum_{\substack{l\in \firstn{n} \\ l \notin I\cup J}}|\setbuild{k\in\firstn{d}}{i_k<l<j_k}|\\
	&=|\setbuild{(k,l)\in \firstn{d}\times\firstn{n}}{l\not\in I\cup J, i_k<l<j_k}|.
\end{align*}
To prove the claim let $w\in \basechangeweyl{I}{J}$ such that $A\in \basechangeimagestrat{I}{J}{w}$ and write $A=L_1wDL_2$ as in the proof of \Cref{thm:richstrateverything}~\Cref{enum:stratparaworks}. Then the condition that $X\in \richparanaive{I}{J}$ is equivalent to $X'=L_1\inv X\in \richstratpara{I}{J}{w}.$ Hence the nullity of $(A\inv)_{F_l\times E_l}$ corresponds to the number $r_l$ of entries $X'_{k,l}$ in $X'$ on which equations \eqref{eq:elementsofUMIminus} and \eqref{eq:elementsofUMJplus} put no condition. This number is given by
\begin{align*}
	r_l&=|\setbuild{k\in\firstn{d}}{l<j_{w\inv(k)}}|-|\setbuild{k\in\firstn{d}}{l\geq i_k}|\\
	&=|\setbuild{k\in\firstn{d}}{l<j_{k}}|-|\setbuild{k\in\firstn{d}}{l\geq i_k}|=|E_l|-|F_l|
\end{align*}
as claimed.
\end{proof}

\subsection{Removing intersections}
Let $i_k=j_{k'}\in I\cap J$ and set $I'=I-{i_k}$ and $J'=J-\{j_{k'}\}.$ The goal is to compare $\rich{I}{J}$ and $\rich{I'}{J'}.$
\begin{remark}
	Note that $\rich{I'}{J'}$ lives inside $\Gr(d-1,\firstn{n}\backslash \{i_k\})$, the Grassmanian  of $d-1$-dimensional subspaces in a vector space with basis indexed by $\firstn{n}\backslash\{i_k\}.$ 
\end{remark}
\begin{theorem}\label{thm:reductiontodisjoint}
There is a vector bundle $\rho: \rich{I}{J}\twoheadrightarrow \rich{I'}{J'}$ of rank $k'-k$ which restrict to $\richstrat{I}{J}{w}\to\richstrat{I'}{J'}{w'},$ where $w'$ is the obvious restriction of $w.$ 
\end{theorem}
\begin{proof} We define $\rho(W)=W'$ where $W=[X]$ and $W'=[X']$ for $X\in \richparanaive{I}{J}$ and $X'=X_{(\firstn{n}-\{i_k\})\times(\firstn{d}-\{k\})}.$
It is straightforward to check that $X'\in \richparanaive{I'}{J'}$ and hence $\rho$ is well-defined. 
Let $A=X_{\firstn{d}\times J}.$ Then the matrix $A$ and its inverse are of the shape
	\begin{align}\label{eq:theshapeofyou}
		A=\left[\begin{array}{c|c|c}
				&	\smash{\vdots}	&	\\
		\smash{\raisebox{.5\normalbaselineskip}{$A_1$}}
        		& 0 				&   \smash{\raisebox{.5\normalbaselineskip}{$A_2$}} \\
      	\hline \\[-\normalbaselineskip]
        \cdots 0& 1 				&  v \\ \hline
                &   0   			& \\
	   \smash{\raisebox{.5\normalbaselineskip}{$0$}}       
	   			& \smash{\vdots}     & \smash{\raisebox{.5\normalbaselineskip}{$A_3$}}
    \end{array}\right] \text{ and }
		A\inv=\left[\begin{array}{c|c|c}
				&	\smash{\vdots}	&	\\
		\smash{\raisebox{.5\normalbaselineskip}{$A_1'$}}
        		& 0 				&   \smash{\raisebox{.5\normalbaselineskip}{$A_2'$}} \\
      	\hline \\[-\normalbaselineskip]
        v'		& 1 				&  0\cdots\\ \hline
                &   0   			& \\
	   \smash{\raisebox{.5\normalbaselineskip}{$0$}}       
	   			& \smash{\vdots}     & \smash{\raisebox{.5\normalbaselineskip}{$A_3'$}}
    \end{array}\right]
	\end{align}
	where the entry $1$ is in position $(k,k')$ and $(k',k),$ respectively. 
	
	
	The $k'$-th row in $A\inv$ is a linear function of $v$ with coefficients depending polynomially on $A_1,A_2,A_3,$ while the other entries of $A\inv$ do not depend on $v.$ Moreover, the requirement that $A\inv X\in M_J^+$ uniquely determines all other entries in the $k$-th row of $X$ as a function of $v'$ and the other entries of $X.$

	Hence, the fiber of  $\rho'$ containing $X$ is parametrized by the possible choices of $v$, which we want to describe now. By removing all corresponding rows and columns from $A$ that correspond to $I'\cap J'$ we may assume that $|I\cap J|=1.$ 
	This ensures that the assumption $X\in M_I^-$ does not constrain $v$ and, likewise, $A\inv X\in M_J^+$ does not constrain $v'.$
	However, there are constraints on $v$ coming from $A\inv.$ Namely, The last $d-k$ entries in the $k'$th-row of $A\inv$ are zero. This imposes a system of $d-k$ linear equations (say given by a matrix $B\in k^{d-k,d-k'}$) on $v.$ We claim that $B$ has full rank, which then shows that $\rho'$ is the pullback of the vector bundle $Q_{d-k,d-k'}\to V_{d-k,d-k'}$ discussed in \Cref{prop:affinebundlesvialinearequations}.

	To see that $B$ has full rank or, equivalently, nullity $k-k'.$ For this, it is more convenient to work with a different parametrization of $\rich{I}{J}.$ Write $A=L_1wL_2$ for $w\in \basechangeweyl{I}{J}$ and consider $X'=L_1\inv X\in \richstratpara{I}{J}{w}.$ 
	We note that $L_1$ does not depend on $v$ and that under this reparametrization, the map $\rho'$ corresponds to the map $\rho'':\richstratpara{I}{J}{w}\to \richstratpara{I'}{J'}{w'}$ also deleting the $i_k$-th column and $k$-th row. Here, $w'\in \basechangeweyl{I'}{J'}$ denotes the obvious restriction of $w.$ To determine the fiber of $\rho''$ amounts to counting the number of entries in the $k$-th row of $X'$ on which \eqref{eq:elementsofUMIminus} and \eqref{eq:elementsofUMJplus} impose no restriction. Using that $w(k')=k,$ this number $r$ equals
	\begin{equation*}
		r=|\setbuild{k''\in \firstn{d}}{k''<k, i_k<j_{w\inv(k'')}}|=k-1-|\setbuild{k''\in \firstn{d}}{k''<k, j_{w\inv(k'')}<i_k}|.
	\end{equation*}
	Now property \eqref{eq:monotonicityWIJ} from \Cref{decsec} and the fact that $w$ is a bijection gives us
	\begin{align*}
		\setbuild{k''\in \firstn{d}}{k''<k, j_{w\inv(k'')}<i_k}&=\setbuild{k''\in \firstn{d}}{k''<k', j_{w\inv(k'')}<i_k}\\
		&=\setbuild{k''\in \firstn{d}}{k''<k'}.
	\end{align*}
	Thus $r=k-1-(k'-1)=k-k'$, and the statement follows.
\end{proof} 
\subsection{Shifting} Let $s\in \pS$ such that $I'=s(I)<I$ and $s(J)=J.$ In other words, $s=(i_k-1,i_k)$ for $i_k\in I-J$ and $i_k-1\not\in I\cup J.$
\begin{theorem}[Shifting Isomorphism]\label{thm:reductionshifting}
	There is an isomorphism of varieties $\rho: \rich{I}{J}\times \A^1\to \rich{I'}{J}$ inducing the identity map  $\basechangeimage{I}{J}\to \basechangeimage{I'}{J}.$
\end{theorem}
\begin{proof} Define $\rho(([X],a))=[X']$ where $X\in\richparanaive{I}{J}$ and $X'$ is defined via
	\begin{align*}
		X'_{[d]\times l}&=
		\begin{cases}
		X_{[d]\times l}& \text{if} l\not\in\{i_k-1,i_k\},\\
		X_{[d]\times i_k}=e_k &\text{if }l=i_k-1,\\
		X_{[d]\times (i_k-1)}+ae_k&\text{if }l=i_k.
		\end{cases}
	\end{align*}
	We leave it to the reader to check that this fulfills the claimed properties.
\end{proof}
\begin{remark}
	There is an analogous result for $s(I)=I$ and $J<s(J).$
\end{remark}
\subsection{Crossing} Let $s\in \pS$ with $I'=s(I)<I$ and $J<J'=s(J).$ In other words, 
 $s=(j_{k},j_{k}+1)$ for $j_k\in J-I$ and $i_{k'}=j_{k}+1\in I-J.$
\begin{theorem}\label{thm:reductioncrossing}
	There is a map $\rho: \rich{I}{J}\times\A^2 \hookrightarrow \rich{I'}{J'}$ which yields isomorphisms $\richstrat{I}{J}{w}\times \A^2\iso \richstrat{I'}{J'}{w}$ and $\basechangeimagestrat{I}{J}{w}\iso \basechangeimagestrat{I'}{J'}{w}$ for $w\in \basechangeweyl{I}{J}.$
\end{theorem}
\begin{proof}
We define $\rho: \rich{I}{J}\times\A^2 \to \rich{I'}{J'}, (W, a, b) \mapsto W'.$
Let $W=[X]$ for $X\in \richparanaive{I}{J}$ and define $W'=[X']$ for $X'\in k^{d \times n}$ via
\begin{align*}
	X'_{[d]\times l}&=
	\begin{cases}
	X_{[d]\times l} &\text{if }l\not\in\{j_k,j_{k}+1\},\\
	a X_{[d]\times j_k}+ e_{k'}&\text{if }l=j_{k}, \\
	X_{[d]\times j_k}+b e_{k'}&\text{if }l=j_{k}+1.
	\end{cases}
\end{align*}
The definition of $X'$ and the fact $X\in M_I^+$ imply $W'=[X']\in Z_{I'}^-.$ With $A=X_{[d]\times J}=\basechange{I}{J}(W)$, the element $Y=A\inv X$ lies in $M^+_J$, and $Y'=A\inv X'$ is given by
\begin{align*}
	Y'_{[d]\times l}&=
	\begin{cases}
	Y'_{[d]\times l}&\text{if }l\not\in\{j_k,j_{k}+1\},,\\
	ae_{k} + Y_{[d]\times (j_k+1)} &\text{if }l=j_{k}, \\
	e_{k} + bY_{[d]\times (j_k+1)}&\text{if }l=j_{k}+1.
	\end{cases}
\end{align*}
This show that $W'=[Y']\in Z_{J'}^+$ using $Y\in M^+_J.$ Hence $\rho$ is  well-defined.

Setting $A'=\basechange{I'}{J'}(W')=(X'_{[d]\times I})\inv X'_{[d]\times J}$ we get 
 $A'=L_aAL_b$ for $L_a,L_b\in \bU$ with 
\begin{align*}
	L_a\inv=I_d+ (e_{k'})\tr a X_{[d]\times j_k} \text{ and } L_b=I_d+(e_{k})\tr b Y_{[d]\times (j_k+1)}.
\end{align*}
Hence $\rho$ induces the map $$\basechangeimage{I}{J}\times \A^2 \to \basechangeimage{I'}{J'}, (A,a,b)\mapsto L_aAL_b.$$ This shows that $\rho$ is compatible with the $\bB\times\bB$ stratification of $\bG.$ To show that $\rho$ yields the promised isomorphism, one uses the explicit parametrizations $\richstratpara{I}{J}{w}\to \richstrat{I}{J}{w}$ and $\richstratpara{I'}{J'}{w}\to \richstrat{I'}{J'}{w}$
\end{proof}

\begin{proof}[Missing reduction argument in proof of  \Cref{thm:richstrateverything}]
 By \Cref{thm:reductiontodisjoint}, there is an affine bundle $\rich{I}{J}\to\rich{I'}{J'}$ where $I'=I-J$ and $J'=J-I$ preserving the \Deodhartype\ decomposition. We hence reduced to the case that $I$ and $J$ are disjoint which we assume from now on.
	
	By applying \Cref{thm:reductioncrossing,thm:reductionshifting} inductively, we may shift elements from $I$ (to the left) and $J$ (to the right) to obtain disjoint intervals $I'$ and $J'$ and a diagram 
\[\begin{tikzcd}
	{\rich{I}{J}\times \A^{k}} & {\rich{I'}{J'}} \\
	{\bigcup_{w\in \basechangeweyl{I}{J}}\rich{I}{J}_w\times \A^{k}} & {\bigcup_{w\in \basechangeweyl{I}{J}}\rich{I'}{J'}_w}
	\arrow["\rho", hook, from=1-1, to=1-2]
	\arrow[Rightarrow, no head, from=2-1, to=1-1]
	\arrow["\sim", from=2-1, to=2-2]
	\arrow[hook, from=2-2, to=1-2]
\end{tikzcd}\]
	for some $k\geq 0$. This induces $\rho$, which respects the \Deodhartype\ decomposition. Since we know that the \Deodhar\ for $\rich{I}{J}$ is a stratification. It suffices to show that $\rho$ is locally closed. But $\rho$ is even closed, since $\basechangeweyl{I}{J}$ is a lower subset of $\basechangeweyl{I'}{J'}=\bW$ for the Bruhat order by \Cref{makedisjoint}. 
\end{proof}


\section{A model for the cohomology of (subsets of) $\operatorname{GL}_d$}\label{sec:cohomologygld}
The reduction techniques from \Cref{sec:reduction} show that, up to factors of $\A^1,$ an open Richardson variety $\rich{I}{J}\subset \Gr(d,n)$ corresponds to a union of Bruhat cells $\bB w\bB$ in $\bG=\GL_d.$ We thus reduced the computation of $H_c^\bullet(\rich{I}{J})$ to a computation of the cohomology with compact support of unions of Bruhat cells in $\GL_d$ which we will undertake in this section. 

The cohomology with compact support of $\GL_d$ (and similarly for other split reductive groups) is an exterior algebra generated by the $d$ elementary symmetric polynomials in degrees $2i-1$ for $i=1,\dots, d$, \cite{kramerHomogeneousSpacesTits2002}. This description however does not restrict to unions of Bruhat cells in $\GL_d.$ We will therefore introduce a model $(\model,\modeld)$ whose cohomology yields cohomology with compact support of $\GL_d$  and which restricts to locally closed unions of Bruhat cells. This then allows to compute $H_c^\bullet(\rich{I}{J}),$ see \Cref{sec:modelforrichardsson}.

We start with the case of tori and $\bG=\GL_d$ for $d=2$, since many aspects for general $d$ can be reduced to these cases. We keep the section mostly self-contained, since it may be of independent interest. 

\subsection{Cohomology of tori}
We identify the cohomology of a split torus $\bT$ of rank $d$ with the exterior algebra of its character lattice 
$H(\bT)=\bigwedge X(\bT).$ The generators on $X(\bT)$ are in bidegree $(1,1).$ 
Via Poincar\'e and algebraic duality, this yields
\begin{align*}
 	H^{k,k}(\bT)&= \mywedge^k X(\bT) = H_{2d-k,d-k}^{BM}(\bT)\text{ and }
 	H_{k,k}(\bT)= \mywedge^k Y(\bT) = H_{c}^{2d-k,d-k}(\bT),
\end{align*}
where $Y(\bT)=X(\bT)^\vee$ denotes the cocharacter lattice. These identifications are functorial with respect to group homomorphism $\bT\to \bT'$ of tori. 
If we have an identification $\bT\cong (\Gm)^d,$ we write $H^\bullet(\bT)=\bigwedge x_1,\dots,x_d$ and $H_\bullet(\bT)=\bigwedge y_1,\dots,y_d$ and abbreviate $X\wedge X'=XX'.$

Given a split cocharacter $\epsilon\in Y(T)$ and let $\bT'=\ker(\epsilon).$ Then, we obtain the following open/closed decomposition
$$\bT\cong \Gm\times \bT'\hookrightarrow \A^1\times \bT' \hookleftarrow \{0\}\times\bT'\cong \bT'.$$
The boundary map $\partial$ in the associated localization seqeuence in cohomology with compact support is given by
\begin{equation}\label{eq:gysinintorus}
	\begin{tikzcd}
	{H_c^k(\bT')} && {H_c^{k+1}(\bT)} \\
	{\mywedge^kY(\bT')} && {\mywedge^{k+1}Y(\bT).}
	\arrow["\partial", from=1-1, to=1-3]
	\arrow[Rightarrow, no head, from=1-1, to=2-1]
	\arrow[Rightarrow, no head, from=1-3, to=2-3]
	\arrow["{\epsilon\wedge-}", from=2-1, to=2-3]
\end{tikzcd}
\end{equation}

\subsection{The case $\GL_2$}
We consider $\bG=\GL_2$ and want to compute $H(\bG)$ using the Gysin long exact sequence associated to the Bruhat stratification $\bG=\bB\cup \bB s \bB.$
For this, let $\bW=\{e,s=(1,2)\}$ and  identify, in the obvious way, 
\begin{align}\label{eq:parambruhatgl2}
	(\Gm)^2\iso \bT,\; \bU\times \{e\}\times\bT\to \bB\text{ and }\bU\times \{s\} \times\bT \times \bU\to \bB s\bB.
\end{align}
This yields identifications of the cohomology with compact support as
$$H_c^\bullet(\bB)=\bigwedge y_1,y_2=\dd{e}\bigwedge y_1,y_2\, \text{ and }\, H_c^\bullet(\bB s\bB)=\bigwedge y_1,y_2=\dd{s}\bigwedge y_1,y_2,$$
where in each case the second equality is just notation; we treat $\dd e$ and $\dd s$ as formal symbols. Moreover, let $\alpha^\vee=y_1-y_2$ be the positive coroot.
\begin{proposition}\label{prop:boundarymapgl2} Under the above identifications, the boundary map 
	$$\partial: H_c^\bullet(\bB)\to H_c^\bullet(\bB s\bB)$$
	in the localisation sequence associated to $\bG=\bB\cup \bB s \bB$ is given by $\delta_eX\mapsto \delta_s\alpha^\vee X.$
\end{proposition}
\begin{proof}
	Since $s\bB s\bB$ is an open neighborhood of $\bB$ we may compute $\partial$ from the boundary map $\partial'$ associated to $s\bB s\bB=(s\bB s\bB\cap \bB s\bB)\cup \bB.$ This decomposition is isomorphic to  
	$\A^1\times (\Gm)^2\times \A^1=(\Gm)^3\times \A^1\cup \{0\}\times(\Gm)^2\times \A^1.$ Indeed, with  $u(a)=I+aE_{2,2}$ we can write elements in $A\in s\bB s\bB$ as
	\begin{align}\label{eq:paramsbsb}
		A=su(a)s\operatorname{diag}(z_1,z_2)u(b)\text{ for }(a,z_1,z_2,b)\in \A^1\times \Gm^2\times \A^1.
	\end{align}
	Correspondingly, we write the cohomology with compact support of $s\bB s\bB\cap \bB s\bB\cong (\Gm)^3\times \A^1$ as $\bigwedge y_a,y_{z_1},y_{z_2}.$
	Using \eqref{eq:gysinintorus}, we obtain a commutative diagram

\[\begin{tikzcd}[row sep=normal]
	{H_c^\bullet(\bB)} & {H_c^\bullet(s\bB s\bB\cap \bB s\bB)} & {H_c^\bullet(\bB s\bB)} \\
	{\dd e\bigwedge y_1,y_2} & {\bigwedge y_a,y_{z_1},y_{z_2}} & {\dd s\bigwedge y_1,y_2.}
	\arrow["{j_!}"', from=1-2, to=1-3]
	\arrow["{\partial'}"', from=1-1, to=1-2]
	\arrow[from=2-2, to=2-3]
	\arrow["{y_a\wedge-}", from=2-1, to=2-2]
	\arrow["\partial", curve={height=-12pt}, from=1-1, to=1-3]
	\arrow[Rightarrow, no head, from=1-3, to=2-3]
	\arrow[Rightarrow, no head, from=1-2, to=2-2]
	\arrow[Rightarrow, no head, from=1-1, to=2-1]
\end{tikzcd}\]
To compute the effect of $j_!$ we need to relate the parametrization \eqref{eq:paramsbsb} of  $s\bB s\bB\cap \bB s \bB$ to the chosen parametrization \eqref{eq:parambruhatgl2} of $\bB s \bB.$ For this, one computes
	$$s u(a)s\diag(z_1,z_2)u(b)=u(a\inv)s\diag(az_1,-a\inv z_2)u((az_1)\inv z_2+b).$$
	Thus, $j_!$ induces $y_a\mapsto y_1-y_2=\alpha^\vee,$ $y_{z_1}\mapsto y_1$, $y_{z_2}\mapsto y_2$, and $\partial=\alpha^\vee\wedge-$.
\end{proof}
The cohomology with compact support fits then into the short exact sequence
$$0\to \coker(\partial)\cong\langle \dd{s}, \dd{s}y_1=\dd{s}y_2   \rangle \to H_c^\bullet(\GL_2)\to \ker(\partial)=\langle \dd{e}(y_1-y_2), \dd{e}y_1y_2 \rangle\to 0.$$
\begin{remark}
By Poincar\'e and algebraic duality, this also yields computations for of the (Borel--Moore) homology and cohomology.
\end{remark}
\subsection{The case of adjacent Bruhat cells}\label{sec:adjacentbruhatcells}
Now let $\bG=\GL_d.$ We use the standard isomorphism $(\Gm)^d\cong \bT$ to parametrize the Bruhat cell  for $w\in \bW$ as
$$\bU\times\{w\}\times \bT \times \bU_w\iso \bB w\bB.$$
This yields a natural identification
$H_c^\bullet(\bB w\bB)=\dd w \bigwedge Y(\bT)$
where, again, $\dd w$ is a formal symbol. 

Let $w<w'$ such that $\ell(w')=\ell(w)+1.$ Then there is a reflection $t=(k,l)$ with $k<l$ such that $w'=wt.$ We denote the corresponding positive coroot by $\alpha_t^\vee=y_k-y_l\in Y(\bT).$ We obtain the locally closed subset $\bB wt\bB\cup \bB w\bB\subset \bG$ with an open/closed decomposition generalising \eqref{eq:parambruhatgl2}. We compute the boundary map in the associated localization sequence:
\begin{proposition}\label{prop:generalboundarymap} Under the above identifications, the boundary map 
	$$\partial: H_c^\bullet(\bB w\bB)\to H_c^\bullet(\bB wt\bB)$$
	in the localization long exact sequence associated to $\bB w\bB\cup \bB wt\bB$ is given by
	\begin{align}\label{eq:generalboundarymap}
	\partial(\delta_{w}Y)=\delta_{wt}\alpha^\vee_tY.
	\end{align}
\end{proposition}
\begin{proof} 
	Write $\bG_{k,l}\subseteq \bG$ for the copy of $\GL_2$ acting on the coordinates $k,l.$ Then
	$\bB w\bB\cup \bB wt\bB=\bB w\bG_{k,l}\bB$ such that the decomposition on the left corresponds to the Bruhat decomposition \eqref{eq:parambruhatgl2} of $\bG_{k,l}$ on the right-hand side. We decompose $\bT=\bT_{k,l}\times\bT_{\neq k,l}$ where $\bT_{k,l}\subset \bG_{k,l}$ and $\bT_{\neq k,l}$ is in the centralizer of $\bG_{k,l}.$ 
	For appropriate $\bU_1,\bU_2,\bU_3\subset \bU$ multiplication yields compatible isomorphisms
	
\[\begin{tikzcd}
	{\bU_1 \times w \times \bU_2\bT_{k,l}\times \bT_{\neq k,l} \times \bU_3} & {\bB w\bB} \\
	{\bU_1 \times w \times \bG_{k,l}\times \bT_{\neq k,l} \times \bU_3} & {\bB w\bG_{k,l}\bB} \\
	{\bU_1 \times w \times \bU_2t\bT_{k,l}\bU_2\times \bT_{\neq k,l} \times \bU_3} & {\bB wt\bB}
	\arrow["\sim", from=2-1, to=2-2]
	\arrow[hook', from=1-2, to=2-2]
	\arrow[hook, from=3-2, to=2-2]
	\arrow[hook, from=3-1, to=2-1]
	\arrow[hook', from=1-1, to=2-1]
	\arrow["\sim", from=3-1, to=3-2]
	\arrow["\sim", from=1-1, to=1-2]
\end{tikzcd}\]
This show that $\partial=\partial'\otimes\id$ where $\partial'$ is the boundary map associated to the Bruhat decomposition of $\bG_{k,l}.$ We hence reduced the statement to the $\GL_2$-case. This is proven in \Cref{prop:boundarymapgl2}.
\end{proof}

\subsection{The combinatorial model}\label{sec:combinatorialmodel} We next introduce our combinatorial model for the cohomology with compact support of $\bG=\GL_d.$ 
The \emph{Nil-Coxeter algebra} $\cox,$ see \cite{fominSchubertPolynomialsNilcoxeter1994}, is the algebra spanned by the symbols $\ddd{w}$ for $w\in \bW=S_d$ with the following deformed multiplication (with identity element $\ddd e\in \cox$):
$$\ddd{w}\ddd{w'}=\begin{cases}
	\ddd{ww'}, & \text{if $\ell(w)+\ell(w')=\ell(ww')$, }\\
	0, & \text{otherwise.}\end{cases}$$
\begin{definition}\label{def:exteriornilcoxeter}{\rm
We define the \emph{exterior Nil-Coxeter algebra} as the semidirect product of $\cox$ with the exterior algebra
$$\model=\cox\#\bigwedge y_1,\dots,y_d$$
via the multiplication rule
$Y\ddd{w}=\ddd{w}w\inv(Y)$ for $w\in S_d$ and $Y\in \bigwedge y_1,\dots,y_d$.
}
\end{definition}
The exterior Nil-Coxeter algebra $\model$ is obviously an $(\cox,\bigwedge y_1,\dots,y_d)$-bimodule. Let $\cox_w=\ddd w\bigwedge y_1,\dots,y_d$ and define a $\bigwedge y_1,\dots,y_d$-linear differential $\modeld$ on $\model$ via
\begin{align}\label{eq:differentialinmodel}
    \modeld(\ddd w)=\sum_{\mathclap{\substack{t\in \bR\\ \ell(wt)=\ell(w)-1}}}\ddd{wt}\alpha^\vee_t.
\end{align}
Here $\bR=\setbuild{wsw\inv}{s\in\bS, w\in\bW}$ denotes the set of reflections/transpositions.

We will also use the basis $\dd{w}=\ddd{w_0w}$ for $w_0$ the longest element, for which the differential has the form
\begin{align}\label{eq:differentialinmodeldualbasis}
    \modeld(\dd w Y)=\sum_{\mathclap{\substack{t\in \bR\\ \ell(wt)=\ell(w)+1}}}\dd{wt}\alpha^\vee_tY.
\end{align}
\begin{proposition}
	The map $d$ is indeed a differential, that is, $d^2=0.$
\end{proposition}
\begin{proof}
	Each interval of length $2$ in the Bruhat order is of the form
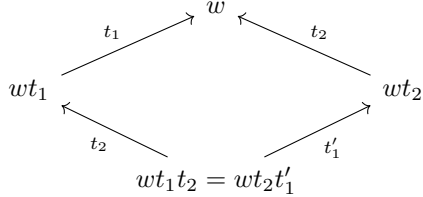
\[\begin{tikzcd}[scale=0.005]
	& {w} \\
	{wt_1} && {wt_2} \\
	& wt_1t_2=wt_2t_1'
	\arrow["{t_2}", from=3-2, to=2-1]
	\arrow["{t_1'}"', from=3-2, to=2-3]
	\arrow["{t_1}", from=2-1, to=1-2]
	\arrow["{t_2}"', from=2-3, to=1-2]
\end{tikzcd}\]
for some $w\in \bW$ such that
$\ell(w)=\ell(wt_1t_2)+2$ and $t_1'=t_2t_1t_2$ for reflections $t_1,t_2\in \bR$, see \cite[Lemma 2.7.3]{bjornerCombinatoricsCoxeterGroups2005}.
It is enough to show that the coefficient $C$ of $\ddd{wt_1t_2}$ in $\modeld(\modeld(\ddd{w}))$. But $C=\alpha^\vee_{t_2}\alpha^\vee_{t_1}+\alpha^\vee_{t_1'}\alpha^\vee_{t_2}=\alpha^\vee_{t_2}\alpha^\vee_{t_1}-\alpha^\vee_{t_2}t_2(\alpha^\vee_{t_1})=\alpha_{t_2}^\vee((1-t_2)(\alpha^\vee_{t_1}))=0$ since the image of $Y(\bT)$ under $(1-t_2)$ is spanned by $\alpha^\vee_{t_2}.$
\end{proof}
\begin{proposition} [dg-algebra structure]\label{prop:dgstructure} The following statements hold.
    \begin{enumerate}
        \item Let $w\in \bW$ and $w=s_1\dots s_{\ell(w)}$ a reduced expression with $s_i\in \bS.$ Then
        \begin{align*}
            \modeld (\ddd w)=\sum_{i=1}^{\ell(w)}\ddd{s_1\cdots s_{i-1}}\alpha_{s_i}^\vee\ddd{s_{i+1}\cdots s_{\ell(w)}}.
        \end{align*}
        \item \label{item:dgstructure:properties} The map $d:\model\to \model$ is uniquely determined by the following properties.
            \begin{enumerate}[(i)]
                \item $d$ is linear with respect to the right action of $\bigwedge y_1,\dots, y_n.$
                \item $d$ fulfills the Leibniz rule $d(ab)=ad(b)+d(a)b.$
                \item $d(\ddd{s})=\alpha_s^\vee$ for simple reflections $s\in \bS.$
            \end{enumerate}
    \end{enumerate}
\end{proposition}
\subsection{Combinatorics and geometry of Bruhat cells}
Recall that $\dd{w}=\ddd{w_0w}$ for $w_0\in \bW$ the longest element. We think of $\dd{w}$ and $\ddd{w}$ as Poincar\'e dual basis elements.
We may identify the algebra $\model$ with the direct sum of the cohomology with compact support (or homology) of the Bruhat cells of $\bG=\GL_d$ via
\begin{equation}\label{A}
\model=\bigoplus_{w\in\bW}\model_w=\bigoplus_{w\in\bW} \delta_w\bigwedge Y(\bT)= \bigoplus_{w\in\bW} H_c^\bullet(\bB w\bB).
\end{equation}
This puts a bigrading on $\model$ inherited from the bigrading on the cohomology with compact support such that for $Y\in \bigwedge^kY(\bT)$ and $w\in \bW$ we have
$$\begin{array}{c|c}
    &\text{bidegree}\\
    \hline
    \dd{w}Y & (2(\ell(w)+\dim \bB)-k,\ell(w)+\dim \bB-k),\\
    \ddd{w}Y & (2(\dim \bG- \ell(w))-k,\dim \bG- \ell(w)-k).
\end{array}$$
We write $\model^\bullet$ and $\model^{\bullet,\bullet}$ to refer to the first and to both gradings, respectively.

\begin{remark}
The multiplication in $\model$ is not strictly compatible with the grading, but one rather has
$$\model^{p,q}\model^{p',q'}\subset \model^{p+p'-2\dim \bG,q+q'-\dim\bG},$$
so that multiplication is a map of bidegree $(2\dim \bG, \dim \bG).$ The reason for this comes from geometry. The multiplication map $m\colon\bG\times \bG\to \bG$ is smooth of codimension $\dim \bG$ and hence induces, via integration along fibers, a map $$H^{p,q}_c(\bG)\otimes H^{p',q'}_c(\bG)\stackrel{m_!}{\longrightarrow}H_c^{p+p'-2\dim \bG,q+q'-\dim\bG}(\bG).$$
To turn $\model$ into a bigraded dg-algebra one may use the Poincar\'e dual grading by changing bidegrees $(p,q)\mapsto (p-2\dim\bG,q-\dim\bG)$ so that $\ddd{w}$ is in degree $(-2\ell(w),- \ell(w)).$
\end{remark}

\subsection{A spectral sequence for $\model=H_c^\bullet(\GL_d)$}

We equip the group $\bG$ with a decreasing filtration $F$ by the open subsets $$F^p\bG=\bigcup_{\ell(w)\geq p}\bB w\bB.$$ 
This filtration induces a filtration on the complex of compactly supported cochains via their support and yields a spectral sequence
\begin{align}\label{eq:spectralsequencegln}
	E_1^{p,n-p}=\bigoplus_{\mathclap{\ell(w)=p}} H^n_c(\bB w\bB)=\bigoplus_{\mathclap{\ell(w)=p}} \model_w^n \Rightarrow H_c^n(\bG).
\end{align}
The boundary maps on the first page of the spectral sequence are given by the differential $\modeld$ of  $\model.$
\begin{theorem}\label{thm:spectralsequencehasmodel} With the above notations, we have
	$(E_1,d_1)=(\model,\modeld).$
\end{theorem}
\begin{proof}
The boundary map $d_1$ on $E_1$ come from the boundary map in the long exact sequences in the cohomology with compact support associated to the open/closed decompositions $F^p\bG-F^{p-2}\bG=(F^p\bG-F^{p-1}\bG)\sqcup( F^{p-1}\bG-F^{p-2}\bG).$
We computed the boundary maps $\partial: H_c^n(\bB w\bB)\to H^{n+1}_c(\bB wt\bB)$ in Section~\ref{sec:adjacentbruhatcells}.
The statement follows by comparing the definition \eqref{eq:differentialinmodeldualbasis} of the differential $\modeld$ on $\model$, and the formula for $\partial,$ see \eqref{eq:generalboundarymap}.
\end{proof}
\subsection{Locally closed unions of Bruhat cells in $\GL_d$}\label{sec:locallyclosedsubsetsgld} The spectral sequence and the combinatorial model discussed above can be generalized to locally closed subsets of $\bG=\GL_d$ which are unions of double cosets for $\bB.$ 

Let $\bI\subseteq \bW$ be a subset that is locally closed (i.e. interval closed) with respect to the Bruhat order, that is, if $w_1,w_3\in \bI$ and $w_1\leq w_2\leq w_3$ then $w_2\in  \bI.$ 
We write $\bI\leq w$ if there is a $w'\in \bI$ with $w'\leq w$ and $\bI<w$ if $w'<w$ for all $w'\in \bI.$

We can then form the subquotient
$$\model_\bI=\model_{\bI\leq }/\model_{\bI<}=\bigoplus_{\bI\leq w} \model_w / \bigoplus_{\bI<w} \model_w$$
which inherits a differential from that of $\model$, \eqref{eq:differentialinmodeldualbasis}. The underlying module is simply
$$\model_\bI=\bigoplus_{w\in\bI} \model_w.$$

With $\bG_\bI=\bigcup_{w\in \bI} \bB w\bB$ we obtain as in \eqref{eq:spectralsequencegln} the spectral sequence
\begin{align}\label{eq:spectralsequencelocallyclosedsubset}
	E_1^{p,n-p}=\bigoplus_{\mathclap{\substack{w\in \bI\\\ell(w)=p}}} H^n_c(\bB w\bB)=\bigoplus_{\mathclap{\substack{w\in \bI\\\ell(w)=p}}} \model_w^n \Rightarrow H_c^n(\bG_\bI).
\end{align}
Similar to Theorem~\ref{thm:spectralsequencehasmodel} we have $(E_1,d_1)=(\model_\bI, \modeld).$
\subsection{Degeneration via weights}\label{sec:degenerationviaweights} The goal of this section is to show the following.
\begin{theorem}\label{thm:spectralsequencecohomologyglddegenerates}
	The spectral sequences \eqref{eq:spectralsequencegln} and \eqref{eq:spectralsequencelocallyclosedsubset} degenerate on page $2.$
\end{theorem}
\begin{proof}
	It suffices to consider \eqref{eq:spectralsequencelocallyclosedsubset}, since the spectral sequence \eqref{eq:spectralsequencegln} is the special case $\bI=\bW$. 
	Let $r\geq 2.$ We show that the differential $d_r:E_r^{p,n-p}\to E_r^{p+r,(n+1)-(p+r)}$ is zero. The domain and codomain are subquotients of 
	$$ A=\bigoplus_{\mathclap{\substack{w\in \bI\\\ell(w)=p}}} H^n_c(\bB w\bB)\, \text{ and }\, B=\bigoplus_{\mathclap{\substack{w\in \bI\\\ell(w)=p+r}}} H^{n+1}_c(\bB w\bB).$$
	The differential $d_r$ preserves weights. So it suffices to show that $A$ and $B$ have no common weights. Recall that $H^{\bullet,\bullet}_c(\bB w\bB)$ is supported in bidegrees $(2(\ell(w)+\dim \bB)-k),\ell(w)+\dim \bB-k)$ for $k=0,\dots,d.$ Hence, $A$ and $B$ are supported in bidegrees
	$(2(\dim \bB+p)-k, \dim \bB+p-k)$ and $(2(\dim \bB+p+r)-k', \dim \bB+p+r-k'),$ respectively. The difference of bidegrees is $(k-k'+2r,k-k'+r)$. On the other hand, the difference of cohomological degrees is $k-k'+2r=(n+1)-n=1,$ and thus the difference in weight, that is  $k-k'+r$, is zero if and only if $r=1.$ 
\end{proof}
\begin{corollary}
	The filtration $F$ on $\bG_\bI$ induces a filtration on $H_c^\bullet(\bG_\bI)$ such that
	$$H^\bullet(\model_\bI,\modeld)=\gr_F^\bullet H^\bullet_c(\bG_\bI).$$
\end{corollary}
\subsection{Degeneration via induction} In this section we sketch a different and combinatorial proof that the spectral sequence for $H^\bullet_c(\bG)$ degenerates. It however requires the knowledge of the dimension of $H^\bullet_c(\bG)$.
\begin{theorem}\label{thm:spectralsequencecohomologyglddegenerates2}
	The spectral sequence \eqref{eq:spectralsequencegln} degenerates on page $2.$
\end{theorem}
The cohomology of $\bG=\GL_d$ can be computed by induction on $d$ using the Serre spectral sequence associated to the fibration
\[\begin{tikzcd}
	{\bG_P\ltimes \A^{d-1}:=\GL_{d-1}\ltimes \A^{d-1}} & {\bG=\GL_d} & {\A^d-\{0\}}
	\arrow[from=1-2, to=1-3]
	\arrow[from=1-1, to=1-2]
\end{tikzcd}\]
mapping a matrix $A\in \GL_d$ to its last column $Ae_d\in \A^d-\{0\}.$ The spectral sequence degenerates on page $2$, see e.g. \cite{kramerHomogeneousSpacesTits2002}.

We construct an algebraic analog of this spectral sequence for $\model.$ For this, let $S_{d-1}\times S_1=\bW_P\subset \bW$ be the Weyl group of $\bG_P$, and let $\bW^P\iso \bW/\bW_P$ be the set of shortest coset representatives.

We first construct a model for the cohomology of $\A^d-\{0\}.$
Let $v_i=(d-i,\dots,d)$  and $t_i=(d-i,d)$ for $i=1,\dots,d-1.$ Then the Bruhat order on $\bW^P$ is the chain
\begin{equation}\label{eq:chainofreflections}\begin{tikzcd}
	{e=v_0} & {v_1} & \dots & {v_{d-1}}
	\arrow["{t_1}", from=1-1, to=1-2]
	\arrow[from=1-2, to=1-3]
	\arrow["{t_{d-1}}", from=1-3, to=1-4].
\end{tikzcd}\end{equation}
The $\delta_{v_i}$ span $\model^P=\cox^P\#\bigwedge y_d\subset\model$ with  differential $\modeld^P(\dd{v_i}Y)=\dd{v_{i+1}}y_dY.$  The cohomology of $\model^P$ splits off and is spanned by $\delta_ey_d$ and $\delta_{v_{d-1}}.$
We identify
 $$\model^P_{v_i}=\delta_{v_i}\bigwedge y_d=H^{\bullet}_c(\A^{i}\times \Gm\times \{0\}^{d-i-1})$$
 so that $\delta_{v_i}$ and $y_d$ have bidegree $(2(i+1), i+1)$ and $(-1,-1),$ respectively. The differential $\modeld^P$ preserves the second degree induced by the weights. We consider the descending weight filtration $W^q\model^P=(\model^P)^{\bullet,\geq q}$ associated to the weight grading.

Let $\model_P=\cox_P\#\bigwedge y_1,\dots, y_{d-1}$ be the algebra corresponding to $\bW_P.$
Multiplication yields an isomorphism of vector spaces
$$\mu: \model^P\otimes \model_P\to \model$$
which is, however, not an isomorphism of complexes since the Leibniz formula
$$\modeld\mu=\mu(\modeld^P\otimes\id) + \mu(\id \otimes \modeld_P)$$
does not hold. To fix this, we consider the filtration $$\overline{W}^q\model=\mu(W^q\model^P\otimes \model_P)$$ induced by the weight filtration of $\model^P.$
\begin{proposition}
	The differential $\modeld$ is compatible with the filtration $\overline{W}^q.$ Moreover, the map $\mu$ induces an isomorphism of complexes
	$$\mu: \model^P\otimes \model_P\to \gr_{\overline{W}}\model.$$
\end{proposition}
\begin{proof} Let $v_i\in \bW^P,$ $Y\in \bigwedge y_d$, $w\in \bW_P$ and $Y'\in \bigwedge y_1,\dots,y_{d-1}.$ We have $d(\dd{v_i}Y\dd{w}Y')=d(\dd{v_i}\dd{w})YY',$ so we may assume that $Y=Y'=1.$
	The set of reflections is a disjoint union $\bR=\bR_P\sqcup \{t_1,\dots,t_{d-1}\}$ where $\bR_P$ is the set of reflections/transpositions in $\bW_P$ and the $t_i$ are from \eqref{eq:chainofreflections}.
	This decomposition is preserved by conjugation with $\bW_P.$ For all $t\in \bR_P$ and $w\in \bW_P$ one has $\ell(v_iwt)=\ell(v_iw)+1$ if and only if $\ell(wt)=\ell(w)+1.$ Moreover, for all $t\in \{t_1,\dots,t_{d-1}\}$ one has $\ell(v_iwt)=\ell(v_iw)+1$ if and only if $wtw\inv=t_{i+1}.$ Using this, we compute the differential of the product as
\begin{align*}
	d(\dd{v_i}\dd{w})
		&=\sum_{\mathclap{\substack{t\in \bR\\ \ell(v_iwt)=\ell(v_iw)+1}}}\dd{v_iwt}\alpha^\vee_t=(\dd{v_i}\sum_{\mathclap{\substack{t\in \bR_P\\ \ell(wt)=\ell(w)+1}}}\dd{wt}\alpha^\vee_t)
		+ \dd{v_{i+1}}\alpha^\vee_{t_{i+1}}\dd{w}\\
		&=\delta_{v_i}\modeld_P(\delta_w)+\modeld^P(\dd{v_i})\dd{w}-\dd{v_{i+1}}y_{d-i+1}\dd{w}.
\end{align*}
While the first two terms of the result have the same degree with respect to $\overline W$, the last term is higher by $1$. It follows that the differential $d$ is compatible with $\overline W.$ Moreover, since $d$ satisfies the Leibniz rule in $\gr_{\overline{W}}\model$ the second claim follows.
\end{proof}
\begin{proof}[Sketch of proof of Theorem~\ref{thm:spectralsequencecohomologyglddegenerates2}]
	The filtration $\overline{W}$ on $A$ yields a spectral sequence $E'_1=\model^P\otimes \model_P\Rightarrow \model.$ Hence, we have two spectral sequences
	\begin{align} \label{eq:trapped}
		\model^P\otimes \model_P\Rightarrow \gr_{\overline{W}}\model\quad \text{ and }\quad \model \Rightarrow  H_c^\bullet (\bG)
	\end{align}

	The spectral sequence $E''_1=\model_P\Rightarrow H_c^\bullet (\bG_P)$ degenerates on $E_2$-page by induction. We obtain
	\begin{align*}
		E_2'&=H^\bullet(\model^P\otimes \model_P)=H^\bullet(\model^P)\otimes H^\bullet(\model_P)\\
		&\cong H_c^\bullet (\A^d-\{0\})\otimes H_c^\bullet (\bG_P)\cong H_c^\bullet (\bG).
	\end{align*}
	Thus, the two spectral sequences in \eqref{eq:trapped} trap the second page $E_2=H^\bullet(\model)$: Namely, the first spectral sequence shows that $E_2$ is (isomorphic to) a subquotient of $H_c^\bullet (\bG)$, while the second one shows that $H_c^\bullet (\bG)$ is (isomorphic to) a subquotient of $E_2.$ It follows that $E_2\cong H_c^\bullet (\bG)$ and both sequences degenerate on page $2.$
	
	For details on `composing' spectral sequences for two filtrations see \cite{deligneTheorieHodgeII1971}.
\end{proof}

\section{A model for the cohomology of open Richardson varieties}
\label{sec:modelforrichardsson}
 We connect now the model for locally closed subsets of $\GL_d$ from  \Cref{sec:cohomologygld} with the cohomology with compact support for open Richardson varieties in $\Gr(d,n)$.
\subsection{Combinatorial model} Let $I=\{i_1<\dots<i_d\}, J=\{j_1<\dots<j_d\}\in \setbinom{\firstn{n}}{d}$ such that $I\leq J.$ We assume that $I\cap J=\emptyset,$ for ease of notation\footnote{For $I'=I-J$ and $J'=J-I$ there is a vector bundle $\rich{I}{J}\to\rich{I'}{J'},$ see \Cref{thm:reductiontodisjoint}.}.
There is a chain of vector bundles compatible with the \Deodhar\ of the open Richardson variety $\rich{I}{J}\subset Z=\Gr(d,n)$ as follows
\[\begin{tikzcd}
	{\rich{I}{J}} & {\basechangeimage{I}{J}} & {\bigsqcup_{w\in \basechangeweyl{I}{J}}\bB w\bB\subseteq\bG=\GL_d} \\
	{\richstrat{I}{J}{w}} & {\basechangeimagestrat{I}{J}{w}} & {\bB w\bB.}
	\arrow[two heads, from=1-1, to=1-2]
	\arrow[two heads, from=1-3, to=1-2]
	\arrow[hook, from=2-1, to=1-1]
	\arrow[hook, from=2-3, to=1-3]
	\arrow[hook, from=2-2, to=1-2]
	\arrow[two heads, from=2-1, to=2-2]
	\arrow[two heads, from=2-3, to=2-2]
\end{tikzcd}\]
By \Cref{makedisjoint}, $\basechangeweyl{I}{J}\subset \bW$ is closed with respect to the Bruhat order. Using \Cref{thm:richstrateverything}, one computes the difference of ranks of the two vector bundles as $$r=\beta-\dim\bU=\ell(J)-\ell(I)-|\setbuild{(k,k')}{i_k<j_{k'}<j_{k}}|-\binom{d}{2}.$$

We hence obtain a model for the cohomology with compact support of 
$\rich{I}{J}$ from the results of Sections \ref{sec:locallyclosedsubsetsgld} and \ref{sec:degenerationviaweights}. For $p,q$ let 
$$\richmodel{I}{J}^{p,q}=\model_{\basechangeweyl{I}{J}}^{p-2r,q-r}$$
and equip $\richmodel{I}{J}$ with the differential $\modeld$ of $\model_{\basechangeweyl{I}{J}}.$ We obtain the following result.
\begin{theorem}[Model]\label{thm:modelrichardsongrassmannian}
    There is an isomorphism of graded vector spaces 
    \begin{align*}H^{\bullet}(\richmodel{I}{J},\modeld)\cong\gr^\bullet_{F}H^\bullet_c(\rich{I}{J}),
    \end{align*}
    where $F$ is a filtration on $H^\bullet_c(\rich{I}{J})$ induced from the \Deodhar.
\end{theorem}

\subsection{The Fukaya model} The model $\richmodel{I}{J}$ has a convenient diagrammatical basis by cohomological \Deodhartype\ Fukaya diagrams, see \Cref{def:cohomologicalmultiplicativefukayadiagram}.
\begin{definition}{\rm
The cohomological \Deodhartype\ Fukaya diagram associated to a basis vector $\dd{w}y_{k_1}\cdots y_{k_l}\in\richmodel{I}{J}$ for $k_1<\dots<k_l$ and $w\in \basechangeweyl{I}{J}$ is constructed from the \Deodhartype\ Fukaya diagram associated to $w$ in the following way. All $\Gm$-nodes in the columns $j_{k_1},\dots,j_{k_l}$ are replaced by circles (representing the generator of $H^1_c(\Gm)$) and all other $\Gm$-nodes by crosses (representing the generator of $H^2_c(\Gm)$).}
\end{definition}
The bidegree of a cohomological Fukaya diagram is given by $(2a+b,a+b)$ where $a$ is the number of crosses and $b$ is the number of circles in the diagram.

Moreover, the differential $\modeld$ of $\richmodel{I}{J}$, see \eqref{eq:differentialinmodel} and \eqref{eq:differentialinmodeldualbasis}, can be expressed by a simple set of local rules on the diagrams which we sketch now.

For each rectangle spanned by a matched node in the upper left and lower right corner \emph{which does not contain any other matched node}\footnote{This reflects the condition that the differential decreases the length of $\ddd{w}$ by one.},we obtain, depending on the decoration of the matched nodes, the following local rules for the differential $d$, describing the linear combinations of modified diagrams in the image of $d$.
\begin{eqnarray}\label{eq:differentialdiagone}\begin{tikzpicture}[baseline={(0, 0)}]
	\def\n{4}
	\def\d{3}
	\pgfmathtruncatemacro{\dpp}{\d}
	\pgfmathtruncatemacro{\npp}{\n}

	\let\mymatrixcontent\empty
		\foreach \row in {1,...,\dpp}{
			\foreach \col in {0,...,\n}{%
				\xappto\mymatrixcontent{ \expandonce{\&}}
	}%
			\gappto\mymatrixcontent{\\}
	}
	\matrix (m) [matrix of nodes, nodes in empty cells, ampersand replacement=\&,
				matrixstyle ] (m){
	\mymatrixcontent 
	};

	
	\foreach \x in {1,4}{
	  \draw[verticallines] (m-1-\x.north) -- (m-1-\x.south);
      \draw[verticallines, dotted] (m-2-\x.north) -- (m-2-\x.south);
      \draw[verticallines] (m-3-\x.north) -- (m-3-\x.south);
    }
	\foreach \k in {1,3}{ 
		\draw[horizontallines] (m-\k-1.west) -- (m-\k-1.east);
        \draw[horizontallines,dotted] (m-\k-1.east) -- (m-\k-4.west);
		\draw[horizontallines] (m-\k-4.west) -- (m-\k-4.east);

    }
	 \foreach \x/\y in {1/1,3/4}{
    	\node[matchingnodes] at (m-\x-\y.center) {};
        \node[aonenodes] at (m-\x-\y.center) {};
     }
\end{tikzpicture}
\hspace{-30pt}
&\mapsto(-1)^{c_2}&
\hspace{-5pt}
\begin{tikzpicture}[baseline={(0, 0)}]
	\def\n{4}
	\def\d{3}
	\pgfmathtruncatemacro{\dpp}{\d}
	\pgfmathtruncatemacro{\npp}{\n}

	\let\mymatrixcontent\empty
		\foreach \row in {1,...,\dpp}{
			\foreach \col in {0,...,\n}{%
				\xappto\mymatrixcontent{ \expandonce{\&}}
	}%
			\gappto\mymatrixcontent{\\}
	}
	\matrix (m) [matrix of nodes, nodes in empty cells, ampersand replacement=\&,
				matrixstyle ] (m){
	\mymatrixcontent 
	};

	
	\foreach \x in {1,4}{
	  \draw[verticallines] (m-1-\x.north) -- (m-1-\x.south);
      \draw[verticallines, dotted] (m-2-\x.north) -- (m-2-\x.south);
      \draw[verticallines] (m-3-\x.north) -- (m-3-\x.south);
    }
	\foreach \k in {1,3}{ 
		\draw[horizontallines] (m-\k-1.west) -- (m-\k-1.east);
        \draw[horizontallines,dotted] (m-\k-1.east) -- (m-\k-4.west);
		\draw[horizontallines] (m-\k-4.west) -- (m-\k-4.east);

    }
	 \foreach \x/\y in {1/4,3/1}{
    	\node[matchingnodes] at (m-\x-\y.center) {};
     }
     \foreach \x/\y in {1/4}{
    	\node[aonenodes] at (m-\x-\y.center) {};
     }
     \foreach \x/\y in {3/1}{
    	\node[gmnodes] at (m-\x-\y.center) {};
     }
\end{tikzpicture}
\hspace{-30pt}
-(-1)^{c_1}
\hspace{-5pt}
\begin{tikzpicture}[baseline={(0, 0)}]
	\def\n{4}
	\def\d{3}
	\pgfmathtruncatemacro{\dpp}{\d}
	\pgfmathtruncatemacro{\npp}{\n}

	\let\mymatrixcontent\empty
		\foreach \row in {1,...,\dpp}{
			\foreach \col in {0,...,\n}{%
				\xappto\mymatrixcontent{ \expandonce{\&}}
	}%
			\gappto\mymatrixcontent{\\}
	}
	\matrix (m) [matrix of nodes, nodes in empty cells, ampersand replacement=\&,
				matrixstyle ] (m){
	\mymatrixcontent 
	};

	
	\foreach \x in {1,4}{
	  \draw[verticallines] (m-1-\x.north) -- (m-1-\x.south);
      \draw[verticallines, dotted] (m-2-\x.north) -- (m-2-\x.south);
      \draw[verticallines] (m-3-\x.north) -- (m-3-\x.south);
    }
	\foreach \k in {1,3}{ 
		\draw[horizontallines] (m-\k-1.west) -- (m-\k-1.east);
        \draw[horizontallines,dotted] (m-\k-1.east) -- (m-\k-4.west);
		\draw[horizontallines] (m-\k-4.west) -- (m-\k-4.east);

    }
	 \foreach \x/\y in {1/4,3/1}{
    	\node[matchingnodes] at (m-\x-\y.center) {};
     }
     \foreach \x/\y in {3/1}{
    	\node[aonenodes] at (m-\x-\y.center) {};
     }
     \foreach \x/\y in {1/4}{
    	\node[gmnodes] at (m-\x-\y.center) {};
     }
\end{tikzpicture}
\end{eqnarray}
\begin{eqnarray}\label{eq:differentialdiagtwo}
    \hspace{-30pt}
(-1)^{c_1}
\hspace{-5pt}
    \begin{tikzpicture}[baseline={(0, 0)}]
	\def\n{4}
	\def\d{3}
	\pgfmathtruncatemacro{\dpp}{\d}
	\pgfmathtruncatemacro{\npp}{\n}

	\let\mymatrixcontent\empty
		\foreach \row in {1,...,\dpp}{
			\foreach \col in {0,...,\n}{%
				\xappto\mymatrixcontent{ \expandonce{\&}}
	}%
			\gappto\mymatrixcontent{\\}
	}
	\matrix (m) [matrix of nodes, nodes in empty cells, ampersand replacement=\&,
				matrixstyle ] (m){
	\mymatrixcontent 
	};

	
	\foreach \x in {1,4}{
	  \draw[verticallines] (m-1-\x.north) -- (m-1-\x.south);
      \draw[verticallines, dotted] (m-2-\x.north) -- (m-2-\x.south);
      \draw[verticallines] (m-3-\x.north) -- (m-3-\x.south);
    }
	\foreach \k in {1,3}{ 
		\draw[horizontallines] (m-\k-1.west) -- (m-\k-1.east);
        \draw[horizontallines,dotted] (m-\k-1.east) -- (m-\k-4.west);
		\draw[horizontallines] (m-\k-4.west) -- (m-\k-4.east);

    }
    \foreach \x/\y in {1/1,3/4}{
    	\node[matchingnodes] at (m-\x-\y.center) {};
     }
	\foreach \x/\y in {1/1}{
    	\node[aonenodes] at (m-\x-\y.center) {};
     }
     \foreach \x/\y in {3/4}{
    	\node[gmnodes] at (m-\x-\y.center) {};
     }
\end{tikzpicture}
\hspace{-30pt}
, -(-1)^{c_2}
\hspace{-5pt}
\begin{tikzpicture}[baseline={(0, 0)}]
	\def\n{4}
	\def\d{3}
	\pgfmathtruncatemacro{\dpp}{\d}
	\pgfmathtruncatemacro{\npp}{\n}

	\let\mymatrixcontent\empty
		\foreach \row in {1,...,\dpp}{
			\foreach \col in {0,...,\n}{%
				\xappto\mymatrixcontent{ \expandonce{\&}}
	}%
			\gappto\mymatrixcontent{\\}
	}
	\matrix (m) [matrix of nodes, nodes in empty cells, ampersand replacement=\&,
				matrixstyle ] (m){
	\mymatrixcontent 
	};

	
	\foreach \x in {1,4}{
	  \draw[verticallines] (m-1-\x.north) -- (m-1-\x.south);
      \draw[verticallines, dotted] (m-2-\x.north) -- (m-2-\x.south);
      \draw[verticallines] (m-3-\x.north) -- (m-3-\x.south);
    }
	\foreach \k in {1,3}{ 
		\draw[horizontallines] (m-\k-1.west) -- (m-\k-1.east);
        \draw[horizontallines,dotted] (m-\k-1.east) -- (m-\k-4.west);
		\draw[horizontallines] (m-\k-4.west) -- (m-\k-4.east);

    }
	 \foreach \x/\y in {1/1,3/4}{
    	\node[matchingnodes] at (m-\x-\y.center) {};
     }
     \foreach \x/\y in {3/4}{
    	\node[aonenodes] at (m-\x-\y.center) {};
     }
     \foreach \x/\y in {1/1}{
    	\node[gmnodes] at (m-\x-\y.center) {};
     }
\end{tikzpicture}
\hspace{-30pt}
&\mapsto&
\hspace{-5pt}
\begin{tikzpicture}[baseline={(0, 0)}]
	\def\n{4}
	\def\d{3}
	\pgfmathtruncatemacro{\dpp}{\d}
	\pgfmathtruncatemacro{\npp}{\n}

	\let\mymatrixcontent\empty
		\foreach \row in {1,...,\dpp}{
			\foreach \col in {0,...,\n}{%
				\xappto\mymatrixcontent{ \expandonce{\&}}
	}%
			\gappto\mymatrixcontent{\\}
	}
	\matrix (m) [matrix of nodes, nodes in empty cells, ampersand replacement=\&,
				matrixstyle ] (m){
	\mymatrixcontent 
	};

	
	\foreach \x in {1,4}{
	  \draw[verticallines] (m-1-\x.north) -- (m-1-\x.south);
      \draw[verticallines, dotted] (m-2-\x.north) -- (m-2-\x.south);
      \draw[verticallines] (m-3-\x.north) -- (m-3-\x.south);
    }
	\foreach \k in {1,3}{ 
		\draw[horizontallines] (m-\k-1.west) -- (m-\k-1.east);
        \draw[horizontallines,dotted] (m-\k-1.east) -- (m-\k-4.west);
		\draw[horizontallines] (m-\k-4.west) -- (m-\k-4.east);

    }
	 \foreach \x/\y in {1/4,3/1}{
    	\node[matchingnodes] at (m-\x-\y.center) {};
     }
     \foreach \x/\y in {3/1,1/4}{
    	\node[gmnodes] at (m-\x-\y.center) {};
     }
\end{tikzpicture}
\end{eqnarray}

Here, $c_1$ and $c_2$ are the number of circled edges in columns of the top left diagram which are to the left of the first and second red strand, respectively.

\section{Standard extension algebra: Explicit geometric description}

\subsection{Parametrizations and gluing matrices}\label{sec:gluingofmatrices}
Let $I=\{i_1<\dots<i_d\}\in \setbinom{\firstn{n}}{d}.$
Recall that $e_I=(I_n)_{I\times \firstn{n}}\in k^{d\times n}$ denotes the matrix whose $k$-th row is the $i_k$-th standard basis vector. There are isomorphism

\[\begin{tikzcd}[row sep=tiny]
	{U_I,\opp{U}_I} & {M_I^\pm} & {Z_I^\pm} \\
	u & {e_Iu\tr;\, X} & {[X]}
	\arrow["\sim", from=1-1, to=1-2]
	\arrow["\sim", from=1-2, to=1-3]
	\arrow[maps to, from=2-1, to=2-2]
	\arrow[maps to, from=2-2, to=2-3]
\end{tikzcd}\]
parametrizing the Bruhat cells, see \Cref{sec:notationsflagvarieties} and \Cref{sec:parabruhat}. Moreover, in the Grassmannian case $[\opp U_I,U_I]=1$ and we obtain a parametrization of the open affine neighborhood $V_I$ from \eqref{eq:parabruhat} of the fixed point $[e_I]$ via

\[\begin{tikzcd}[row sep=tiny]
	{\opp U_I\times U_I} & {M_I} & {V_I} \\
	{(u_1,u_2)} & {e_Iu_1\tr u_2\tr;\, X} & {[X].}
	\arrow["\sim", from=1-1, to=1-2]
	\arrow["\sim", from=1-2, to=1-3]
	\arrow[maps to, from=2-1, to=2-2]
	\arrow[maps to, from=2-2, to=2-3]
\end{tikzcd}\]
Combining these isomorphisms, we obtain a commutative diagram

\begin{equation}\label{eq:gluingmatrices}
	\begin{tikzcd}[row sep=0]
	{} & {M_I^-\times M_I^+} & {Z_I^-\times Z_I^+} \\
	{\opp U_I\times U_I} \\
	{} & {M_I} & {V_I}
	\arrow["\gamma_I", from=1-2, to=3-2]
	\arrow[from=2-1, to=1-2]
	\arrow[from=1-2, to=1-3]
	\arrow[from=1-3, to=3-3]
	\arrow[from=3-2, to=3-3]
	\arrow[from=2-1, to=3-2]
\end{tikzcd}
\end{equation}
The isomorphism $\gamma_I$ `glues together two matrices'. It maps $(X,Y)$ to $Z$ where
$$Z_{k,l}=\begin{cases}
    X_{k,l} & \text{ if } l\leq i_k,\\
    Z_{k,l} & \text{ if } l\geq i_k,
\end{cases}$$
in other words, $Z=X+Y-e_I.$
\subsection{Multiplication and normal forms} Consider now $H=\{h_1<\dots<h_d\}$ and $J=\{j_1<\dots<j_d\}\in \setbinom{\firstn{n}}{d}$ with $H\leq I\leq J.$

Using \eqref{eq:gluingmatrices}, we obtain the commutative diagram
\begin{equation}\label{eq:multiplicationinnormalform}
\begin{tikzcd}
	{\rich{H}{I}\times \rich{I}{J}} & {\rich{H}{J}} \\
	{(\bG M_H^-\cap M_I^+)\times (M_I^-\cap\bG M_J^-)}
	\arrow["{\psi_I}", from=1-1, to=1-2]
	\arrow["{[_-]\times [_-]}"', from=2-1, to=1-1]
	\arrow["{[_-]\circ\gamma_I}"', from=2-1, to=1-2]
\end{tikzcd}
\end{equation}
where $\psi_I$ is the map defined in \eqref{eq:multiplicationofrichardsons}. We can therefore describe the composition of extensions of the standard object $\Delta\in \D(Z)$ in terms of the gluing map $\gamma_I.$
\begin{corollary}[Composition of Extensions]
There is a commutative diagram
\[\begin{tikzcd}
	{\Hom(\Delta_H,\Delta_I[\bullet])\otimes\Hom(\Delta_I,\Delta_J[\bullet])} & {\Hom(\Delta_H,\Delta_J[\bullet])} \\
	{H^\bullet_c(\rich{H}{I})\otimes H^\bullet_c(\rich{I}{J})} & {H^\bullet_c(\rich{H}{J})} \\
	{H^\bullet_c((\bG M_H^-\cap M_I^+)\times (M_I^-\cap\bG M_J^-))}
	\arrow["\wr", from=1-1, to=2-1]
	\arrow["\wr"', from=3-1, to=2-1]
	\arrow["\circ^{\operatorname{op}}", from=1-1, to=1-2]
	\arrow["{([_-]\circ\gamma_I)_!}"', from=3-1, to=2-2]
	\arrow["\wr"', from=1-2, to=2-2]
	\arrow["{(\psi_{I})_!}", from=2-1, to=2-2]
\end{tikzcd}\]
\end{corollary}
\subsection{Multiplication and base change}
For simplicity, assume that $H\cap J=\emptyset.$ In this case, the explicit description of the map $\gamma_I$ and \eqref{eq:multiplicationinnormalform} allows to show that the following diagram is commutative
\begin{equation}\label{eq:multiplicationbasechange}
	\begin{tikzcd}
	{\rich{H}{I}\times\rich{I}{J}} & {\rich{H}{J}} \\
	{\basechangeimage{H}{I}\times \basechangeimage{I}{J}} & {\basechangeimage{H}{J}}
	\arrow["{\psi_{I}}", from=1-1, to=1-2]
	\arrow["{\basechange{H}{I}\times \basechange{I}{J}}", two heads, from=1-1, to=2-1]
	\arrow["m",from=2-1, to=2-2]
	\arrow["{\basechange{H}{J}}", two heads, from=1-2, to=2-2]
\end{tikzcd}
\end{equation}
where the $\varphi$'s are the respective base change maps introduced in  \Cref{basechange} and $m$ is the multiplication map in the group $\bG.$
The base change maps $\varphi$ are vector bundles by \Cref{thm:reductiontobasechangeimage}. Denote their respective ranks by $r,r'$ and $r''.$ We obtain the following explicit description of the compositions of extensions.
\begin{corollary}[Exts via cohomology of $\bG$]
	There is a commutative diagram
\[\begin{tikzcd}
	{\Hom(\Delta_H,\Delta_I[\bullet])\otimes\Hom(\Delta_I,\Delta_J[\bullet])} & {\Hom(\Delta_H,\Delta_J[\bullet])} \\
	{H^{2r+\bullet}_c(\basechangeimage{H}{I})\otimes H^{2r'+\bullet}_c(\basechangeimage{I}{J})} & {H^{2r''+\bullet}_c(\basechangeimage{H}{J})}
	\arrow["\wr", from=1-1, to=2-1]
	\arrow["\circ^{\operatorname{op}}", from=1-1, to=1-2]
	\arrow["\wr"', from=1-2, to=2-2]
	\arrow["{m_!}", from=2-1, to=2-2]
\end{tikzcd}\]
\end{corollary}
\begin{example}
Write $I\ll J$ if $i<j$ for all $i\in I$ and $j\in J.$ If $H\ll I\ll J$, then $\basechangeimage{H}{I},\basechangeimage{I}{J},\basechangeimage{H}{J}=\bG.$ Denoting by $e_1=y_1+\dots+y_d, \dots, e_d=y_1\cdots y_d$ the elementary symmetric polynomials, one can identify
$$H_c^{2\dim\bG-\bullet}(\bG)=H_\bullet(\bG)=\mywedge^\bullet e_1,\dots,e_d$$
where $e_i$ sits in homological degree $2i-1.$ Since the group multiplication induces the wedge product on the exterior algebra, we can describe the composition of extensions of standard objects (up to appropriate shifts) simply via the wedge product
\[\begin{tikzcd}
	{\Hom(\Delta_H,\Delta_I[\bullet])\otimes\Hom(\Delta_I,\Delta_J[\bullet])} & {\Hom(\Delta_H,\Delta_J[\bullet])} \\
	{\mywedge^\bullet e_1,\dots,e_d\otimes \mywedge^\bullet e_1,\dots,e_d} & {\mywedge^\bullet e_1,\dots,e_d}.
	\arrow["\wr", from=1-1, to=2-1]
	\arrow["\circ^{\operatorname{op}}", from=1-1, to=1-2]
	\arrow["\wr"', from=1-2, to=2-2]
	\arrow["{-\wedge-}", from=2-1, to=2-2]
\end{tikzcd}\]
\end{example}
\begin{remark} In case that $H\cap J\neq \emptyset,$ the diagram \eqref{eq:multiplicationbasechange} may not commute. Rather, one needs to twist the multiplication map by certain factors depending on the intersecting columns indexed by $H\cap J.$ We will not go into details here.
\end{remark}
\subsection{A $dg$-model for sheaves on Grassmannians?}
Recall from \Cref{sec:modelforrichardsson} the model $\richmodel{I}{J}=\model_{\basechangeweyl{I}{J}}$ describing the cohomology with compact support of $\rich{I}{J}.$
We define 
$\model_{d,n}=\bigoplus_{I\leq J}\richmodel{I}{J}.$ Then $\model_{d,n}$ inherits a differential $d$ from each $\richmodel{I}{J}.$ Using \Cref{thm:modelrichardsongrassmannian} we finally connect to the standard extension algebra  $\op{E}(\Delta)$ from \Cref{Extalg}:
\begin{corollary}\label{dg}
	There is a filtration $F$ on $\op{E}(\Delta)$ such that
	$$H^\bullet(\model_{d,n},d)=\gr_F^\bullet \op{E}(\Delta).$$
\end{corollary}
\begin{remark}
	This leads to the natural questions of how to turn $\model_{d,n}$ into a $dg$-algebra whose cohomology yields the standard extension algebra $\op{E}(\Delta)$.  Clearly, the multiplication should satisfy $\richmodel{H}{I}\richmodel{I'}{J}=0$ if $I\neq I'.$ In the case $H\cap J=\emptyset,$ one may define a multiplication via the commutative diagram
\[\begin{tikzcd}
	\model\otimes\model & \model \\
	{\richmodel{H}{I}\otimes \richmodel{I}{J}} & {\richmodel{H}{J}}
	\arrow[two heads, from=1-1, to=2-1]
	\arrow[from=1-1, to=1-2]
	\arrow[two heads, from=1-2, to=2-2]
	\arrow[from=2-1, to=2-2]
\end{tikzcd}\]
where the top row is the multiplication for $\model$ introduced in  \Cref{sec:combinatorialmodel}. 
\end{remark}

\section{Extensions of Verma modules in category $\pazocal O$}\label{catO}
In the section, we connect the geometric and diagrammatical results with extensions of Verma modules in parabolic category $\cO.$  Since geometric and algebraic conventions often differ in subtle ways, we start with grading conventions. 
\subsection{Grading conventions}\label{conventions}
If  $\Aa$ is a category of modules over an algebra $A$ which carries a $\Z$-grading we denote by $\hat{\Aa}$ the category of graded $A$-modules with degree preserving morphisms. This is the a {\emph{graded version} of $\Aa$. We denote by $v^i$ the grading shift functor up by $i$, that means $(v^iM)_i=M_0$.  If $M,N\in \hat{\Aa}$ are considered as objects in $\Aa$, then $\Hom_{\Aa}(M,N)$ is $\Z$-graded with degree $i$ part given by 
$\Hom_{\hat\Aa}(M,v^{-i}N)=\Hom_{\hat\Aa}(v^{i}M,N)$. Analogously statements hold for the derived categories and homotopy categories of projective modules. We denote by $t^i$ the shift in the homological grading, e.g. $(t\Ext^\bullet)^i=\Ext^{i-1}$. If $M,N\in \hat{\Aa}$ are considered as objects in $\Aa$, then $\Ext_{\Aa}(M,N)$ is bigraded with $\Ext^{i,j}_{\Aa}(M,N)=\Ext^{i-1}_{\hat\Aa}(M,v^{-i}N).$ 

\subsection{Category $\cOp$}
   Let  $\mg\supset\mp\supset\mb\supset\mh$ be a semisimple complex Lie algebra with fixed compatible parabolic, Borel and Cartan Lie subalgebras. We denote by $\rho$ the half-sum of positive roots. Let $\pW$ the Weyl group with parabolic subgroup\footnote{The notation is chosen to agree with the previous ones if $\mb$, $\mp$ are the Lie algebras of $B$, $P$.} $\pW_P$ associated with $\mp$. We write $w_0$ and $w_{P}$ for the  longest elements in $\pW$ and $\pW_P$ respectively. Let $\WP$ and $\WPl$ be the set of minimal length coset representatives for $\pW/\pW_P$ respectively $\pW_P \backslash \pW$. 
   
  We denote by $\cOp$ the parabolic category $\cO$ associated to  $(\mg,\mp,\mh)$, that means the category of representations $M$ of $\mg$ which have a weight space decomposition and are a union of finite dimensional subrepresentations when restricted to $\mp$.  We focus on the principal block $\cOpo$ that means we additionally require the representations to have generalized trivial central character\footnote{As usual: regular blocks are equivalent,  nonintegral blocks are treated by passing to integral Weyl groups, parabolic-singular blocks are equivalent to regular blocks for smaller rank.}, see \cite{humphreysRepresentationsSemisimpleLie2008} for more details. Representatives for the isomorphism classes of simple objects in $\cOpo$ are the irreducible highest weight modules $L(x\cdot 0)$ of highest weight $x\cdot0:=x(\rho)-\rho$ for $x\in \WPl$.  
   
    \subsection{Lie theoretic extension algebra}  
   Category $\cOpo$ is a highest weight category with standard objects the parabolic Verma modules $\Delta^\mp(x\cdot0)$ with $x\in \WPl$ and the usual ordering on weights. This is opposite to the Bruhat ordering in $\WPl$:   
    \begin{lemma}\label{stupidlemma}
 For  $y,x\in \WPl$ and $z\in \pW$ the following holds:
 \begin{enumerate}
  \item $x\cdot 0<y\cdot0$ $\Leftrightarrow$  $x > y $  $\Leftrightarrow$ $w_P xw_0<w_P yw_0$.
 \item $z\in \WP$ $\Leftrightarrow$ $z^{-1}\in \WPl$ and $w_P zw_0\in\WP$ $\Leftrightarrow$ $z\in\WP$.
 \end{enumerate}
    \end{lemma}
 \begin{definition}{\rm
   We define the \emph{Lie theoretic standard extension algebra} as
 \begin{eqnarray*}
  \op{E}(\Delta^\mp)&=&\Ext^\bullet(\Delta^\mp,\Delta^\mp)=
\bigoplus_{x,y\in\WPl}\Ext^\bullet_{\cO^\mp}(\Delta^\mp(x\cdot0),\Delta^\mp(y\cdot0))
   \end{eqnarray*}}
  \end{definition}
  By the highest weight structure, $\Ext^\bullet_{\cO^\mp}(\Delta^\mp(x\cdot0),\Delta^\mp(y\cdot0))=0$ unless $x\cdot0\leq y\cdot0$. 
  To make the connection with \Cref{sec:notationsflagvarieties} let $G\supset P\supset B\supset T$ be complex reductive algebraic groups with Lie algebras $\mg\supset\mp\supset\mb\subset\mh$ and let $Z=G/P.$  We can realise $\cOpo$ geometrically, via the localisation theorem and the Riemann--Hilbert correspondence, as the category of perverse sheaves on  $Z$ with the  stratification by Bruhat cells, see \cite[\S{12}]{hottaDmodulesPerverseSheaves2008}. The parabolic Verma module $\Delta^\mp(x\cdot0)$ corresponds to the constant perverse sheaf supported on $Z_\psi(x)$, where we use the isomorphism of posets $\psi:\WPl\iso Z^0$, $\psi(x)=w_0 x^{-1}w_P$ from \Cref{stupidlemma}. We obtain directly from \Cref{thm:geometricdescriptionofcomposition} (keeping track of perversity shifts)
\begin{proposition}
  There is an isomorphism of graded algebras $ E(\Delta^\mp)\cong\op{E}(\Delta)$.
  which induces isomorphisms  of graded vector spaces
  \begin{eqnarray} \label{graded}
  \Ext^\bullet_{\cOp}(\Delta^\mp(x\cdot0),\Delta^\mp(y\cdot0))&\cong&H^{\bullet+\ell(y)-\ell(x)}_c(Z^-_{\psi(x)}\times_ZZ_{\psi(y)}^+).
  \end{eqnarray}
  \end{proposition}
\begin{remark}
The isomorphism \eqref{graded} was already established in \cite{richeModularKoszulDuality2014}. Our constructions,  \Cref{thm:twocompositionsagree}, provide two alternative definitions of the multiplication. 
\end{remark}
 As mentioned already in \Cref{sec:bigradings}, it is most natural to consider all spaces as bigraded vector spaces. For category $\cO$ the second grading is the extra mixed Hodge or Koszul grading from \cite{beilinsonKoszulDualityPatterns1996}. It was made algebraically explicit in \cite{Strgradings}, \cite{stroppelCategorificationTemperleyLiebCategory2005} using Soergel modules and a graded categorification of the antispherical parabolic Hecke module. Important for us is that a graded version $\hat\cO_0^\mp$ of $\cO_0^\mp$ exists and $\Delta^\mp(x\cdot0)$ has a graded lift $\Verma{x}$ which is unique up to isomorphism and overall grading shift. We choose $\Verma{x}$ with head  concentrated in degree zero. Then we obtain, using the grading conventions from \Cref{conventions}: 
 \begin{theorem}\label{thm:catO}
There is an isomorphism of  bigraded vector spaces
  \begin{eqnarray*}
    \Ext_{\cO_0^\mp}^{i,j}(\Delta^\mp(x\cdot0),\Delta^\mp(y\cdot0))\cong H^{\ell(x)-\ell(y)+i,\frac{1}{2}(\ell(x)-\ell(y)-j)}_c(Z^-_{\psi(x)}\times_ZZ_{\psi(y)}^+)
       \end{eqnarray*}
which gives rise to an isomorphism of bigraded algebras   $E(\Delta^\mp)\cong\op{E}(\Delta)$.  
\end{theorem}
In particular, the mixed Hodge polynomial on the geometric side corresponds to the bigraded Poincar\'e polynomial of the Lie-theoretic standard extension algebra,
$$\mixPol^\mp(x,y;v,t)=\sum_{i,j}\dim \Ext_{\cO_0^\mp}^{i,j}(\Delta^\mp(x\cdot0),\Delta^\mp(y\cdot0))v^jt^i.$$ 
\subsection{Translation and shuffling functors and projective resolutions}
 For a simple reflection $s$ let $\theta_s:\cO_0^\mp\rightarrow\cO_0^\mp$ be the \emph{translation functor} through the $s$-wall and $\Theta_s:\hat\cO_0^\mp\rightarrow\hat\cO_0^\mp$ its graded lift from \cite{Strgradings}. The usual adjunctions of translation to, off and through walls lift to the graded setting,  \cite{Strgradings}:
\begin{lemma}\label{transfunc}
Let $s$ be a simple reflection and $z\in\WPl$. Then it holds in $\hat\cO_0^\mp$:
\begin{enumerate}
\item The pair $(\Theta_s, \Theta_s)$ is a pair of adjoint functors.
\item Let $zs>z$. If $zs\notin\WPl$, then $\Theta_s\Verma{z}=0$, otherwise $\Theta_s\Verma{z}\cong v\Theta_s\Verma{zs}$ and there are short exact sequences
\begin{eqnarray}
v\Verma{z} \ra{\adj{s}}&\Theta_s\Verma{z}&\ra{\phantom{\adjov{s}}} \Verma{zs},\label{ses}\\
\Verma{z}\ra{\phantom{\adjov{s}}} &{\Theta_s\Verma{zs}}&\ra{\adjov{s}}v^{-1}\Verma{zs}.\label{ses2}
\end{eqnarray}
where $\adj{s}$ and ${\adjov{s}}$ are the graded adjunction morphism. 
\end{enumerate}
\end{lemma} 

Consider $\Sh_s:=\left(v\op{I}\xrightarrow{a_s}\Theta_s\right)$, the {\it Rouquier complex} of functors. The {\it shuffling functor} induced by $\Sh_s$ on the bounded derived category of $\gcO$ or the homotopy category of projectives in $\gcO$ (by taking the total complex of the tensor product with $\Sh_s$)
is an autoequivalence, and \eqref{ses} implies that $\Verma{x}$ is sent to $\Verma{xs}$ (respectively $\PP{x}$ to $\PP{xs}$) if $xs>x$, see \cite[Theorem 5.7]{mazorchukTranslationShufflingProjectively2005a}. 
\begin{proposition}\label{projres}
Given a projective resolution $\PP{x}$ of $\Verma{x}$ with $x<xs$, the total complex of $\Sh_s \PP{x}$ is a projective resolution of $\Verma{xs}$. 
\end{proposition}
 for a detailed argument see \cite[Theorem 5.3.]{brundanHighestWeightCategories2010}. Since $\Verma{e}\in \gcO$ is projective, this allows to construct inductively minimal projective resolutions $\PP{x}$ of all  $\Verma{x}$'s. This is crucial for computing the standard extension algebra itself and even more, its standard dg-model $(\Hom^\bullet(\oplus_{x}\PP{x},\oplus_y\PP{y}),d)$ or $A_\infty$-model.

\subsection{The Grassmannian case} 
From now on we will assume, parallel to the Grassmannian case, that we are in the special hermitian symmetric case, that means we assume that the parabolic Weyl group is of the form $\pW_P=
S_d\times S_{n-d}\subset S_n= \pW.$ 

It is a special feature of the Grassmannian permutations  that two reduced expressions of $w\in\WPl$ differ only by commutativity braid relations. In particular, the inductively defined projective resolutions $\PP{x}$ in \Cref{projres} are up to canonical isomorphism independent  of the chosen reduced expression of $x$. 

  \subsection{Deodhar combinatorics via translation functors} 
We now mimic Lie theoretically the combinatorics of the \Deodhar\  in the hermitian symmetric case. For this abbreviate 
\begin{eqnarray*}
 \stdE(x,y):=\Ext^\bullet_{\cOpo}(\Verma{x},\Verma{y}),&\;& \stdE(x,y)^r:=\Ext^r_{\cOpo}(\Verma{x},\Verma{y}). 
  \end{eqnarray*}
  We build inductively complexes which compute all $\stdE(x,y)$. We show that the total spaces agree with the geometrically defined complexes by proving recursion formulas analogously to \Cref{prop:recursiondistexpression}  and \Cref{prop:recursionW}.

\begin{theorem}\label{Es}
Let $x,y\in\WPl$ and let $s$ be a simple reflection with $x>xs$. 
\begin{enumerate}[(i)]
\item \label{catI}  Let $y>ys$. Then  
$\stdE(xs,ys)\cong \stdE(x,y)$ as bigraded vector spaces. The isomorphism is given by $f\mapsto F_s(f)$.
\item \label{catII} Let $ys\notin\WPl$. Then $ \stdE(x,y)\cong q^{-1}t \stdE(xs,y)$ as bigraded vector spaces.
\item \label{catIII} Let $y<ys$, $ys\in\WPl$. Then 
$E(x,y)^r$ fits into a short exact sequence 
\begin{equation}\label{tea1}
0\rightarrow v^{-1} \stdE(xs,y)^{r-1}\rightarrow E(x,y)^{r}\rightarrow M^r\rightarrow0, 
\end{equation}
where $M^r=\Ext^\bullet_{\cOpo}(\Theta_s\Verma{xs},\Verma{y})$. It is the $(3r-1)th$ term of a long exact sequence with terms
\begin{equation}\label{tea2}
\ldots\rightarrow v\stdE(xs,y)^{r}\rightarrow M^r\rightarrow \stdE(xs,ys)^{r}\rightarrow\ldots, 
\end{equation}
In particular, $\stdE(x,y)$ is the cohomology of a complex of graded vector spaces with total space bigraded isomorphic to  
\begin{equation}
\label{fritz}
(v^{-1}t+v)\stdE(xs,y)\oplus \stdE(xs,ys),
  \end{equation} 
  \end{enumerate}
\end{theorem}

\begin{proof}
Statement \ref{catI} holds, since $F_s$ induces a derived equivalence sending $\Verma{xs}$ to $\Verma{x}$ and $\Verma{ys}$ to $\Verma{y}$. Applying $\Hom_{\hat\cO_0^\mp}(_-,\Verma{y})$ to the short exact sequence \eqref{ses} with $z=xs$ gives a long exact sequence of the form
\begin{equation}\label{teales}
\ldots\rightarrow v^{-1} \stdE(xs,y)^{r-1}\rightarrow E(x,y)^{r}\rightarrow M^r\ra{f_r} v^{-1} \stdE(xs,y)^{r}\rightarrow\dots, 
\end{equation}
By \Cref{transfunc} $M^r=0$ if $ys\notin\WPl$, thus \ref{catII} holds.  
By \cite[4.3]{sheltonExtensionsGeneralizedVerma1988}, the maps $f_r$ are in fact zero, and \eqref{tea1} holds with $M^r=\Ext_{\hat\cO_0^\mp}(\Theta_s\Verma{xs},\Verma{y})\cong \Ext_{\hat\cO_0^\mp}(\Verma{xs},\Theta_s\Verma{y})$. 
Applying $\Hom_{\hat\cO_0^\mp}(\Verma{xs},_-)$ to \eqref{ses} with $z=y$ gives \eqref{tea2} of \ref{catIII}. For \eqref{fritz}, note that $\stdE(x,y)$ is the cohomology of the complex  $\oplus_{r}M^r\rightarrow\stdE(xs,y)$ by \eqref{tea1} and the claim follows.
\end{proof}
\begin{remark}\label{vanilla}{\rm
The proof implies that the non-zero part of the differential in \eqref{tea2} is
\begin{eqnarray}
\stdE(xs,ys)&\ra{{b_1}\circ_-,{b_2}\circ_-}&(v^{-1}t+v)\stdE(xs,y),
\end{eqnarray}
with homogeneous vectors ${b_1}$, $b_2$ from $\stdE(y,ys)$ of indicated degree.\footnote{Since this space is 2-dimensional, see \Cref{Shelton}, they form in fact a basis.}
}
\end{remark}
We obtain a recursive formula for the bigraded Poincar\'e polynomial $\mixPol^\mp$:
\begin{corollary}[Bigraded Recursion formulas] \label{cor:recursiveformulabigradedpoincare}
  Let $x,y\in\WPl$ and let $s$ be a simple reflection with $x>xs$. Then 
  \begin{eqnarray}
      \mixPol^\mp(x,y)&=&
  \begin{cases}
      \mixPol^\mp(xs,ys)&\text{if $y>ys$},\\
  v^{-1}t \mixPol^\mp(xs,y)&\text{if $ys\notin\WPl$},\\
  (v^{-1}t+v)\mixPol^\mp(xs,y)+\mixPol^\mp(xs,ys)&\text{otherwise.}\\
  \end{cases}
  \end{eqnarray}
\end{corollary}
This is a graded refinement of Shelton's formulas \cite{sheltonExtensionsGeneralizedVerma1988}, see also \cite{biagioliClosedProductFormulas2004}:
\begin{corollary}[Shelton's recursion formulas]\label{Shelton}
For $x,y\in\WPl$ set $e(x,y)^r=\dim\stdE(x,y)^r$. If now $s$ is a simple reflection with $xs<x$, then the following hold. 
\begin{eqnarray*}
e(x,y)^r&=&
\begin{cases}
e(xs,ys)^r&\text{if $y<ys\in\WPl$},\\
e(xs,y)^{r-1}&\text{if $ys\notin\WPl$},\\
e(xs,y)^{r-1}+e(xs,y)^{r}&\text{if $y<ys\in\WPl$, $xs\not\geq ys$},\\
e(xs,y)^{r-1}-e(xs,y)^{r+1}+e(xs,ys)^{r}&\text{if $y<ys<xs<x$}.
\end{cases}
\end{eqnarray*}
\end{corollary}
\begin{proof}
The first two formulas and the last one are direct consequences of \Cref{Es} and its proof, see also \Cref{vanilla}. The second last one is a special case of the last, with $E(xs,ys)=0$, and thus $e(xs,ys)^{r}=0$.
\end{proof}
\subsection{Bases of the standard extension algebra}
\Cref{Es} can be used to inductively control the spaces of extensions and construct  morphisms between chain complexes $\PP{x}$ which give rise to a basis of the standard extension algebra. 
\begin{enumerate}[(c-i)]
\item The isomorphism from \Cref{Es}\ref{catI} can explicitly be constructed on the level of projective resolutions by sending $f: \PP{xs}\rightarrow \PP{ys}$ first to $(f,\Theta_s(f)):\;\left(v\PP{xs}\ra{\adj{s}}\Theta_s\PP{xs}\right) \rightarrow \left(v\PP{ys}\ra{\adj{s}}\Theta_s\PP{ys}\right)$ and then take the induced map $\PP{x}\rightarrow\PP{y}$ between the cones.
\item  The isomorphism in \Cref{Es}\ref{catII} is constructed in exactly the same way, but now $\Theta_s\PP{ys}=0$, so that the target complex is just $v\PP{y}[1]$. 
\item In the case \ref{catIII} we construct lifts (to the chain level) of elements in $\stdE(x,y)$ from those of  $(v^{-1}t+v)\stdE(xs,y)$ and  $\stdE(xs,ys)$. First, we have an inclusion $\Psi_1:v^{-1}t\stdE(xs,y)\rightarrow \stdE(x,y)$ induced from the chain level by sending $f: \PP{xs}\rightarrow \PP{y}$ to the morphism of complexes 
$$(f,0):\left(v\PP{xs}\ra{\adj{s}}\Theta_s\PP{xs}\right)\rightarrow \left(v\PP{y}\rightarrow0\right)$$ 
and then taking the induced morphism $\PP{x}\rightarrow\PP{y}$. 

Now {\it given} a morphism of chain complexes $g:v\PP{ys}\rightarrow\PP{y}$ which induces (the unique up to scalar) nonzero element in $\stdE(ys,y)$, there is an inclusion  $\Psi_2: \stdE(xs,y)\rightarrow v^{-1}\stdE(x,y)$ induced by sending $f: \PP{xs}\rightarrow \PP{y}$ to the morphism of complexes $$g\circ (f,\Theta_s(f)):\;\left(v\PP{xs}\xrightarrow{\adj{s}}\Theta_s\PP{xs}\right)\rightarrow v^{-1}\left(0\rightarrow \PP{y}\right)$$ 
and then again passing to the cone to get a morphism $\PP{x}\rightarrow\PP{y}$.

Moreover, we have a map $\Psi_3:\stdE(x,y)\rightarrow \stdE(xs,ys)$ which is induced from the chain level as follows: any chain map from $\PP{x}$ to $\PP{y}$ provides a chain map from $\Theta_s\PP{xs}$ to $\PP{y}$, or equivalently from $\PP{xs}$ to $\Theta_s\PP{y}$ by adjunction, which then in turn gives a chain map from $\PP{xs}$ to $\PP{ys}$. The first and last step is pre-composition with the obvious quotient map.

On the other hand, the differential is, by \Cref{vanilla} given by
\begin{eqnarray}
\stdE(xs,ys)&\ra{\overline{h}\circ_-,\overline{g}\circ_-}&(v^{-1}t+v)\stdE(xs,y),
\end{eqnarray} 
with homogeneous basis vectors $\overline{g}$, $\overline{h}$ from $\stdE(ys,y)$. These two basis vector can be chosen with lifts $h:v^{-1}\PP{ys}\rightarrow \PP{y}[1]$, $g:v\PP{ys}\rightarrow \PP{y}$ where $h$ projects onto the  $\PP{y}$-part in $\PP{ys}$ and $g$ is as above. 

Putting everything together, one can use the embeddings $\Psi_1,\Psi_2$ and the map $\Psi_3$ to construct inductively lifts to the chain level of basis vectors for $\stdE(x,y)$. A missing (and non-trivial) part is however the construction of $g$.  
\end{enumerate}
A more detailed treatment of a dg-lift (including the missing part) will be dealt with in article II of this series.

\begin{observation}{\rm 
We observe a remarkable behaviour in \Cref{Es} and \Cref{Shelton}: Interesting differentials only appear in the last case of \Cref{Es} and there in a very controlled way, see \eqref{vanilla}. Correspondingly, minus signs appear only in the last of Shelton's formulas and this can be tracked down using the locality of the differential. }
\end{observation}

\begin{observation}{\rm 
The embeddings $\Psi_1$, $\Psi_2$ and the map $\Psi_3$ are the crucial step to connect  \Cref{Es} with the diagrammatics we developed. For an explicit translation we refer to \Cref{appendix}.
}
\end{observation}

\begin{observation}{\rm 
The non-trivial differentials in \Cref{Es} as well as in $\richmodel{I}{J}$ only occur when square configurations appear in terms of Fukaya diagrams.  It would be interesting to give a precise interpretation of this phenomenon in terms of Floer homology of the Fukaya-Seidel category \cite{makFukayaSeidelCategories2021} associated with $\cOpo$.
}
\end{observation}

\appendix
\section{Example $\Gr(2,4)$}\label{appendix}
We discuss some of the techniques and results in the example $\Gr(2,4).$
\subsection{Open Richardson varieties and diagrams}
All open Richardson varieties are listed in \Cref{fig:allrichardsonsingrtwofour}. The corresponding  \Deodhartype\  Fukaya diagrams are listed in \Cref{fig:alldiagramsgrtwofour}; for the open Richardson variety with two strata the two strong matchings are shown on the bottom left. Cohomological \Deodhartype\  Fukaya diagrams with the corresponding basis elements are then obtained by choosing a symbol $\times$ or $\circ$ out of $\otimes$. In the basis vector $a^{\otimes,\otimes,w}$ the first and second upper index encodes the symbol on the first and second vertical strand respectively.
\begin{figure}[h]
\begin{tabular}[t]{c||c|c|c|c|c|c|c}
\diagbox[innerwidth=0.7cm]{$I$}{$J$}&\{1,2\}&\{1,3\}&\{1,4\}&\{2,3\}&\{2,4\}&\{3,4\}\\
\hline
\hline
$\{1,2\}$&$\op{pt}$& $\Gm$ &$\Gm\times \A^1$&$\Gm\times \A^1$&$\Gm\times \A^2$&$\GL_2$\\
$\{1,3\}$& $\emptyset$ & $\op{pt}$ & $\Gm$ & $\Gm$ & $(\Gm)^2$ & $\Gm\times \A^2$\\
$\{1,4\}$& $\emptyset$ & $\emptyset$ & $\op{pt}$ &$\emptyset$&$\Gm$&$\Gm\times \A^1$\\
$\{2,3\}$&$\emptyset$&$\emptyset$&$\emptyset$&$\op{pt}$&$\Gm$&$\Gm\times \A^1$\\
$\{2,4\}$&$\emptyset$&$\emptyset$&$\emptyset$&$\emptyset$&$\op{pt}$&$\Gm$\\
$\{3,4\}$&$\emptyset$&$\emptyset$&$\emptyset$&$\emptyset$&$\emptyset$&$\op{pt}$
\end{tabular}
\caption{Open Richardson varieties $\rich{I}{J}$ in $\Gr(2,4).$ The open Richardson $\rich{\{1,2\}}{\{3,4\}}\cong \GL_2$ decomposes into two strata $(\Gm)^2\times \A^2 \sqcup(\Gm)^2\times \A^1$ in the \Deodhar.}\label{fig:allrichardsonsingrtwofour}
\end{figure}
\input{alltwofourdiagrams}
\newcommand{\grpictures}[4]{
    \begin{tikzpicture}
	\def\n{4}
	\def\d{2}
	\pgfmathtruncatemacro{\dpp}{\d}
	\pgfmathtruncatemacro{\npp}{\n}

	\let\mymatrixcontent\empty
		\foreach \row in {1,...,\dpp}{
			\foreach \col in {1,...,\n}{%
				\xappto\mymatrixcontent{ \expandonce{\&}}
	}%
			\gappto\mymatrixcontent{\\}
	}
	\matrix (m) [matrix of nodes, nodes in empty cells, ampersand replacement=\&,
				matrixstyle ] (m){
	\mymatrixcontent 
	};

	\foreach \x in {3,4}
	  \draw[verticallines] (m-1-\x.north) -- (m-\d-\x.south);
	
	\foreach \ik [count=\k] in {1,2} 
		\draw[horizontallines]  (m-\d-\ik.south) -- (m-\k-\ik.center) -- (m-\k-\n.east);

	 Add matching 
	 \foreach \x/\y in {#1}
    	\node[matchingnodes] at (m-\x-\y.center) {};

	 \foreach \x/\y in {#2}
     	\node[aonenodes] at (m-\x-\y.center) {};
	 \foreach \x/\y in {#3}
        \node[gmnodes] at (m-\x-\y.center) {};
    \node [below,xshift=-0.1cm] at (current bounding box.south) {#4};
\end{tikzpicture}
}
\subsection{Computing cohomology} All (but one) open Richardson varieties are simply of the form $(\Gm)^a\times\A^b.$ For these, we can directly compute the cohomology with compact support 
$H_c^\bullet(\rich{I}{J})=\richmodel{I}{J}.$ For example, for $I=\{2,3\}$ and $J=\{3,4\}$ we obtain a two-dimensional cohomology with compact support 
\begin{align}\label{eq:middlerichardson}
    H^\bullet_c(\basechangeimage{I}{J})=\richmodel{I}{J}=\dd{s}\bigwedge\nolimits^\bullet y_1.
\end{align}

The only non-trivial case is $I=\{1,2\}$ and $J=\{3,4\}.$ In this case, $\rich{I}{J}\iso\basechangeimage{I}{J}=\bG=\GL_2.$ We compute the cohomology with compact support in two ways. First we can use that we know the cohomology of $\GL_2$ and directly compute
\begin{align}\label{eq:firstcomputationcohomologyofgl2}
    H^\bullet_c(\rich{I}{J})&=H^\bullet_c(\GL_2)=\bigwedge\nolimits^\bullet e_1,e_2
\end{align}
where $e_1=y_1+y_2$ and $e_2=y_1y_2$ are elementary symmetric polynomials.
On the other hand, we can also compute the cohomology using the model 
$$\richmodel{I}{J}=\dd e\bigwedge\nolimits^\bullet y_1,y_2\oplus \dd s\bigwedge\nolimits^\bullet  y_1,y_2$$ 
visualised in \Cref{fig:cohomologyofgltwoinrichardsontwofour}. In fact, this just comes from the long exact sequence associated to the decomposition $\bG=\bB\cup \bB s\bB.$ We obtain 
$$H^\bullet_c(\rich{I}{J})\cong H^\bullet(\richmodel{I}{J})=\langle[\dd{e}y_1y_2],[\dd{e}(y_1-y_2)],[\dd{s}y_1=\dd{s}y_2],[\dd{s}]\rangle.$$
It is convenient to use the Poincar\'e dual basis $\ddd{s}=\dd{e}$ and $\ddd{e}=\dd{s}.$
The model $\richmodel{I}{J}$ is a $dg$-algebra with unit $\ddd{e}$ and $\ddd{s}^2=0$ (these are the Nil-Coxeter relations) with a linear right action of $\bigwedge\nolimits^\bullet y_1,y_2.$ Hence, $H^\bullet(\richmodel{I}{J})$ inherits a multiplication and we can write 
\begin{align}\label{eq:secondcomputationcohomologyofgl2}
    H^\bullet_c(\rich{I}{J})\cong \bigwedge\nolimits^\bullet [\ddd{s}(y_1-y_2)],[\ddd{e}y_1].
\end{align}
Via \eqref{eq:firstcomputationcohomologyofgl2} and \eqref{eq:secondcomputationcohomologyofgl2}, $e_1$ and $e_2$ correspond to $[\ddd{e}y_1]$ and $[\ddd{s}(y_1-y_2)],$ respectively. 

\subsection{The differential}\label{fig:cohomologyofgltwoinrichardsontwofour}
Consider for example the complex $\richmodel{I}{J}$ for $I=\{1,2\}$ and $J=\{3,4\}$ with corresponding cohomological \Deodhartype\ Fukaya diagrams. In this basis, the differential sends a basis vector to a basis vector or a difference of two basis vectors (as indicated by the arrows), compare with\  \eqref{eq:differentialdiagone}, \eqref{eq:differentialdiagtwo}:
\[
\xymatrix@C=0.3pc @R=.2pc{
&\!\!\!\!\!\!\!\!\!\dd{e}y_1=a_{{12},{34},e}^{\circ,\times}\ar[rr]&&\dd{s}y_1y_2=a_{{12},{34},s}^{\circ,\circ}
&&\dd{s}y_1=a_{{12},{34},s}^{\times,\circ}\\
\dd{e}y_1y_2=a_{{12},{34},e}^{\circ,\circ}
&&&&&
&&\!\!\!\!\!\!\!\!\!\!\!\dd{s}=a_{{12},{34},s}^{\circ,\circ}
\\
&\!\!\!\!\!\!\!\!\!\dd{e}y_2=a_{{12},{34},e}^{\times,\circ}\ar[uurr]&&\dd{e}=a_{{12},{34},s}^{\times,\times}\ar[uurr]\ar[rr]^{-}&&\dd{s}y_2=a_{{12},{34},e}^{\circ,\times}\\
}
\]

\subsection{Computing compositions}
Let $H=\{1,2\},$ $I=\{2,3\}$ and $J=\{3,4\}.$ The cohomology with compact support of $\rich{H}{I}\cong \rich{I}{J}$ and $\rich{H}{J}$ is given by \eqref{eq:middlerichardson} and \eqref{eq:firstcomputationcohomologyofgl2}, respectively. With the notation $u(a)=I_2+aE_{2,2},$ we have
$$\basechangeimage{H}{I}=\basechangeimage{I}{J}=
    \setbuild{\begin{bmatrix}
    0 & 1 \\
    x & a
    \end{bmatrix}=u(a)sy_1(x)}{x\in \Gm, a\in \A^1}$$
and matrix multiplication $m$ yields an isomorphism with an open subset
$$\basechangeimage{H}{I}\times \basechangeimage{I}{J}\iso s\bB s\bB\subset \basechangeimage{H}{J}=\bG.$$ Using this, one can compute that $m_!(\dd sy_1^{\epsilon_1}\otimes \dd sy_1^{\epsilon_2})=e_1^{\epsilon_1+\epsilon_2},$ where $\epsilon_i=0,1.$

This can also be seen our model $\model.$ In fact, multiplication on $\model$ induces a map
$$\richmodel{H}{I}\otimes \richmodel{I}{J}\to \richmodel{H}{J}$$
such that $\dd sy_1^{\epsilon_1}\otimes\dd sy_1^{\epsilon_2}=\ddd ey_1^{\epsilon_1}\otimes \ddd ey_1^{\epsilon_2}\mapsto \ddd ey_1^{\epsilon_1}\ddd ey_1^{\epsilon_2}=\ddd e y_1^{\epsilon_1+\epsilon_2}.$
Passing to cohomology, the class $[\ddd e y_1^{\epsilon_1+\epsilon_2}]$ corresponds exactly to $e_1^{\epsilon_1+\epsilon_2}.$

\subsection{The complex using category $\hat{\cO}_0^\mp$}
To compute the Lie theoretic extension algebra $\stdE(\Delta^\mp)$, we first construct projective resolutions of all $\Verma{x}$, $x\in\WPl$ using \Cref{projres} and \Cref{transfunc} and then describe chain maps between these complexes which induce a basis in $\stdE(\Delta^\mp)$. We denote by $\P{x}=\P{i_1i_2\cdots i_r}$ the projective cover of $\Verma{x}$ in $\OP$ for $x=s_{i_1}\cdots s_{i_r}$ a reduced expression. The resolutions $\PP{x}$ look as follows, where we underline the cohomological degree zero:
\begin{eqnarray*}
\PP{e}=\udot{$\P{e}$},
&&\PP{21}=\left(q^2\P{e}\ra{-\op{a}}q\P{2}\ra{\op{b}}\udot{$\P{21}$}\right)\\ 
\PP{2}=\left(q\P{e}\ra{\op{a}}\udot{$\P{2}$}\right)&&
\PP{23}=\left(q^2\P{e}\ra{-\op{a}}q\P{2}\ra{\op{b}'}\udot{$\P{23}$}\right)
\end{eqnarray*}
\begin{equation*}
\PP{213}=\left(q^3\P{e}\ra{\op{a}}q^2\P{2}\ra{\left(\!\begin{smallmatrix}-\op{b}'& \op{b}\end{smallmatrix}\!\right)^t}q{\P{23}\oplus q\P{21}}\ra{\left(\!\begin{smallmatrix}\op{b}& \op{c}\end{smallmatrix}\!\right)}\udot{$\P{231}$}\right)
\end{equation*}
The inductively defined resolution $\PP{2312}$  has underlying minimal resolution 
\begin{equation*}
q^4\P{e}\ra{-\op{a}}q^3\P{2}\ra{\left(\!\begin{smallmatrix}\op{b}'& -\op{b}\end{smallmatrix}\!\right)^t} q^2{\P{23}}\oplus q^2{\P{21}}\ra{\op{A}}q{\P{2}}\oplus q{\P{231}} \ra{\left(\!\begin{smallmatrix}d,-e\end{smallmatrix}\!\right)}\udot{$\P{2312}$}
\end{equation*}
\begin{remark}
If we abbreviate by $\adj{i}=\adj{s_i}$ and  $\adjov{i}=\adjov{s_i}$ the adjunction morphisms from \Cref{transfunc}, then the maps are given by $\op{a}=\adj{2}$, $\op{b}=\adj{1}$, $\op{b}'=\adj{3}$, $\op{c}=\Theta_{s_1}(\op{b}')$, $\op{A}=\left(\!\begin{smallmatrix}\op{a}&\op{a}\\ \op{b}&\op{c}\end{smallmatrix}\!\right)$ using that $\theta_{s_2}\theta_{s_1}\theta_{s_2}=\theta_{s_2}=\theta_{s_2}\theta_{s_3}\theta_{s_2}$, $e=\overline{\op{a}}$. Then there is a unique $d$ turning this into a complex. Some readers might prefer to work (equivalently) with the category of modules over the corresponding extended Khovanov arc algebra $K_2^2$ from \cite{BS3}. Then the maps are chosen such that all letters represent a basis element in $K_2^2$ which is uniquely determined by the source and target.
\end{remark}
\begin{example}
The maps $(-\op{id},\udot{0})$ and $(0,\udot{$\overline{\op{a}}$}=(0,\udot{$\adjov{2}$})$ define morphisms of complexes $\PP{2}\rightarrow \PP{e}$  of bidegrees $(1,-1)$ and $(0,1)$. They induce a basis of $E(s_2,e)$  corresponding to the \Deodhartype\ diagrams $a_{{24},{34}}^{\times}$, $a_{{23},{34}}^{\circ}$ respectively.
\end{example}
\begin{example}
The maps $(\op{id},-\op{id},\left(0\; -\op{id}\right),0)$ and 
$(0,\adjov{2},\left(\adjov{2}\;0\right),\adjov{2})$ define  morphisms $\PP{231}\rightarrow \PP{21}$ of complexes and induce a basis of $E(s_2s_3s_1,s_2)$ 
corresponding to the \Deodhartype\ diagrams $a_{{13},{23}}^{\times}$ and  $a_{{13},{23}}^{\circ}$ respectively.
\end{example}
\begin{example}
The maps $(\op{id},0,0,\udot{0})$ and $(0,\overline{\op{a}}=\adjov{2},0,0,\udot{0})$ define morphisms of complexes $\PP{231}\rightarrow \PP{e}$.  They induce a basis of $E(s_2s_1s_3,e)$ in bidegrees $(3,-3)$ and $(2,-1)$ corresponding to the \Deodhartype\ diagrams $a_{{13},{34}}^{\circ}$, $a_{{13},{34}}^{\times}$ respectively.
\end{example}
\begin{example}
The maps $(-\op{id},-\op{id},0,\udot{0})$, 
$(0,-\adjov{2},\left(0\;-{\adjov{1}}\right),\udot{0})$ together with 
$(0,\adjov{2},\left(-\adjov{3}\; 0\right),\udot{0})$, 
$(0,0,0,\udot{$\adjov{1}\adjov{3}$})$ define morphisms of complexes  $\PP{213}\rightarrow \PP{2}$ which induce a basis of $E(s_2s_1s_3,s_2)$ in bidegrees $(2,-2)$, $(1,0)$, $(1,0)$, $(0,2)$
corresponding to the \Deodhartype\ diagrams $a_{{13},{24}}^{\times,\times}$, $a_{{13},{34}}^{\times,\circ}$, $a_{{13},{34}}^{\circ,\times}$, $a_{{13},{34}}^{\circ,\circ}$ respectively.
\end{example}
From \Cref{Es} and \Cref{vanilla} we obtain that $E(s_2s_1s_3s_2,e)$ can be computed as the cohomology of the complex  $$E(s_2s_1s_3,s_2)\ra{(a_{{13},{34}}^{\circ}\circ_-,a_{{13},{34}\circ_-}^{\times})} (q^{-1}t+q)E(s_2s_1s_3,e).$$ We recover the complex from \Cref{fig:cohomologyofgltwoinrichardsontwofour}.

\subsection{Quiver description}
Altogether, we obtain a description of the extension algebra $\stdE(\Delta^\mp)=\E(\Delta)$ as a quotient of a path algebra of a quiver modulo relations. 
The arrows from $\{i_1,i_2\}$ to $\{i_1',i_2'\}$ correspond to the elements  $a_{i_1i_2\;{j_1j_2}}^\circ$ (for the blue arrows) and $a_{i_1i_2\;{j_1 j_2}}^\times$ (for the black arrow) and the relations are as follows:
\begin{center}
\begin{minipage}[t][5cm][c]{6cm}
\begin{tikzcd}
	{\{1,2\}} & {} & {}\\
	{\{1,3\}} & {\{1,4\}} &&\\
	{\{2,3\}}& {\{2,4\}} & {\{3,4\}}\\
	\arrow[shift left=1, color={rgb,255:red,92;green,92;blue,214}, from=2-1, to=3-1]
	\arrow[shift left=1, color={rgb,255:red,92;green,92;blue,214}, from=2-2, to=3-2]
	\arrow[shift left=1, color={rgb,255:red,92;green,92;blue,214}, from=2-1, to=2-2]
	\arrow[shift left=1, color={rgb,255:red,92;green,92;blue,214}, from=3-1, to=3-2]
	\arrow[shift left=1, color={rgb,255:red,92;green,92;blue,214}, from=1-1, to=2-1]
	\arrow[shift left=1, color={rgb,255:red,92;green,92;blue,214}, from=3-2, to=3-3]
	\arrow[shift right=1, from=2-1, to=2-2]
	\arrow[shift right=1, from=1-1, to=2-1]
	\arrow[shift right=1, from=3-1, to=3-2]
	\arrow[shift right=1, from=3-2, to=3-3]
	\arrow[shift right=1, from=2-1, to=3-1]
	\arrow[shift left=1, color={rgb,255:red,92;green,214;blue,92}, curve={height=-18pt}, from=1-1, to=3-3]
	\arrow[shift right=1, from=2-2, to=3-2]
	\arrow[shift right=1, color={rgb,255:red,214;green,92;blue,92}, curve={height=-18pt}, from=1-1, to=3-3]
\end{tikzcd}
\end{minipage}
\begin{minipage}[t][4cm][c]{5cm}
\begin{eqnarray*}
a_{13\;{14}}^\circ a_{13\;{14}}^\circ&=&a_{23\;{14}}^\circ a_{13\;{23}}^\circ\\
a_{13\;{14}}^\times a_{13\;{14}}^\circ&=&a_{23\;{14}}^\circ a_{13\;{23}}^\times\\
a_{13\;{14}}^\circ a_{13\;{14}}^\times&=&a_{23\;{14}}^\times a_{13\;{23}}^\circ\\
a_{13\;{14}}^\times a_{13\;{14}}^\times&=&-a_{23\;{14}}^\times a_{13\;{23}}^\times\\
a_{24\;{34}}^\circ a_{14\;{24}}^\times&=&a_{24\;{34}}^\times a_{14\;{24}}^\circ\\
a_{24\;{34}}^\circ a_{23\;{24}}^\times&=&a_{24\;{34}}^\times a_{23\;{24}}^\circ\\
a_{24\;{34}}^\circ a_{14\;{24}}^\circ\;=&0&=\;a_{24\;{34}}^\circ a_{23\;{24}}^\circ\\
a_{13\;{23}}^\circ a_{12\;{13}}^\circ\;=&0&=\;a_{13\;{14}}^\circ a_{12\;{13}}^\circ
\end{eqnarray*}
\end{minipage}
\end{center}
We have realised the standard extension algebra via generators and relations. 


\bibliographystyle{amsalpha} 
\bibliography{main}

\end{document}